\documentclass[english]{article}

\usepackage{geometry}
\usepackage[T1]{fontenc}
\usepackage[latin9]{inputenc}
\usepackage{bm}
\usepackage{amsmath,mathtools}
\usepackage{amssymb}
\usepackage[unicode=true,
 bookmarks=false,
 breaklinks=false,pdfborder={0 0 1},colorlinks=false]
 {hyperref}
\hypersetup{
 colorlinks,citecolor=blue,filecolor=blue,linkcolor=blue,urlcolor=blue}
\makeatletter
\usepackage{amsthm}
\usepackage{cite}  
\usepackage{comment}
\usepackage{natbib}
\usepackage{booktabs}
\usepackage{graphicx}
\usepackage[linesnumbered,ruled,vlined]{algorithm2e}
\usepackage{algorithmic}
\usepackage{float}
\usepackage{float}
\usepackage{multirow}
\usepackage{dsfont}
\usepackage{tcolorbox}

\usepackage{color}
\definecolor{yxc}{RGB}{255,0,0}
\definecolor{yjc}{RGB}{125,0,0}
\definecolor{ytw}{RGB}{255,69,0}
\definecolor{gen}{RGB}{0,0,200}
\definecolor{wei}{RGB}{0,125,125}

\allowdisplaybreaks

\DeclareMathOperator{\ind}{\mathds{1}}  

\newcommand{\second}{\prime\prime}
\newcommand{\dx}{\mathrm{d} x}
\newcommand{\lt}{\left}
\newcommand{\rt}{\right}

\newcommand{\vstar}{v^\star}
\newcommand{\third}{\prime\prime\prime}
\newcommand{\Xbayes}{\widehat{X}^{\textrm{bayes}}}

\newcommand{\ltwo}[1]{\|#1\|_2}

\newcommand{\real}{\ensuremath{\mathbb{R}}}

\newcommand{\inprod}[2]{\ensuremath{\langle #1 , \, #2 \rangle}}

\newcommand{\mprob}{\ensuremath{\mathbb{P}}}

\newcommand{\defn}{\coloneqq}

\newcommand{\Exs}{\ensuremath{\mathbb{E}}}

\long\def\comment#1{}

\newcommand{\HACKPROOF}{\begin{proof}}
\newcommand{\HACKENDPROOF}{\end{proof}}

\newlength{\widebarargwidth}
\newlength{\widebarargheight}
\newlength{\widebarargdepth}

\makeatletter
\long\def\@makecaption#1#2{
        \vskip 0.8ex
        \setbox\@tempboxa\hbox{\small {\bf #1:} #2}
        \parindent 1.5em  
        \dimen0=\hsize
        \advance\dimen0 by -3em
        \ifdim \wd\@tempboxa >\dimen0
                \hbox to \hsize{
                        \parindent 0em
                        \hfil 
                        \parbox{\dimen0}{\def\baselinestretch{0.96}\small
                                {\bf #1.} #2
                                } 
                        \hfil}
        \else \hbox to \hsize{\hfil \box\@tempboxa \hfil}
        \fi
        }
\makeatother

\theoremstyle{plain}

\newtheorem{theo}{Theorem}[section]

\newtheorem{lem}{Lemma}[section]
\newtheorem{prop}{Proposition}[section]
\newtheorem{cor}{Corollary}[section]

\theoremstyle{definition} 

\newtheorem{nota}{Notation}[section]
\newtheorem{de}{Definition}[section]
\newtheorem{exa}{Example}[section]
\newtheorem{as}{Assumption}[section]
\newtheorem{alg}{Algorithm}[section]

\newcommand{\btheo}{\begin{theo}}
\newcommand{\bde}{\begin{de}}
\newcommand{\ble}{\begin{lem}}
\newcommand{\bpr}{\begin{prop}}
\newcommand{\bno}{\begin{nota}}
\newcommand{\bex}{\begin{exa}}
\newcommand{\bcor}{\begin{cor}}
\newcommand{\spro}{\begin{proof}}
\newcommand{\bas}{\begin{as}}
\newcommand{\balg}{\begin{alg}}

\newcommand{\etheo}{\end{theo}}
\newcommand{\ede}{\end{de}}
\newcommand{\ele}{\end{lem}}
\newcommand{\epr}{\end{prop}}
\newcommand{\eno}{\end{nota}}
\newcommand{\eex}{\end{exa}}
\newcommand{\ecor}{\end{cor}}
\newcommand{\fpro}{\end{proof}}
\newcommand{\eas}{\end{as}}
\newcommand{\ealg}{\end{alg}}

\theoremstyle{plain}

\newtheorem{theos}{Theorem}
\newtheorem{props}{Proposition}
\newtheorem{lems}{Lemma}
\newtheorem{cors}{Corollary}

\theoremstyle{definition}
\newtheorem{exas}{Example}
\newtheorem{algs}{Algorithm}
\newtheorem{asss}{Assumption}
\newtheorem{defns}{Definition}

\newcommand{\btheos}{\begin{theos}}
\newcommand{\etheos}{\end{theos}}
\newcommand{\bprops}{\begin{props}}
\newcommand{\eprops}{\end{props}}
\newcommand{\bdes}{\begin{defns}}
\newcommand{\edes}{\end{defns}}
\newcommand{\blems}{\begin{lems}}
\newcommand{\elems}{\end{lems}}
\newcommand{\bcors}{\begin{cors}}
\newcommand{\ecors}{\end{cors}}
\newcommand{\bexs}{\begin{exas}}
\newcommand{\eexs}{\end{exas}}
\newcommand{\balgs}{\begin{algs}}
\newcommand{\ealgs}{\end{algs}}
\newcommand{\bass}{\begin{asss}}
\newcommand{\eass}{\end{asss}}

\theoremstyle{plain} \newtheorem{lemma}{\textbf{Lemma}}\newtheorem{theorem}{\textbf{Theorem}}\newtheorem{assumption}{\textbf{Assumption}}

\theoremstyle{remark}\newtheorem{remark}{\textbf{Remark}}

\title{Approximate message passing from random initialization \\ with applications to $\mathbb{Z}_{2}$ synchronization}

\date{\today}
\makeatother

\begin{document}

\author{Gen Li ~~~~~ Wei Fan ~~~~~ Yuting Wei \\[.2in]
	Department of Statistics and Data Science, the Wharton School \\
	University of Pennsylvania, Philadelphia, PA
}
\maketitle

\begin{abstract}

This paper is concerned with the problem of reconstructing an unknown rank-one matrix with prior structural information from noisy observations.
While computing the Bayes-optimal estimator seems intractable in general due to its nonconvex nature, Approximate Message Passing (AMP)  emerges as an efficient first-order method to approximate the Bayes-optimal estimator. However, the theoretical underpinnings of AMP remain largely unavailable when it starts from random initialization, a scheme of critical practical utility. Focusing on a prototypical model called $\mathbb{Z}_{2}$ synchronization, we characterize the finite-sample dynamics of AMP from random initialization, uncovering its rapid global convergence. Our theory provides the first non-asymptotic characterization of AMP in this model without requiring either an informative initialization (e.g., spectral initialization) or sample splitting. 

\end{abstract}


\medskip


\setcounter{tocdepth}{2}
\tableofcontents

\section{Introduction}

The problem of estimating an unknown structured signal $\vstar\in \mathbb{R}^n$, when given  access to noisy observations 
\begin{align}
\label{eqn:spike-Wigner}
	M = \lambda \vstar v^{\star\top} + W \in \mathbb{R}^{n\times n} 
	\qquad \text{with }\lambda > 0
\end{align}
is of fundamental interest and has been investigated in a diverse array of contexts \citep{singer2011angular,abbe2020entrywise,johnstone2001distribution,keshavan2009matrix,zhong2018near,candes2010matrix,chi2019nonconvex}. 
This model is commonly referred to as a deformed Wigner model or spiked Wigner model 
when the entries of the noise matrix $W=[W_{ij}]_{1\leq i,j\leq n}$ are independently drawn from Gaussian distributions --- more precisely, $W_{ii} \overset{\mathrm{i.i.d.}}{\sim} \mathcal{N}(0,\frac{2}{n})$ and $W_{ji} = W_{ij} \overset{\mathrm{i.i.d.}}{\sim} \mathcal{N}(0,\frac{1}{n})$ for $i \neq j$ ---   
which serves as a prototypical model towards understanding the feasibility and fundamental limits of low-rank matrix estimation.  


The spectral properties of the observed matrix $M$ has been extensively studied  (see, e.g.~\cite{peche2006largest,baik2005phase,feral2007largest,capitaine2009largest,cheng2021tackling}), motivating the design of spectral methods when there is no structural information associated with $\vstar$ \citep{singer2011angular,keshavan2009matrix,chen2021spectral,cai2018rate,yan2021inference}. 
In practice, there is no shortage of applications where additional structural information about $\vstar$ is available {\em a priori}, 
 examples including finite-group structure \citep{perry2018message}, cone constraints \citep{deshpande2014cone,lesieur2017constrained}, sparsity \citep{johnstone2009consistency,berthet2013optimal}, among others. 
The presence of prior structure further exacerbates the nonconvexity issue when computing the maximum likelihood estimate or Bayes-optimal estimate, 
thereby presenting pressing needs for designing algorithms that can be executed efficiently.  
Remarkably, the approximate message passing (AMP) algorithm emerges as an efficient nonconvex paradigm that rises to the aforementioned challenge \citep{feng2022unifying,donoho2009message}. 
Originally proposed in the context of compressed sensing, 
AMP has served as not only a family of first-order iterative algorithms that enjoy rapid convergence \citep{donoho2010message,bayati2011dynamics,fan2022approximate,rangan2011generalized,celentano2022fundamental}, but also a powerful machinery that assists in determining the performance of other statistical procedures in high-dimensional asymptotics \citep{bayati2011lasso,li2021minimum,sur2019likelihood,zhang2022precise,donoho2016high,ma2018optimization,lelarge2019fundamental,javanmard2016phase,bu2020algorithmic,sur2019modern}.


Nevertheless, while the existing suite of AMP theory covers a wealth of applications, it remains inadequate in at least two aspects. 
To begin with, a dominant fraction of existing AMP theory is asymptotic in nature, in the sense that it predicts the AMP dynamics in the large-$n$ limit for any {\em fixed} iteration $t$. For this reason, prior AMP theory falls short of describing how AMP behaves after a growing number of iterations, 
which stands in contrast to other optimization-based procedures that often come with non-asymptotic analysis accommodating a large number of iterations \citep{chen2015fast,chi2019nonconvex,keshavan2009matrix,ma2020implicit}. 
Another issue stems from the requirement of an informative initialization, that is, existing AMP theory for low-rank estimation often requires starting from a point that already enjoys non-vanishing correlation with the true signal \citep{montanari2021estimation,celentano2021local,zhong2021approximate}. 
While an informative initial estimate like spectral initialization is sometimes plausible and analyzable, 
this requirement presents a hurdle to understanding the effect of other widely adopted alternatives like random initialization. 
As shall be made clear shortly, tackling this issue might also necessitate a new non-asymptotic framework for AMP, due to the difficulty of tracking the AMP dynamics when the iterates exhibit only extremely weak correlation with the truth.


Inspired by the aforementioned issues, 
there has been growing interest in understanding the finite-sample performance of AMP. 
The first work of this kind was \citet{rush2018finite} (see also its follow-up work \citet{cademartori2023non}), 
which studied AMP for sparse regression and permitted the total number of iterations to be as large as $o\big( \frac{\log n}{\log \log n}\big)$. 
A recent work \citet{li2022non} developed a non-asymptotic framework for the spiked Wigner models, 
which characterized the AMP behavior for up to $O\big(\frac{n}{\mathsf{poly}(\log n)}\big)$ iterations. 
Although the theory therein is well-suited to spectrally initialized AMP, 
it remains in mystery whether randomly initialized  AMP would be able  to achieve the same performance as AMP with informative initialization.

%
%

\subsection{This paper: randomly initialized AMP for $\mathbb{Z}_2$ synchronization}

In this work, we attempt to address the above challenges by studying a concrete model called $\mathbb{Z}_2$ synchronization. 
To be precise, $\mathbb{Z}_2$ synchronization is a special case of the spiked Wigner model  
when the ground truth is known to have a discrete structure obeying $\vstar \in \{\pm\frac{1}{\sqrt{n}}\}^{n}$. 
Here and throughout, we impose a prior distribution on $\vstar=[\vstar_i]_{1\leq i\leq n}$ such that 
$$v^\star_{i} \stackrel{\mathrm{i.i.d.}}{\sim} \mathsf{Unif}\Big(\pm \frac{1}{\sqrt{n}}\Big), \qquad 1\leq i\leq n.$$
The goal is to reconstruct $\vstar$ on the basis of the measurements $M$ (see \eqref{eqn:spike-Wigner}). 
This problem can be viewed as a basic example of a more general problem --- synchronization over compact groups \citep{abbe2020entrywise,singer2011angular,chen2018projected,perry2018message,zhong2018near,gao2022sdp}, and has an intimate connection to stochastic block models \citep{deshpande2017asymptotic,lelarge2019fundamental}.

\paragraph{The AMP algorithm.}
Due to the combinatorial nature of the underlying optimization problem, it is in general intractable to calculate the Bayes-optimal solution directly.  
This motivates the search for computationally feasible alternatives, for which AMP emerges as a natural and successful option \citep{fan2021tap,deshpande2017asymptotic,celentano2021local,li2022non}. 
More concretely, given the initialization points $x_0, x_{1}\in \real^{n}$, AMP tailored to $\mathbb{Z}_2$ synchronization adopts the following update rule: 
\begin{align} 
\label{eqn:AMP-updates}
	x_{t+1} = M\eta_t(x_{t}) - \langle\eta_t^{\prime}(x_{t})\rangle\eta_{t-1}(x_{t-1}), 
	\qquad t\geq 1,
\end{align}
where we denote $\langle x\rangle \defn \frac{1}{n}\sum_{i = 1}^n x_i$ for any vector $x=[x_i]_{1\leq i\leq n}$, and the denoising function is given by\footnote{Note that for ease of analysis, we adopt a slightly different scaling from that of \cite{deshpande2017asymptotic}, but they are equivalent up to global scaling.}
\begin{equation}
\label{defi:eta}
\begin{aligned} 
\eta_t(x) = \gamma_t\tanh(\pi_tx),&\qquad\text{for } t \ge 1\\
\text{with }~
\pi_t \defn \sqrt{\max\big\{n(\|x_t\|_2^2 - 1), 1\big\}}& \quad\text{and}\quad \gamma_t \defn \|\tanh(\pi_tx_t)\|_2^{-1}.
\end{aligned}
\end{equation}
Here, it is understood that $\eta_t(\cdot)$, $\eta_t^{\prime}(\cdot)$ and $\tanh(\cdot)$ are applied entrywise if the input argument is a vector.

Thus far, there are two analysis strategies that accommodate a growing number of iterations  in the most challenging regime (i.e.~$\lambda > 1$). 
One attempt was made by \cite{celentano2021local}, which proposed a three-stage hybrid algorithm that runs spectrally initialized AMP followed by natural gradient descent (NGD). 
It was conjectured therein that the third stage (i.e.~NGD) is unnecessary. 
Recently, \cite{li2022non} put forward another strategy to address this conjecture,  showing that spectrally-initialized AMP alone is sufficient without the need of a third refinement stage. 
Despite the nonconvex nature of the underlying optimization problem, AMP with spectral initialization performs nearly the same as the Bayes-optimal estimate.



\begin{figure}[t]
\begin{center}
\begin{tabular}{cc}
\includegraphics[width=0.48\textwidth]{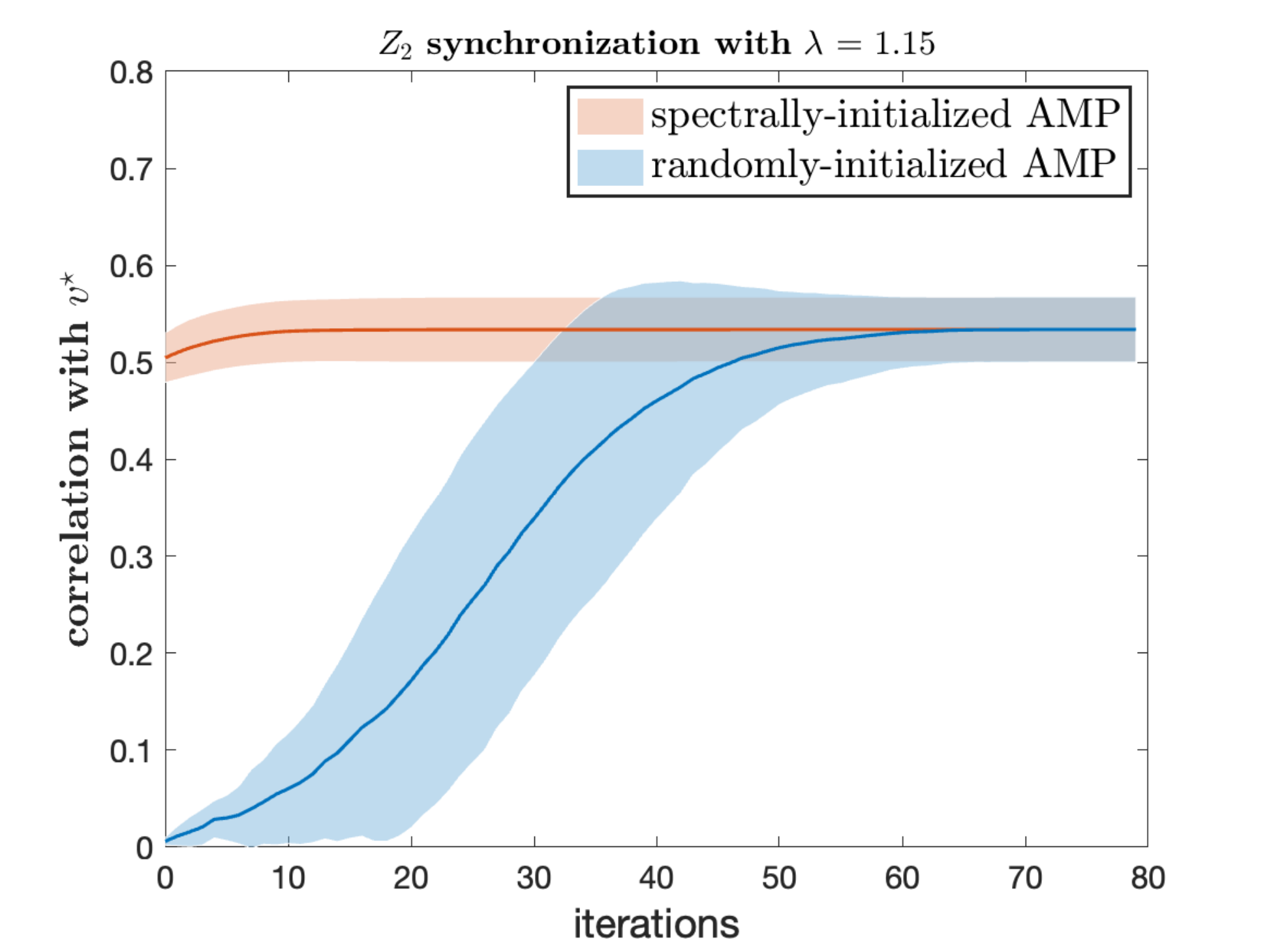} & \includegraphics[width=0.48\textwidth]{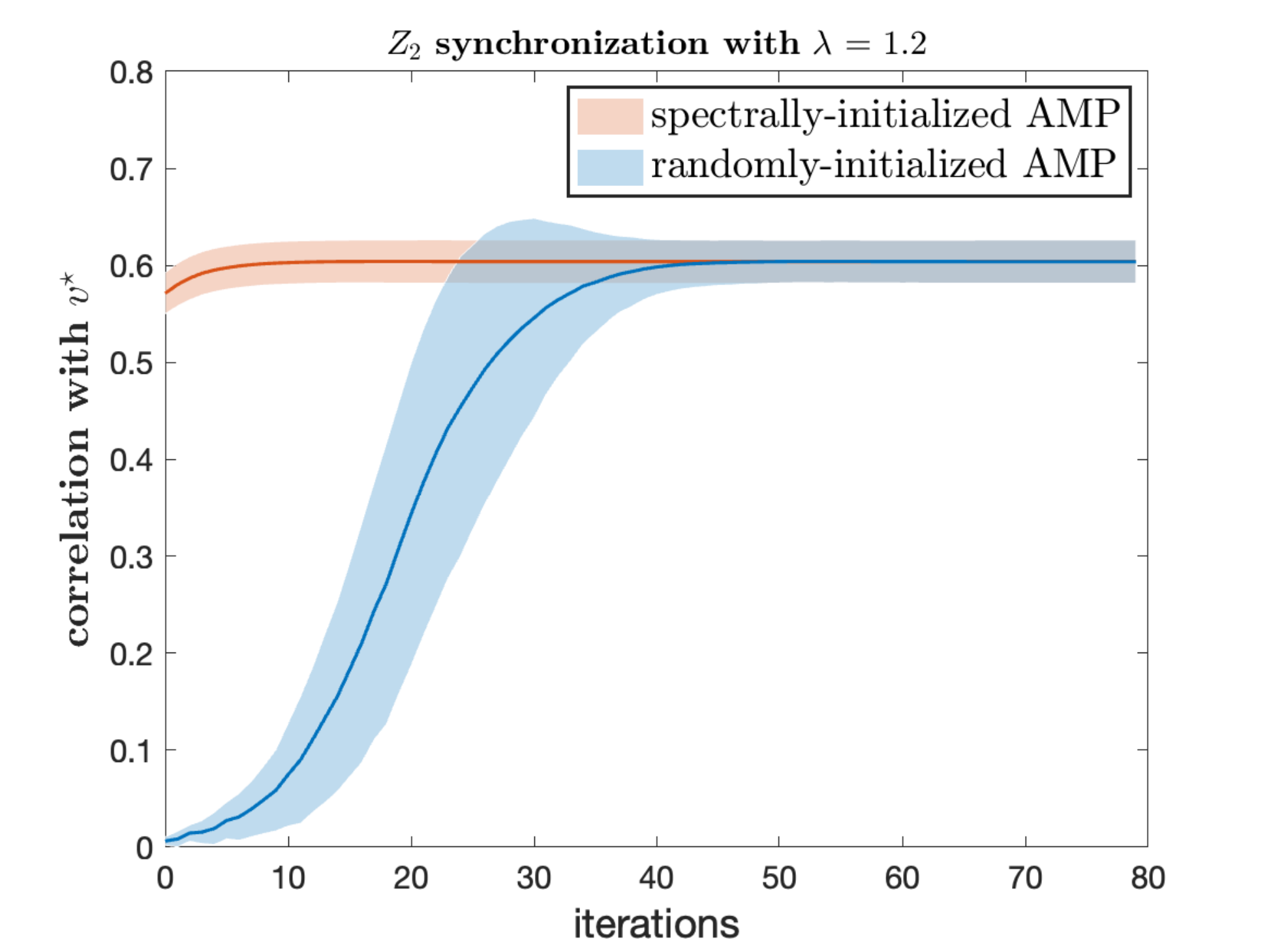} \\
(a) AMP with $\lambda = 1.15$ & (b) AMP with $\lambda = 1.2$ \\
\end{tabular}
\end{center}
\caption{
	The correlation of $\eta_{t}(x_t)$ and $\vstar$ (i.e.~$\frac{|\inprod{\eta_t(x_t)}{\vstar}|}{\ltwo{\eta_t(x_t)}}$) vs.~iteration count $t$ 
	for AMP with both random and spectral initialization. 
	Here, $n = 10000$ and $v^\star_{i} \stackrel{\mathrm{i.i.d.}}{\sim} \mathsf{Unif}(\pm \frac{1}{\sqrt{n}})$ ($1\leq i\leq n$). 
	We generate 20 independent copies of $M$ according to \eqref{eqn:spike-Wigner} and report the averaged results,  with the width of the shaded region reflecting (twice) the standard deviation. 
Figure (a) and (b) correspond to $\lambda  = 1.15$ and $\lambda = 1.2$,  respectively. 
%
}
\label{fig:amp-plot}
\end{figure}

\paragraph{Random initialization.}
%
%
As alluded to previously, all existing AMP theories for this problem \citep{montanari2021estimation,mondelli2021pca,celentano2021local,li2022non} require informative initialization obtained by, for example,  spectral methods. This leads to the following natural questions: 
\begin{center}
	\emph{Is spectral initialization necessary for the success of AMP? \\ Can we start with a simpler initialization scheme but still work equally well as spectral initialization?} 
\end{center}
\noindent 
One initialization strategy that enjoys widespread adoption is to initialize AMP randomly; 
for instance, 
\begin{align}
\label{eqn:simulation-setting}
	x_1 \sim \mathcal{N}\bigg(0, \frac{1}{n}I_n\bigg) \qquad \text{ (independent of } M)
	\qquad \text{and} \qquad \eta_0(x_0)=0.
\end{align}
In order to answer the above questions, let us first conduct a series of numerical experiments using \eqref{eqn:simulation-setting}, 
as reported in Figure~\ref{fig:amp-plot}.  
Encouragingly,  AMP with random initialization seems to work surprisingly well: 
it only takes several tens of iterations to achieve nearly the same performance as spectrally initialized AMP (note that spectral initialization also consists of several tens of power iterations). 
Such encouraging numerical results motivate us to pursue in-depth theoretical understanding about the effect of random initialization upon AMP convergence, 
which was poorly understood in the literature.


\subsection{Main contributions and technical challenges}
In the present paper, we provide a non-asymptotic analysis that allows one to predict how AMP evolves over time from random initialization. 
Our theory is able to track the correlation of the AMP iterates and the truth $\vstar$.  
In particular, we demonstrate in Theorem~\ref{thm:main} that the signal component in the AMP iterates increases exponentially fast at the initial stage, 
taking no more than $O\big(\frac{\log n}{\lambda - 1}\big)$ iterations to grow from $\widetilde{O}(\frac{1}{\sqrt{n}})$ to $O(\sqrt{\lambda^2 - 1})$ (the latter of which coincides with the correlation of spectral initialization and the truth). 
Furthermore, once the signal component surpasses $\frac{1}{2}\sqrt{\lambda^2 - 1}$ in magnitude, 
the finite-sample AMP dynamics are very well predicted by the state evolution recursion derived in the asymptotics. 
To the best of our knowledge, 
our paper delivers the first result to characterize the performance of AMP from random initialization, 
justifying and advocating the use of random initialization in practice.


Built upon the analysis recipe recently developed by \cite{li2022non}, 
the development of our theory requires novel ideas beyond this framework in order to track AMP from random initialization. 
Before continuing, we take a moment to single out the key challenges that need to be overcome. 
\begin{itemize}
\item 
	Prior theory based on state-evolution analysis falls short of providing ``fine-grained'' understanding about the iterates when they  have vanishing correlation with the truth. More precisely, past theory fails to measure the progress of AMP during the initial stage when its signal component is $o(1)$ (in fact, as small as $\widetilde{O}(\frac{1}{\sqrt{n}})$ when initialized), but instead treats the signal strength as 0 in the large-$n$ limit.

\item Another technical challenge results from the complicated statistical dependency across iterations, which is particularly difficult to cope with when the algorithm starts with random initialization and when the number of iterations grows with the dimension $n$. While prior literature tackles this issue for other nonconvex optimization methods by resorting to either delicate leave-one-out decoupling arguments (see, e.g.~\citet{chen2019gradient}) or global landscape analysis (see, e.g.~\cite{ge2017no}), these approaches remain unavailable when analyzing AMP. 
\end{itemize}


%
%
%

\subsection{Notation} 

Throughout this paper, we use $\varphi(\cdot)$ (resp.~$\varphi_n(\cdot)$) to denote the probability density function (p.d.f.) of a standard Gaussian random variable (resp.~a Gaussian random vector $\mathcal{N}(0,I_n)$). 
For any matrix $M$, we let $\|M\|$ and $\|M\|_{\mathsf{F}}$ denote the spectral norm and the Frobenius norm of $M$, respectively. 
If not noted otherwise, for any vector $x\in [x_i]_{1\leq i\leq n} = \real^n$, we denote by $|x|_{(i)}$ (resp. $x_{(i)}$) the $i$-th largest absolute value (resp. value) of $x$ in magnitude. 
We write $\mathcal{S}^{d-1}=\{x\in \real^d \mid \|x\|_2=1\}$ as the unit sphere in $\real^d$. 
Moreover, for any two vectors $x,y\in \real^n$, we write $x \circ y $ for their Kronecker product, namely, $x \circ y = (x_1y_1,\ldots, x_ny_n)^{\top} \in \real^{n}.$
When a function is applied to a vector, it should be understood as being applied in a component-wise fashion; 
for instance,  for any vector $x=[x_i]_{1\leq i\leq n}$, we let $x + 1\coloneqq [x_i + 1]_{1\leq i\leq n}$. 

In addition, given two functions $f(n)$ and $g(n)$, we write $f(n)\lesssim g(n)$ or $f(n)=O(g(n))$ to indicate that $|f(n)|\leq c_1 g(n)$ for some universal constant $c_1>0$ independent of $n$, and similarly, $f(n)\gtrsim g(n)$ means that $f(n)\geq c_2 |g(n)|$ for some universal constant $c_2>0$. 
We write $f(n)= \widetilde{O}(g(n))$ if $f(n)=O(g(n))$ up to logarithm factors. 
We also adopt the notation $f(n)\asymp g(n)$ to indicate that both $f(n)\lesssim g(n)$ and $f(n)\gtrsim g(n)$ hold simultaneously. 
Moreover, when we write $f(n) \ll g(n)$ or $f(n)=o(g(n))$, it means $f(n)/g(n)\to 0$ as $n\to \infty$;  we also write $f(n) \gg g(n)$, if $g(n)/f(n)\to 0$ as $n\rightarrow \infty$. 
We use $c, C$ to denote universal constants that do not depend on $n$; note that constants may change from line to line. 
%

\section{Main results}
\label{sec:main-results}

In this section, 
we provide precise statements of our main theoretical guarantees for randomly initialized AMP. 
For notational convenience, let us introduce 
\begin{equation}
	\alpha_{t+1} \defn \lambda v^{\star\top} \eta_{t}(x_t),
	\label{eq:defn-alpha-t-signal}
\end{equation}
which captures the projection of the $t$-th iterate (after denoising) onto the direction of the truth $\vstar$. 
In some sense, this quantity captures the size of the signal component carried by the $t$-th iterate. 
With this notation in place, we single out a key threshold as follows:
\begin{align}
	\varsigma := \min\bigg\{ t: \alpha_t \geq \frac{1}{2}\sqrt{\lambda^2 - 1} \bigg\}, 
	\label{eq:defn-varsigma}
\end{align}
which reflects the time taken for the AMP iterate to carry a significant signal component (note that a random initial guess obeys $|v^{\star\top} x_1|\lesssim \widetilde{O}\big(\frac{1}{\sqrt{n}}\big)$, meaning that the initial signal component is exceedingly small).  
Additionally, we define the state-evolution recursion starting from the $\varsigma$-th iteration as follows:
\begin{align} 
\label{eq:SE}
\alpha_{\varsigma}^{\star} = |\alpha_{\varsigma}|
	\qquad \text{and} \qquad
\alpha_{t+1}^{\star} &= \lambda \lt[\int \tanh\lt(\alpha_t^{\star}\lt(\alpha_t^{\star} + x\rt)\rt)\varphi(\dx)\rt]^{1/2}
	\quad (\forall t\geq \varsigma).
\end{align}
Equipped with the above definitions, our main results are summarized in the following theorem.

\begin{theorem} 
\label{thm:main}
Consider the $\mathbb{Z}_2$ synchronization problem with $n^{-1/9}\log n \lesssim \lambda - 1 \le 0.2$. 
Suppose we run  AMP (cf.~\eqref{eqn:AMP-updates} and \eqref{defi:eta}) with random initialization \eqref{eqn:simulation-setting}.  
Consider any $t$ obeying $1\le t\le \frac{c n(\lambda-1)^5}{\log^2 n}$, where $c > 0$ is some universal constant. 
Then with probability at least $1 - O(n^{-10})$, the following results hold:  
\begin{itemize}
	\item \textbf{(Decomposition and error bound)} The AMP iterates admit the decomposition 
\begin{subequations}
\begin{align}
\label{eqn:xt-decomposition}
	x_t = \alpha_{t} v^{\star} + \sum_{k = 1}^{t-1} \beta_{t-1}^k\phi_k + \xi_{t-1},
\end{align}
	where $\alpha_t$ is defined in \eqref{eq:defn-alpha-t-signal}, 
	the $\phi_k$'s are i.i.d.~Gaussian vectors obeying $\phi_k \overset{\text{i.i.d.}}{\sim} \mathcal{N}(0, \frac{1}{n}I_n)$, 
	and 
\begin{align}
	\ltwo{\eta_t(x_t)} &= \|\beta_t\|_2 \defn  \lt\|(\beta_t^1,\beta_t^2,\ldots,\beta_t^t)\rt\|_2 = 1, \\
	\|\xi_t\|_2  &\lesssim 
	\sqrt{\frac{t\log n}{n(\lambda-1)^2}} + \sqrt{\frac{\log^4 n}{n(\lambda-1)^3}}; \label{eqn:xi-norm-fin}
\end{align}
\end{subequations}

	\item \textbf{(Crossing time)} The threshold $\varsigma$ defined in \eqref{eq:defn-varsigma} satisfies
\begin{equation}
	\varsigma = O\bigg(\frac{\log n}{\lambda - 1}\bigg); 
\end{equation}
	\item \textbf{(Non-asymptotic state evolution)} For any $t$ obeying $\varsigma \leq t \le \frac{cn(\lambda-1)^5}{\log^2 n}$, we have 
\begin{align}
\label{eqn:cadence}
\alpha^2_{t} = \lt(1+O\Bigg(\sqrt{\frac{\big(t+ \frac{\log^3n}{\lambda-1}\big)\log n}{n(\lambda-1)^5}}\Bigg)\rt)\alpha_{t}^{\star 2},
\end{align}
where $\{\alpha_t^{\star}\}$ stand for the asymptotic state evolution parameters defined in \eqref{eq:SE}. 
\end{itemize}
\end{theorem}
\begin{remark}[Range of $\lambda$]
Theorem~\ref{thm:main} only covers the regime where $\lambda$ is larger than but close to $1$. 
In fact, $\lambda=1$ represents the phase transition point for $\mathbb{Z}_{2}$ synchronization \citep{deshpande2017asymptotic}, in the sense that (i) when $\lambda < 1$, no estimator performs better than the $0$ estimator asymptotically, and (ii) when $\lambda$ is strictly larger than 1, it is possible to achieve non-trivial correlation with $\vstar$. 
We focus on the feasible regime by considering a more refined yet highly challenging case with  $\lambda -1 \gtrsim n^{-1/9}\log n$ (so that $\lambda$ can be very close to 1). 
While it is possible to improve the exponent $1/9$, it is beyond the scope of this paper. 
The upper bound $\lambda \leq 1.2$ is not crucial, and we make this assumption merely to simplify the presentation. 
\end{remark}


In the sequel, we provide some interpretations of Theorem~\ref{thm:main} and discussions about its implications. It is assumed below that $\lambda > 1$.

\paragraph{Gaussian approximation.} 
The first result \eqref{eqn:xt-decomposition} in Theorem~\ref{thm:main} asserts that 
each AMP iterate is composed of three components: (i) a signal component $\alpha_t\vstar$ that aligns with the true signal $\vstar$, 
(ii) a noise component $\sum_{k = 1}^{t-1} \beta_{t-1}^k\phi_k$ that is a linear combination of i.i.d.~Gaussian vectors, 
and (iii)  a residual component $\xi_{t-1}$. 
While this decomposition resembles that of \cite{li2022non}, we justify its validity even in the absence of carefully designed spectral initialization. 
A few remarks are in order. 
\begin{itemize}

	\item Regarding the noise component,  the 1-Wasserstein distance between its distribution (denoted by $\mu\big(\sum_{k = 1}^{t} \beta_{t}^k\phi_k\big)$) and a Gaussian distribution $\mathcal{N}\big(0, \frac{1}{n}I_n\big)$ is at most \citep[Lemma 9]{li2022non}
\begin{align}
	W_{1}\Bigg(\mu\bigg(\sum_{k=1}^{t}\beta_{k}\phi_{k}\bigg), \,\mathcal{N}\bigg(0,\frac{1}{n}I_{n}\bigg)\Bigg) \lesssim 
	\sqrt{\frac{t\log n}{n}}. 
\end{align}
For $t$ not too large, the noise component well approximates a Gaussian vector $\mathcal{N}(0,\frac{1}{n}I_{n})$. 

 \item Regarding the signal component $\alpha_t\vstar$, it is self-evident that $\alpha_t$ governs how effective AMP is in recovering the true signal. 
	 Importantly, once $\alpha_t$ exceeds the threshold $\frac{1}{2}\sqrt{\lambda^2 - 1}$, 
		it follows a non-asymptotic state evolution that closely resembles the asymptotic counterpart $\alpha_t^{\star}$ (see \eqref{eqn:cadence}), 
		a result that is made possible thanks to the non-asymptotic nature of our analysis.

\end{itemize}
To summarize, up to a small error term at most $\widetilde{O}\big(\sqrt{\frac{t}{n(\lambda-1)^2}} + \sqrt{\frac{1}{n(\lambda-1)^3}}\big)$, 
the AMP iterate is approximately
\begin{equation*}
	x_t \approx \alpha_t \vstar + \mathcal{N}\bigg(0,\frac{1}{n} I_n\bigg),
	\qquad t < O\bigg( \frac{n(\lambda - 1)^5}{\log^2 n}\bigg)
\end{equation*}
even when initialized randomly. An asymptotic version of this observation has been made in \cite{montanari2021estimation}, although the result therein required both informative initialization and a fixed $t$ that does not grow with $n$.

%

\paragraph{Dynamics after random initialization.} 
Compared to  \citet[Theorem 4]{li2022non}, the most challenging element of Theorem~\ref{thm:main} lies in analyzing the initial stage after random initialization. 
As shall be made clear from our analysis, we can understand the AMP trajectory by dividing it into three phases. 
\begin{itemize}
	\item {\bf Phase $\#1$: escape from random initialization.} 
		When initialized randomly with $x_{1} \sim \mathcal{N}(0,\frac{1}{n}I_n)$,  AMP starts with an extremely small signal component about the order of $\widetilde{O}(\frac{1}{\sqrt{n}})$, for which the canonical state evolution becomes vacuous.  
		To overcome this technical hurdle, we develop fine-grained characterizations regarding how  $\alpha_{t}$ evolves in this phase (before $|\alpha_t|$ surpasses $\sqrt{\lambda - 1}\, n^{-1/4}$), that is,  
	\begin{align}
		\alpha_{t+1} &\approx \lambda \alpha_t + \lambda g_{t-1}, \qquad \text{with }g_{t-1}\sim \mathcal{N}\Big(0,\frac{1}{n}\Big); 
	\end{align}
	see Section~\ref{sec:claim-alpha-t} for details.
	From this approximate noisy recursion, 
	while the signal component might be initially buried under the noise term, 
		it only takes $O(\frac{\log n}{\lambda - 1})$ iterations for the signal component to rise above  the noise size and reach the order of $\sqrt{\lambda - 1}\, n^{-1/4}$ (see Section~\ref{sec:alpha}). 

	\item {\bf Phase $\#2$: exponential growth.} 
	Once the signal component exceeds  $\sqrt{\lambda - 1} n^{-1/4}$ in size, the AMP iterate correlates non-trivially with the true signal. 
		Interestingly, the signal strength $\alpha_{t}$ starts to grow exponentially until reaching the order of $\sqrt{\lambda^2-1}.$ 
		As we shall justify in Section~\ref{sec:alpha}, $\alpha_{t+1}$ obeys  
	\begin{align} 
		\alpha_{t+1} \geq \sqrt{1 + \frac{1-o(1)}{3}(\lambda-1)} \,\, \alpha_t 
	\end{align}
	in this phase, which accounts for at most  $O(\frac{\log n}{\lambda - 1})$ iterations. 
	
	\item {\bf Phase $\#3$: local refinement.} 
	Upon reaching the order of $\sqrt{\lambda^2-1}$, $\alpha_{t}$ enters a local refinement phase, 
		during which randomly initialized AMP behaves similarly as AMP with spectral or other informative initialization. 
		In this phase, the asymptotic state evolution \eqref{eq:SE} also starts to be effective when predicting the evolution of $\alpha_{t}$ (see \eqref{eqn:cadence}).
	As we shall show in Section~\ref{sec:SE}, the signal strength  $\alpha_{t}$ satisfies
	\begin{align}
	\label{eqn:alpha-stage-3}
	  	|\alpha^2_{t+1} - \alpha^{\star 2}| \lesssim  \big(1 - (\lambda - 1)\big)^{t -\varsigma}
	  	+
			\widetilde{O}\Bigg(\sqrt{\frac{t + \frac{1}{\lambda-1}}{n(\lambda-1)^5}}\Bigg),
	  \end{align}  
	  where  $\alpha^\star $ (determined by $\lambda$) denotes the limit of $\alpha_{t}^{\star}$ (cf.~\eqref{eq:SE}) and is unique solution of 
	\begin{align}
	\alpha^{\star 2} &= \lambda^2 \Exs\big[\tanh\lt(\alpha^{\star}\lt(\alpha^{\star} + G\rt)\rt) \big],
	\qquad \text{for } G\sim \mathcal{N}(0,1).
		\label{eq:defn-alpha-star}
\end{align}

\end{itemize}

\paragraph{Asymptotic optimality.}
As was shown previously (see e.g.~\citet[Lemma A.7]{celentano2021local} and \citet[Remark 4]{li2022non}),  we can construct an AMP-based estimator whose risk coincides with that of the Bayes-optimal estimator $\Xbayes \defn \Exs[v^\star v^{\star\top} \mid M]$. More precisely, taking the AMP-based estimator as
\begin{align}
	u_t \defn \frac{1}{\lambda\sqrt{n(\alpha_t^2+1)}} \tanh(\pi_t x_t),
	\label{eq:defn-AMP-based-estimator}
\end{align}
its asymptotic risk satisfies (see Section~\ref{sec:pf-cor} and \citet{deshpande2017asymptotic}): 
\begin{align}
\label{eqn:opt}
	\lim_{t\to \infty} \lim_{n\to \infty} 
	\Exs\big[\big\|\vstar v^{\star \top} - u_t u_t^\top\big\|_{\mathsf{F}}^2\big]
	=
	\lim_{n\to \infty} \Exs\big[\|\vstar v^{\star \top} - \Xbayes\|_{\mathsf{F}}^2 \big]
	=
	1 - \frac{\alpha^{\star 4}}{\lambda^4}, 
\end{align}
where $\alpha^{\star}$ is the fixed point of the limiting state evolution  (cf.~\eqref{eq:defn-alpha-star}).   
This together with the non-asymptotic results in Theorem~\ref{thm:main} leads to a more refined risk characterization, as we shall prove in  Section~\ref{sec:pf-cor}.
\begin{cors}
\label{cor:opt}
With probability at least $1 - O(n^{-10})$, there exists some  $t=O(\frac{\log n}{\lambda - 1})$ such that 
	\begin{align}
	 \big\|\vstar v^{\star \top} - u_t u_t^\top\big\|_{\mathsf{F}}^2 
	 =  1 -  \frac{\alpha^{\star 2}}{\lambda^4} +  O\bigg(\sqrt{\frac{\log^4 n}{n(\lambda - 1)^6}}\bigg). 
	\end{align}

\end{cors}


%
%
%

\section{Proof of Theorem~\ref{thm:main}}

In this section, we present the proof of our main result: Theorem~\ref{thm:main}. 
We find it helpful to introduce the following notation that helps streamline the presentation: 
\begin{align}
\label{defn:v-t-thm1} 
v_t \defn \alpha_t \vstar + \sum_{k=1}^{t-1} \beta_{t-1}^k \phi_k.
\end{align}

\subsection{Preliminaries}

Before we embark on our proof of the main theorem, we collect a couple of useful results that shall be used frequently throughout this proof.

\paragraph{Concentration results.}  We first record several useful concentration results from \cite{li2022non}.   
Here and throughout, we let $|x|_{(i)}$ denote the magnitude of the $i$-th largest entry (in magnitude) of $x\in \real^{n}$.  
%
\begin{lemma}
\label{lem:concentration}
Consider a collection of random vectors $\{\phi_k\}_{1\leq k < t}$ in $\mathbb{R}^n$.  
Suppose that for each $1\leq k\leq t-1 < n$, $\phi_k=(\phi_{k,1},\phi_{k,2},\ldots,\phi_{k,n})$ is i.i.d.~drawn from $\mathcal{N}\left(0,\frac{1}{n}I_n\right)$.  
Consider the following set
\begin{align*}
\mathcal{E}_s& :=
	\left\{\left\{\psi_k\right\}_{k=1}^{t-1}: \max _{1 \leq k \leq t-1}\left\|\psi_k\right\|_2<1+C \sqrt{\frac{\log \frac{n}{\delta}}{n}}\right\} \bigcap
\left\{\left\{\psi_k\right\}_{k=1}^{t-1}: \sup _{a \in \mathcal{S}^{t-2}}\Big\|\sum_{k=1}^{t-1} a_k \psi_k\Big\|_2<1+C \sqrt{\frac{t \log \frac{n}{\delta}}{n}}\right\} \\
&\bigcap\left\{\left\{\psi_k\right\}_{k=1}^{t-1}: \sup _{a\in \mathcal{S}^{t-2}} \sum_{i=1}^s\Big|\sum_{k=1}^{t-1} a_k \psi_k\Big|_{(i)}^2<\frac{C(t+s) \log \frac{n}{\delta}}{n}\right\}
\bigcap\left\{\left\{\psi_k\right\}_{k=1}^{t-1}:\max_{1\leq k< t,1\leq i\leq n}\vert\psi_{k,i}\vert<C\sqrt{\frac{\log\frac{n}{\delta}}{n}}\right\}, \\
\end{align*}
and denote $\mathcal{E} :=\bigcap_{s=1}^{n}\mathcal{E}_s$. 
Then there exists some large enough constant $C>0$  such that, for every $\delta > 0$, 
\begin{align*}
	\mprob\big(\left\{\phi_k\right\}_{k=1}^{t-1}\in\mathcal{E}\big)\geq 1-\delta.
\end{align*}
In particular, by setting $\delta=n^{-11}$, we see that the following event happens with probability at least $1-O(n^{-11})$: 
\begin{align*}
	\left|\max _{1 \leq k \leq t-1}\left\|\phi_k\right\|_2-1\right|\lesssim\sqrt{\frac{\log n}{n}},
	\qquad
	\text{and }~\max_{1\leq k\leq t,1\leq i\leq n}\vert\phi_{k,i}\vert\lesssim\sqrt{\frac{\log n}{n}}.
\end{align*} 
\end{lemma}
\noindent This lemma is a consequence of standard concentration of measure for Gaussian random vectors \citep{massart2007concentration}; its proof can be found in \citet[Section D.1.1]{li2022non} and is hence omitted for brevity.

\paragraph{Properties of $\eta_t$, $\pi_t$, and $\gamma_t$ (cf.~\eqref{defi:eta}).} 
Next, we summarize several basic properties about the three sets of key quantities defined in \eqref{defi:eta}. 
We begin by gathering several basic properties for our choices of $\pi_t$ and $\gamma_t$ defined in \eqref{defi:eta}; the proof is deferred to Section~\ref{sec:proof-parameters}.
\begin{lemma}
\label{lem:parameters}
Suppose the decomposition~\eqref{eqn:xt-decomposition} is valid with $\|\xi_{t-1}\|_2 \lesssim 1$. 
With probability at least $1 - O(n^{-10})$, the following properties hold true:    
\begin{subequations}
\begin{align} 
\label{eqn:pi_t}
\frac{1}{\sqrt{n}}\pi_t &
= |\alpha_t| + O\Bigg(\bigg(\|\xi_{t-1}\|_2 + \sqrt{\frac{t\log n}{n}}\bigg)^{1/2} \,\wedge\, \frac{1}{|\alpha_t|}\bigg(\|\xi_{t-1}\|_2 + \sqrt{\frac{t\log n}{n}}\bigg) \Bigg); \\
\label{eqn:gamma_t}
\gamma_t^{-2} & = n\int \tanh^2\lt(\frac{\pi_t}{\sqrt{n}}(\alpha_t + x)\rt)\varphi(\dx) + \pi_t^2O\bigg(\Vert \xi_{t-1}\Vert_2 + \sqrt{\frac{t\log n}{n}}\bigg).
\end{align}
Additionally, one has 
\begin{align} 
\label{eqn:tanh_pi_t}
\int \tanh^2\lt(\frac{\pi_t}{\sqrt{n}}(\alpha_t + x)\rt)\varphi(\dx) &= \frac{\pi_t^2}{n}(\alpha_t^2 + 1) + O\lt(\frac{\pi_t^4}{n^2}\rt),
\end{align}
which in turn implies that 
\begin{align}
\label{eqn:gamma_t_arrange}
\gamma_t^{-2} & = \pi_t^2\bigg(\alpha_t^2+1+O\bigg(\frac{\pi_t^2}{n}+\Vert\xi_{t-1}\Vert_2+\sqrt{\frac{t\log n}{n}}\bigg)\bigg).
\end{align}
\end{subequations}
\end{lemma}

Next, let us single out several properties about the quantity $\eta_t$ in the lemma below, whose proof is postponed to Section~\ref{sec:proof-eta_property}.
\begin{lemma} 
\label{lem:eta_property}
Consider any $1\leq t\leq n$ and suppose that $\ltwo{\xi_{t-1}}$ satisfies
\begin{align}
\|\xi_{t-1}\|_2  &\lesssim \sqrt{\frac{t^3\log n}{n}}\qquad\text{for all }t \lesssim \frac{\log n}{\lambda - 1}.
	\label{eq:induction-lemma-xi}
\end{align}
\begin{subequations}
Then the following properties hold true with probability at least $1 - O(n^{-11})$:
\begin{itemize}
	\item If $t \lesssim \frac{n(\lambda-1)^4}{\log^2 n}$ and $\|\xi_{t-1}\|_2 \lesssim \frac{1}{\sqrt{\log n}}$, then any $x \in \mathbb{R}$ obeys 
\begin{align}
\label{eq:eta_derivative}
	|\eta_t (x) | \lesssim |x|, \qquad
	|\eta_t^{\prime} (x) | \lesssim 1 \eqqcolon \rho, 
	\qquad |\eta_t^{\second} (x)| \lesssim \sqrt{n} \eqqcolon \rho_1, 
	\qquad  |\eta_{t}^{(\third)} (x)| \lesssim n \eqqcolon \rho_2;
\end{align}

\item If $t \lesssim \frac{\log n}{\lambda-1}$ and $\alpha_t\lesssim \sqrt{\lambda-1} \, n^{-0.1}$, then one has 
$\pi_t\lesssim\sqrt{\lambda-1}\, n^{0.4}$; and for any $x$ obeying $|x| \lesssim \sqrt{\log n/n}$, we have
\begin{align} 
\label{eq:eta_derivative_small}
\eta_{t}^{\prime}(x) = 1 + O\big((\lambda-1)n^{-0.2}\log n \big)\qquad\text{and}\qquad |\eta_{t}^{\second}(x)| \lesssim (\lambda-1)n^{0.8}|x|;
\end{align}

\item If $t \lesssim \frac{\log n}{\lambda-1}$ and $\alpha_t \lesssim \sqrt{\lambda-1}\,n^{-1/4}$, 
then one has $\pi_t \lesssim (\lambda-1)^{-3/4}n^{1/4}\log n$; 
		for any $x$, one can find a quantity $c_0 \lesssim (\log^2 n)/\sqrt{n(\lambda-1)^{3}}$ independent from $x$, and another quantity $|c_x| \lesssim n|x|^5(\log^4 n)/(\lambda-1)^3$ depending on $x$, such that
\begin{align} 
\label{eq:eta_small}
\eta_t(x) = (1 - c_0) \lt(x - \frac{1}{3}\pi_t^2x^3 + c_x\rt).
\end{align}

\end{itemize}
\end{subequations}
\end{lemma}

%

\subsection{Non-asymptotic analysis for the AMP dynamics}

We are now in a position to present the proof of our main theorem. The structure of our proof is outlined in what follows. 
\begin{itemize}
	\item Firstly, focusing on the initial stage obeying $t \le \varsigma\wedge \frac{\log n}{c(\lambda - 1)}$, 
		we develop an upper bound on $\ltwo{\xi_{t}}$ in Section~\ref{sec:xi} as follows: 
		\begin{align}
\|\xi_{t}\|_2  &\lesssim \sqrt{\frac{t^3\log n}{n}};
\label{eq:induction}
\end{align}
		here, $\varsigma$ is a threshold defined in \eqref{eq:defn-varsigma} and $c > 0$ is some constant small enough. 
		This is accomplished by means of an inductive argument.
	\item Secondly,  we investigate in Section~\ref{sec:alpha} how the signal strength $\alpha_{t}$ evolves during the execution of AMP. 
		Crucially, recalling that $\varsigma$ reflects the first time $t$ that satisfies $|\alpha_{t}| \gtrsim \sqrt{\lambda^2 - 1}$ (cf.~\eqref{eq:defn-varsigma}), 
		we demonstrate that  
\begin{align}
\label{claim:alpha}
	\varsigma \lesssim \frac{\log n}{\lambda-1}; 
\end{align}
		in words, in spite of random (and hence uninformative) initialization, it takes AMP at most $O\big(\frac{\log n}{\lambda-1}\big)$ iterations to find an informative estimate.

	\item 
Thirdly, with the above control of $\varsigma$ in place, we go on to develop a more complete upper bound on $\ltwo{\xi_{t}}$ that covers the iterations after $\varsigma$, that is,   
\begin{align}
\|\xi_{t}\|_2  &\lesssim \sqrt{\frac{t\ind(t > \varsigma)\log n}{n(\lambda-1)^2}} + \sqrt{\frac{\min\{t, \varsigma\}^3\log n}{n}}  
\label{eq:induction-main}
\end{align}
for any $t < \frac{c n(\lambda-1)^5}{\log^2 n}$.   
This is the main content of Section~\ref{sec:more-complete-xit}, accomplished again via an inductive argument. 

	\item Finally, after the iteration number exceeds the threshold $\varsigma$, 
		we demonstrate in Section~\ref{sec:SE} that the asymptotic state evolution (the one characterizing large-system limits) becomes fairly accurate in the finite-sample/finite-time regime. In particular, an intimate connection is established between the non-asymptotic state evolution and its asymptotic analog, which plays a critical role in characterizing the finite-sample convergence behavior of AMP. 
\end{itemize}

\noindent These four steps will be explained in detail in the sequel.

\subsubsection{Controlling $\xi_t$ when $t \le \varsigma\wedge \frac{\log n}{c(\lambda - 1)}$ (Proof of Claim~\eqref{eq:induction})}
\label{sec:xi} 




In this subsection, we establish the claimed bound \eqref{eq:induction} for $\|\xi_{t}\|_2$, which leverages on ideas from \citet{li2022non}. 
To begin with, let us restate \citet[Theorem 2]{li2022non} below, with slight simplification tailored to $\mathbb{Z}_2$ synchronization (i.e., through the use of the properties $\|\beta_t\|_2^2 = 1$, $E_t = 0$, and \eqref{eq:eta_derivative}). For notational convenience, define $\kappa_{t} > 0$ such that  
\begin{align}
\begin{aligned}
\label{eqn:kappa-def}
\kappa_t^2 \defn \max\bigg\{\bigg\langle\int\left[ x\eta_t^{\prime}\left(\alpha_tv^\star+\frac{1}{\sqrt{n}}x\right)-\frac{1}{\sqrt{n}}\eta_t^{\prime\prime}\left(\alpha_tv^{\star}+\frac{1}{\sqrt{n}}x\right)\right]^2\varphi_n(\dx)\bigg\rangle, \\
\qquad\qquad\bigg\langle\int\left[\eta_t^{\prime}\left(\alpha_tv^\star+\frac{1}{\sqrt{n}}x\right)\right]^2\varphi_n(\dx)\bigg\rangle\bigg\},
\end{aligned}
\end{align}
where we recall that $\varphi_n(\cdot)$ is the p.d.f.~of $\mathcal{N}\left(0,I_n\right)$ and for any vector $x=[x_i]_{1\leq i\leq n}$, we denote $\langle x\rangle := \frac{1}{n}\sum_{i=1}^{n}x_i$ and $x^2=[x_i^2]_{1\leq i\leq n}$. We shall work with the following assumptions. 

\begin{assumption} 
\label{assump:A-H-eta}
For any $1\leq t \leq n$, consider arbitrary vectors $\mu_t \in \mathcal{S}^{t-1}$, $\xi_{t-1}\in \real^n$, and coefficients $(\alpha_t, \beta_{t-1}) \in \real \times \real^{t-1}$ that might all be statistically dependent on $\phi_{k}$. Let $v_t$ be defined as in \eqref{defn:v-t-thm1}. 
We assume the existence of (possibly random) quantities $A_{t}, B_t, D_{t}$ such that with probability at least $1-O(n^{-11})$, the following inequalities hold: 
\begin{subequations}
\begin{align}
	\label{defi:A}
	\bigg|\sum_{k = 1}^{t-1} \mu_t^k\Big[\big\langle \phi_k, \eta_{t}(v_t)\big\rangle - \big\langle\eta_t^{\prime}(v_t)\big\rangle \beta_{t-1}^k\Big]\bigg| &\,\le\, A_t, \\
	\bigg|v^{\star\top}\eta_{t}(v_t) - v^{\star\top}\int\eta_t\Big(\alpha_t \vstar + \frac{\|\beta_{t-1}\|_2}{\sqrt{n}}x\Big)\varphi_n(\dx)\bigg| &\,\le\, B_t, \label{defi:B}\\
%
%
	\bigg\|\sum_{k = 1}^{t-1} \mu_t^k\phi_k \circ \eta_{t}^{\prime}(v_t) - \frac{1}{n}\sum_{k = 1}^{t-1} \mu_t^k\beta_{t-1}^k\eta_{t}^{\second}(v_t)\bigg\|_2^2 - \kappa_t^2 &\,\le\, D_t. \label{defi:D}
%
\end{align}
\end{subequations}
\end{assumption}

Under these assumptions, \citet[Theorem 2]{li2022non} developed a general non-asymptotic characterization for AMP iterates as follows.
\begin{theorem}{[Adapted from \citet[Theorem 2]{li2022non}]}
\label{thm:main-AMP}
Suppose that Assumption~\ref{assump:A-H-eta} holds, and consider any  $t\leq n$.  
With probability at least $1-O(n^{-11})$, 
the AMP iterates \eqref{eqn:AMP-updates} for $\mathbb{Z}_2$-synchronization satisfy the decomposition~\eqref{eqn:xt-decomposition} with $\|\beta_t\|_2^2 = 1$ and
\begin{align}
\label{eqn:alpha-t-genearl}
	\alpha_{t+1} &= \lambda v^{\star \top} \int{\eta}_{t}\lt(\alpha_t \vstar + \frac{1}{\sqrt{n}}x\rt)\varphi_n(\dx) + \Delta_{\alpha,t} 
\end{align}
where the residual terms obey
\begin{subequations}
\label{eqn:para-general}
\begin{align}
\label{eqn:delta-alpha-general}
|\Delta_{\alpha,t}| &\,\lesssim\, B_t + \big|v^{\star\top}\eta_{t}(x_t) - v^{\star\top}\eta_{t}(v_t)\big| \lesssim B_t + \|\xi_{t-1}\|_{2}, \\
\label{eqn:xi-t-general}
	\|\xi_{t}\|_{2}&\le\sqrt{\kappa_{t}^{2}+D_{t}}\,\|\xi_{t-1}\|_{2}+O\Bigg(\sqrt{\frac{t\log n}{n}}+A_{t}+\sqrt{(1+t\ind_{t \le \varsigma})\log n}\,\|\xi_{t-1}\|_{2}^{2}+\sqrt{\frac{{t\log n}}{n}}\|\xi_{t-1}\|_{2}\Bigg) .
\end{align}
\end{subequations}
\end{theorem}
\begin{remark}
With regards to the above bound \eqref{eqn:xi-t-general} for $\|\xi_{t}\|_{2}$, a direct application of \citet[Theorem 2]{li2022non} 
results in a term $\sqrt{t\log n} \,\|\xi_{t-1}\|_{2}^{2}$ (as opposed to $\sqrt{(1+t\ind_{t \le \varsigma})\log n}\,\|\xi_{t-1}\|_{2}^{2}$ in \eqref{eqn:xi-t-general}).  
We make slight modifications here to make it better-suited for the current setting. 
\begin{itemize}
	\item[(i)] When $t \leq \varsigma$, such a term $\sqrt{t\log n} \,\|\xi_{t-1}\|_{2}^{2}$ works fine for our purpose; 
	\item[(ii)] When $t > \varsigma$ (so that $\alpha_{t}$ exceeds the order of $\sqrt{\lambda^2 - 1}$), 
		one can simply invoke \citet[display (249)]{li2022non} to improve the factor in front of $\|\xi_{t-1}\|_{2}^{2}$ from $\sqrt{t\log n}$ to $\sqrt{\log n}.$ 
\end{itemize}
Putting these together leads to the claimed bound \eqref{eqn:xi-t-general}. 
Notably, this seemingly minor change turns out to be essential in order to push the number of iterations to $O(n/\mathsf{poly}(\log n))$ instead of $O(\sqrt{n}/\mathsf{poly}(\log n)).$ 
\end{remark}

With Theorem~\ref{thm:main-AMP} in mind, in order to control $|\Delta_{\alpha,t}|$ and $\|\xi_{t}\|_{2}$, it boils down to determining $A_t, B_t, D_t$, and $\kappa_t$, respectively.
\begin{itemize}
	\item {\bf Bounding  $A_t, B_t,$ $D_t$.}
Repeating the same analysis as in~\citet[Section D.2]{li2022non}, we obtain 
\begin{align}
	\label{eqn:ABD}
	A_t \lesssim \sqrt{\frac{t\log n}{n}},
	\qquad B_t \lesssim \sqrt{\frac{t\log n}{n}},
	\qquad D_t \lesssim \sqrt{\frac{t\log^2 n}{n}}.
\end{align}
The only term that needs more discussion is $A_t$, 
		as \citet{li2022non} only proved that $A_t \lesssim \frac{1}{\alpha_t} \sqrt{\frac{t\log n}{n}}$ (taking $s = 1$ therein) for AMP with independent initialization.  
To get rid of the prefactor $1/\alpha_t$, we rely on  an improved control of $\eta_t(x)$ (cf.~\eqref{eq:eta_derivative}). 
In particular, property \eqref{eq:eta_derivative} tells us that
\begin{align*}
	\|\eta_t(v_t)\|_2 \lesssim \|v_t\|_2 = 
	\Big\|\alpha_{t} v^{\star{}} + \sum_{k = 1}^{t-1} \beta_{t-1}^k\phi_k\Big\|_2 \lesssim 1, 
\end{align*}
where the last inequality can be found in display~\eqref{eqn:v-norm}. In turn, this leads to 
\begin{align}
	\|\nabla_{\Phi} f_{\theta}(\Phi)\|_2 \lesssim \frac{1}{\sqrt{n}}
\end{align}
through the same analyses as detailed around \citet[Section D.2.1, inequality (229)]{li2022non}. 
Here, $\nabla_{\Phi} f_{\theta}(\Phi)$ is the key quantity to control $A_t$ in~\citet[Section D.2.1]{li2022non}, and our desired bound for $A_t$ follows immediately. 
Given that this only consists of very minor and straightforward changes to \citet[Section D.2.1]{li2022non}, 
we omit the details for brevity and refer the readers to \citet[Section D.2.1]{li2022non} for more details. 

\item {\bf Bounding $\kappa_t$.} 
The main step then comes down to bounding $\kappa_t$. Towards this end, we claim that the following relation holds for $\kappa_t$, whose proof is postponed to Section~\ref{sec:pf-kappa}.
\begin{lemma}
\label{lem:kappa}
With probability at least $1-O(n^{-10})$, the following results hold true: 
\begin{itemize}
	\item Under the inductive assumption~\eqref{eq:induction} for $\xi_{t-1}$, one has 
\begin{subequations}
\label{eq:kappa-claim}
\begin{align}
\kappa_t\le 
1+o\Big(\frac{\lambda-1}{\log n}\Big)  \label{eq:kappa-claim-init}
\end{align}
provided that $t \le \varsigma\wedge \frac{\log n}{c(\lambda - 1)}$; 

\item Under the inductive assumption~\eqref{eq:induction-main} for $\xi_{t-1}$, one has 
\begin{align}
\kappa_{t} \le 1-\frac{1}{15} (\lambda-1),   \label{eq:kappa-claim-main}
\end{align}
provided that $t \le \frac{cn(\lambda-1)^5}{\log^2 n}$ and $|\alpha_t| \gtrsim \sqrt{\lambda^2 - 1}$.
\end{subequations}
\end{itemize}
\end{lemma}

\end{itemize}



With the above estimates of $A_t,B_t,D_t,\kappa_t$ in place, we are ready to apply Theorem~\ref{thm:main-AMP}. 
Under the inductive assumption~\eqref{eq:induction}, the recursive formula~\eqref{eqn:para-general} in Theorem~\ref{thm:main-AMP} taken together with \eqref{eqn:ABD} yields 
\begin{align}
\label{eqn:christmas}
\|\xi_{t}\|_{2} 
&\le \sqrt{\kappa_t + \sqrt{\frac{t\log^2 n}{n}}} \|\xi_{t-1}\|_{2}+O\Bigg(\sqrt{\frac{t\log n}{n}} +\sqrt{(1+t\ind_{t \le \varsigma})\log n}\,\|\xi_{t-1}\|_{2}^{2}+\sqrt{\frac{{t\log n}}{n}}\|\xi_{t-1}\|_{2}\Bigg),
\end{align}
which combined with \eqref{eq:kappa-claim}  further implies that
\begin{align}
\label{eqn:error-final-1}
\|\xi_{t}\|_{2} \le 
\lt(1+o\Big(\frac{\lambda-1}{\log n}\Big) + O\bigg(\sqrt{\frac{t^4\log^2 n}{n}}\bigg)\rt)\|\xi_{t-1}\|_{2} + O\Big(\sqrt{\frac{t\log n}{n}}\Big) ,
%
\end{align}
with the proviso that $t \le \varsigma\wedge \frac{\log n}{c(\lambda - 1)}$.

We are now ready to prove relation~\eqref{eq:induction} via induction.  
To verify its validity for the base case (i.e.~$t=1$), we note that by construction (see, e.g.~\cite[Step 3, Proof of Theorem 1]{li2022non}), $\xi_{1}$ takes the form
\begin{align*}
	\xi_1 = \Big(\frac{\sqrt{2}}{2}-1\Big)z_1 z_1^\top W z_1, \quad \text{where } z_1 = \eta_1(x_1) \text{ is independent of } W.
\end{align*}
Elementary calculations reveal that,  with probability at least $1 - O(n^{-11})$, 
\begin{align}
\label{eqn:initial}
 	\ltwo{\xi_1} = \Big|\frac{\sqrt{2}}{2}-1\Big|\cdot\ltwo{z_1} \cdot|z_1^\top W z_1| \lesssim \sqrt{\frac{\log n}{n}},
 \end{align} 
 given that $\ltwo{z_1} = 1$ and $z_1^\top W z_1 \sim \mathcal{N}(0,\frac{2}{n}I_n)$. This already establishes \eqref{eq:induction} for the base case with $t=1$. 
 Next, consider the case where $t \le \varsigma\wedge \frac{\log n}{c(\lambda - 1)}$ for some small enough constant $c>0$. 
 Given that $\sqrt{\frac{t^4\log^2 n}{n}} = o(\frac{\lambda-1}{\log n})$ under our assumption on $\lambda - 1$, 
 the recursive relation~\eqref{eqn:error-final-1} immediately leads to 
\begin{align}
\label{eq:xi-stage1}
	\|\xi_t\|_2 &\le
	\notag \lt(1+o\Big(\frac{\lambda-1}{\log n}\Big)\rt)\|\xi_{t-1}\|_{2} + O\Big(\sqrt{\frac{t\log n}{n}}\Big) \\
	\notag &\le \lt(1+o\Big(\frac{\lambda-1}{\log n}\Big)\rt)^{t-1} \ltwo{\xi_1} + 
	\sum_{j=0}^{t-2} \lt(1+o\Big(\frac{\lambda-1}{\log n}\Big)\rt)^{j} O\Big(\sqrt{\frac{(t-j)\log n}{n}}\Big)
	\\
	&\lesssim \sqrt{\frac{t^3\log n}{n}}
\end{align}
for all $t \le \varsigma\wedge \frac{\log n}{c(\lambda - 1)}$, as claimed.

\begin{remark}
	Careful readers might note that the recursive formula established in \eqref{eqn:error-final-1} for $t \le \varsigma\wedge \frac{\log n}{c(\lambda - 1)}$ does not rely on the relation~\eqref{claim:alpha} (a relation that shall be established in the next subsection). 
\end{remark}




\subsubsection{Evolution of $\alpha_t$ and a bound on $\varsigma$ (Proof of Claim~\eqref{claim:alpha})}
\label{sec:alpha}

 

We now move on to establish the claim~\eqref{claim:alpha} concerning an upper bound on the threshold $\varsigma$, 
which requires careful analysis about how the signal strength $\alpha_t$ evolves at the initial stage. 
Towards this end,  we divide into two cases based on the magnitude of $\alpha_{t}$, 
which we shall detail after presenting several preliminary facts. 


\paragraph{Preliminary facts.} 
Before proceeding, we first recall some additional preliminary facts already established in \citet{li2022non}. From the analysis of \citet[Theorem 1]{li2022non}, we know that: by construction, 
\begin{align}
	\xi_{t} \in \mathsf{span}(U_{t-1}) = \mathsf{span}\Big\{\eta_1(x_1),\ldots,\eta_{t-1}(x_{t-1})\Big\},
	\label{eq:span-Ut-eta-t-equal}
\end{align}
where $U_{t-1} \in \real^{n\times (t-1)}$ is a matrix whose columns are formed by a set of orthonormal basis $\{z_1, \ldots, z_{t-1}\}$. 
In fact, we can specify $U_{t}$ in a more explicit manner. Following \citet[Section 4.1]{li2022non}, let us define  
\begin{subequations}
\label{eqn:z-w-recursion}
\begin{align}
\label{eqn:z-w-init}
	z_1 \defn \frac{\eta_1(x_1)}{\lt\|\eta_1(x_1)\rt\|_2} \in \real^n
	\qquad\text{and}\qquad 
	W_1 \defn W \in \real^{n\times n},
\end{align}
which are statistically independent from each other; 
and then any $2 \leq t \leq n$, we can define the following objects recursively:  
\begin{equation}
	U_{t-1} \defn [z_k]_{1 \le k \leq t-1} \in \real^{n\times (t-1)}, \label{eq:defn-Ut-minus-1}
\end{equation}
and also 
\begin{align}
	z_t &\defn \frac{\lt(I_n - U_{t-1}U_{t-1}^{\top}\rt)\eta_{t}(x_{t})}{\lt\|\lt(I_n - U_{t-1}U_{t-1}^{\top}\rt)\eta_{t}(x_{t})\rt\|_2}, 
	\label{eqn:zt}\\
	W_t &\defn \lt(I_n - z_{t-1}z_{t-1}^{\top}\rt)W_{t-1}\lt(I_n - z_{t-1}z_{t-1}^{\top}\rt),
	\label{eqn:Wt}
\end{align}
where $\{x_t\}$ is the sequence generated by the AMP updates~\eqref{eqn:AMP-updates}. 
\end{subequations}
This process thus leads to more explicit forms for $\{U_t\}$ and the orthonormal basis $\{z_t\}$ (see \citet[Section 4.1]{li2022non} for the orthonormality of $\{z_t\}$). What is more, the orthonormality of $\{z_t\}$ reveals the decomposition 
\begin{align}
	\eta_t(x_t) = \sum_{k=1}^t \beta_t^k z_k, \qquad \text{with } \beta_t^k := \inprod{\eta_t(x_t)}{z_k},
	\label{eq:eta-t-expand-zt}
\end{align}
which satisfies $\ltwo{\eta_t(x_t)} = \ltwo{\beta_t}$ with $\beta_t=[\beta_t^1,\ldots,\beta_t^{t}]$. 
Additionally, we find it convenient to generate 
\begin{align}
\label{eqn:nutcracker}
	\phi_k &\defn W_kz_k + \zeta_k, \qquad \text{where }\zeta_k \defn \Big(\frac{\sqrt{2}}{2} - 1\Big) z_k^{\top}W_kz_k \cdot z_k + \sum_{i = 1}^{k - 1} g_i^kz_i, 
	\qquad 1\leq k\leq n,
\end{align}
where the $g_i^k$'s are independently drawn from $\mathcal{N}(0, \frac{1}{n})$.
The following properties have been shown in \citet[Lemma 2]{li2022non}, which play a crucial role in our subsequent analysis: 
\begin{itemize}
	\item $\phi_k\stackrel{\mathrm{i.i.d}}{\sim }\mathcal{N}(0, \frac{1}{n} I_n)$, for $1\leq k\leq n$;
	\item The randomness of $\phi_{k}$ only comes from $W_{k}$,  and  $\phi_{k}$ is independent of $x_{1}$ and $\{z_i\}_{i < k}$.
	\item $x_k$ and $z_k$ are conditionally independent from $W_{k}$ given $\{z_i\}_{i < k}$ and $x_{1}$.
	\item $\phi_{k}$ is independent from $\{x_{j}\}_{j\leq k}$ and $\{z_{j}\}_{j\leq k}$.
\end{itemize}

\paragraph{Stage I: small correlation ($|\alpha_t| \lesssim \sqrt{\lambda-1}\,n^{-1/4}$).} 
Let us start from the very beginning when the correlation coefficient $\alpha_t$ is reasonably small. 
Towards this, we define a threshold $\tau_0$ such that
\begin{equation}
	\tau_0 := \max \big\{ \tau: |\alpha_t| \lesssim \sqrt{\lambda-1}\,n^{-1/4} \text{ for all } t\leq \tau \big\} ;
	\label{eq:defn-tau0-threshold}
\end{equation}
in words, $\tau_0+1$ represents the first term that exceeds the level of $\sqrt{\lambda-1}\,n^{-1/4}$. 
In the following, we would like to prove that, with probability at least $1 - O(n^{-10})$, this threshold is not too large in the sense that
\begin{align}
\label{eqn:tau0}
	\tau_0 \lesssim \frac{\log n}{\lambda-1}. 
\end{align}

\begin{proof}[Proof of Claim~\eqref{eqn:tau0}] 
In order to establish this result \eqref{eqn:tau0}, 
we first state an important claim:  the AMP iterates --- when initialized at a random point --- satisfy the following recursive relation with high probability:  
\begin{align}
	\alpha_{t+1} = \lambda^{t-k+1}\alpha_k + \sum_{i = 1}^{t-k+1} \lambda^ig_{t-i} + O\lt(\lambda^{t-k}\frac{\log^4 n}{n^{3/4}(\lambda-1)^{1.5}}\rt) \label{eq:alpha_t}
\end{align}
for any $1 \le k \le t$, 
where we denote 
\begin{equation}
	g_k \defn v^{\star\top}\phi_{k}~~ (1 \le k \le t) \qquad \text{and} \qquad g_0 = 0.
\end{equation}
This claimed relation lies at the heart of the analysis for Stage I, 
in which the correlation between the AMP iterate and $\vstar$ keeps growing to a non-trivial value. 
To streamline the presentation, we defer the proof of this claim to Section~\ref{sec:claim-alpha-t}.

Equipped with the above recursive formula \eqref{eq:alpha_t}, we now turn to proving the relation~\eqref{eqn:tau0}. 
Define $t_i := C' i\log n$ for some quantity $C' = \frac{C''}{\lambda-1}$, where $C''$ is some large enough constant.
Observe that 
\begin{align*}
\mathbb{P}\Big(|\alpha_k| \lesssim \frac{\sqrt{\lambda-1}}{n^{1/4}},\text{ for all }k \le 201C'\log n\Big) 
&\le \mathbb{P}\Big(|\alpha_{t_i + 1}| \lesssim \frac{\sqrt{\lambda-1}}{n^{1/4}},\text{ for all } 1\le i \le 200\Big) \\
&= \prod_{i = 1}^{200} \mathbb{P}\Big(|\alpha_{t_i + 1}| \lesssim \frac{\sqrt{\lambda-1}}{n^{1/4}} ~\Big\vert~ |\alpha_{t_{j} + 1}| \lesssim \frac{\sqrt{\lambda-1}}{n^{1/4}}, ~\forall 1\leq j < i\Big).
\end{align*}
To control the right-hand side of the above relation,  consider the following random variable 
\begin{align}
\label{eqn:xi-dist}
	X_i := \sum_{j = 1}^{C'\log n} \lambda^jg_{t_i-j} \sim \mathcal{N}\lt(0, \frac{\lambda^{2C'\log n+2} - \lambda^2}{n(\lambda^2 - 1)}\rt)
\end{align}
for each $1 \leq i\leq 200$. 
Armed with this piece of notation, invoking relation~\eqref{eq:alpha_t} gives 
\begin{align}
	\alpha_{t_i + 1} = \lambda^{t_i - t_{i-1}} \alpha_{t_{i-1} + 1} + X_{i} + O\lt(\lambda^{t_i - t_{i-1}-1}\frac{\log^4 n}{n^{3/4}(\lambda-1)^{1.5}}\rt).
\end{align}
As mentioned in the above preliminary facts, each $\phi_j$ is independent with the AMP iterate $x_i$ for  $i \le j$ and therefore $\alpha_{i+1}$, given that $\alpha_{i+1} = v^{\star \top} \eta_{i}(x_i).$
As a result, the random variable $X_i$ defined above is independent from $\alpha_{t_{j} + 1}$ for all $j \leq i-1$. 
Taking this together with the relation \eqref{eq:alpha_t} then leads to
\begin{align}
\notag &\mathbb{P}\Big(|\alpha_{t_i + 1}| \lesssim \frac{\sqrt{\lambda-1}}{n^{1/4}} ~\Big\vert~ |\alpha_{t_{j} + 1}| \lesssim \frac{\sqrt{\lambda-1}}{n^{1/4}}, 1\le j < i\Big) \\
\notag &\le \mathbb{P}\Big(|\lambda^{t_i - t_{i-1}} \alpha_{t_{i-1} + 1} + X_{i}| \lesssim \frac{\sqrt{\lambda-1}}{n^{1/4}} + \lambda^{t_i - t_{i-1} - 1}\frac{\log^4 n}{n^{3/4}(\lambda-1)^{1.5}} ~\Big\vert~ |\alpha_{t_{j} + 1}| \lesssim \frac{\sqrt{\lambda-1}}{n^{1/4}} , 1\le j < i\Big) \\
%
%
\notag &= \mathbb{P}\Big(|\lambda^{t_i - t_{i-1}} \alpha_{t_{i-1} + 1} + X_{i}| \lesssim \frac{\sqrt{\lambda-1}}{n^{1/4}} + \lambda^{C'\log n}\frac{\log^4 n}{n^{3/4}(\lambda-1)^{1.5}}~\Big\vert~ |\alpha_{t_{i-1} + 1}| \lesssim \frac{\sqrt{\lambda-1}}{n^{1/4}}\Big)\\
&\le \mathbb{P}\Big(|X_{i}| \lesssim \frac{\sqrt{\lambda-1}}{n^{1/4}} + \lambda^{C'\log n}\frac{\log^4 n}{n^{3/4}(\lambda-1)^{1.5}}\Big) .
\end{align}
Here, the penultimate line follows from the independence relation stated above, whereas the last line follows from the elementary fact that
\begin{align*}
	\mathbb{P}_{X \sim \mathcal{N}(\mu, \sigma^2)}\lt(|X| < x\rt) \le \mathbb{P}_{X \sim \mathcal{N}(0, \sigma^2)}\lt(|X| < x\rt),
	\qquad \forall x > 0.
\end{align*}
Putting these pieces together, we conclude that  
\begin{align}
\notag \mathbb{P}\Big(|\alpha_k| \lesssim \frac{\sqrt{\lambda-1}}{n^{1/4}},\text{ for all }k \le 201C'\log n\Big)
\notag &\le \prod_{i = 1}^{200} \mathbb{P}\lt(|X_i| \lesssim \frac{\sqrt{\lambda-1}}{n^{1/4}} + \frac{\lambda^{C'\log n}\log^{4} n}{n^{3/4}(\lambda-1)^{1.5}}\rt) \\
\notag &= \prod_{i = 1}^{200} \mathbb{P}\lt(\sqrt{\frac{n(\lambda^2 - 1)}{\lambda^{2C'\log n+2} - \lambda^2}}|X_i| \lesssim \frac{(\lambda-1)n^{1/4}}{\lambda^{C'\log n}} + \frac{\log^{4} n}{n^{1/4}(\lambda-1)}\rt) \\
&\lesssim \Big(\frac{\log^{4} n}{n^{1/4}(\lambda-1)}\Big)^{200} \lesssim n^{-11},
\end{align}
where we invoke the distribution of $X_{i}$ in expression~\eqref{eqn:xi-dist}, and the last inequality results from the assumption that $\lambda - 1 \gtrsim n^{- 1/9}$. 
Therefore, the above inequality guarantees that with probability at least $1 - O(n^{-10})$, 
there exists some $k \lesssim \frac{\log n}{\lambda-1}$ such that 
\begin{align}
\label{eqn:stageI-alpha-nice}
	|\alpha_{k}| \gtrsim \sqrt{\lambda - 1}\, n^{-1/4}.
\end{align}
It thus implies that $\tau_0 \lesssim \frac{\log n}{\lambda-1}$ (see the definition \eqref{eq:defn-tau0-threshold}), as claimed in \eqref{eqn:tau0}.  
In other words, after at most $O(\frac{\log n}{\lambda-1})$ iterations, $|\alpha_t|$ shall surpass the order of $\sqrt{\lambda - 1} \,n^{-1/4}$.
\end{proof}


\paragraph{Stage II: moderate-to-large correlation ($\sqrt{\lambda-1}\,n^{-1/4} \lesssim |\alpha_t| \leq \frac{1}{2} \sqrt{\lambda^2 - 1}$).}
Next, let us look at the time interval after $|\alpha_t|$ surpasses the level of  $\sqrt{\lambda-1}\,n^{-1/4}$ but before it reaches the level of $\frac{1}{2} \sqrt{\lambda^2 - 1}$.  Mathematically, this refers to the interval $(\tau_0, \varsigma]$, where $\tau_0$ and $\varsigma$ are defined in \eqref{eq:defn-tau0-threshold} and \eqref{eq:defn-varsigma}, respectively. 
In fact, we shall start by examining 
\begin{equation}
	t\in \Big( \tau_0, \, \varsigma \wedge \frac{c_5\log n}{\lambda - 1} \Big]
	\label{eq:t-range-stage-II}
\end{equation}
for some constant $c_5>0$; we shall demonstrate that  $ \varsigma \lesssim \frac{\log n}{\lambda - 1}$ shortly.

In view of  Theorem~\ref{thm:main-AMP} and the bounds \eqref{eqn:ABD}, we can write 
\begin{align}
\label{eqn:alpha-crude-1}
	\alpha_{t+1} = \lambda v^{\star \top} \int{\eta}_{t}\lt(\alpha_t \vstar + \frac{1}{\sqrt{n}}x\rt)\varphi_n(\dx) +  \Delta_{\alpha,t},
\end{align}
where the residual term obeys
\begin{equation}
 	|\Delta_{\alpha,t}| \lesssim \sqrt{\frac{t\log n}{n}} + \big|v^{\star\top}\eta_{t}(x_{t})-v^{\star\top}\eta_{t}(v_{t})\big|.
	\label{eq:delta-alpha-general-579}
\end{equation}
%
We first make a claim concerned with a refined recursive relation for $\alpha_{t+1}$:   
\begin{align} 
\label{eq:alpha-global}
|\alpha_{t+1} |
&\geq \frac{\lambda |\alpha_t| }{\sqrt{\alpha^2_t + 1 }} 
	+ o\big( (\lambda-1) |\alpha_t| \big) + O(|\Delta_{\alpha,t}|); 
\end{align}
the proof of this result is postponed to Section~\ref{sec:proof-alpha-global}.
Observe that whenever $\alpha_t < \frac{1}{2}\sqrt{\lambda^2 - 1}$, it holds that 
\begin{align}
\label{eqn:algebra-sim}
	\bigg(\frac{\lambda}{\sqrt{1 + \alpha_t^2}}\bigg)^2 
	\geq \frac{\lambda^2}{\frac{1}{4}\lambda^2 + \frac{3}{4}}
	> 1 + \frac{1}{3}(\lambda - 1), \quad \text{for } \lambda \in (1, 1.2],
\end{align}
which when taken together with expression~\eqref{eq:alpha-global}, implies 
\begin{align}
\label{eqn:middle-step}
|\alpha_{t+1} | &\geq \Bigg(\sqrt{1 + \frac{1}{3}(\lambda-1)} + o(\lambda-1)\Bigg) |\alpha_t| + O(|\Delta_{\alpha,t}|).
\end{align}

In addition, we claim that for every $t \leq \tau'$ where $\tau' \defn \min \{t : \alpha_t \geq (\lambda-1)^{-3/4}n^{-1/4}\}$, it satisfies 
\begin{align}
\label{eqn:delta-alpha}
	|\Delta_{\alpha,t}| \ll (\lambda-1) |\alpha_t|, 
\end{align}
which we shall establish in Section~\ref{sec:proof-claim:eqn:delta-alpha}. 
With the relations~\eqref{eqn:middle-step} and \eqref{eqn:delta-alpha} in place, it obeys $|\alpha_{\tau'+1}| \geq |\alpha_{\tau'}|$.  
Moreover, observe that the bound~\eqref{eq:delta-alpha-general-579} taken together with \eqref{eq:xi-stage1} ensures that 
\begin{align}
\label{eqn:alpha-crude}
 |\Delta_{\alpha,\tau'+1}| \lesssim \sqrt{\frac{(\tau'+1)^3\log n}{n}} \ll (\lambda - 1)|\alpha_{\tau'+1}|.
\end{align}
Invoking this argument recursively, we thus arrive at 
\begin{align*}
	|\alpha_{t+1}| \geq \Bigg(\sqrt{1 + \frac{1}{3}(\lambda-1)} + o(\lambda-1)\Bigg) |\alpha_t|. 
\end{align*}

Now taking the above recursive relation collectively with the assumption $\lambda - 1 \gtrsim n^{- 1/9}$ reveals that $|\alpha_t|$ surpasses $\frac{1}{2}\sqrt{\lambda^2 - 1}$ within at most $O\big(\frac{\log n}{\lambda-1}\big)$ iterations.
Therefore, recalling our definition \eqref{eq:defn-varsigma} of $\varsigma$, 
we can readily conclude that 
\begin{equation}
	\varsigma = O\bigg(\frac{\log n}{\lambda-1}\bigg).
	\label{eq:UB-var-sigma-log}
\end{equation}
%


\subsubsection{A more complete bound for $\xi_t$ (Proof of Claim~\eqref{eq:induction-main})}
\label{sec:more-complete-xit}

We now move on to establish claim~\eqref{eq:induction-main} for any $t$ obeying $\varsigma < t < \frac{c n(\lambda-1)^5}{\log^2 n}$ (recall from \eqref{eq:UB-var-sigma-log} that $\varsigma = O\big(\frac{\log n}{\lambda - 1}\big)$), again via an inductive argument. Along the way, we also need to demonstrate that
\begin{align}
	|\alpha_t| \geq \frac{1}{2} \sqrt{\lambda^2 - 1} 
	\label{eq:condition-alpha-t-lower-bound}
\end{align}
within this stage (namely, once $|\alpha_t|$ exceeds $ \frac{1}{2} \sqrt{\lambda^2 - 1}$, the signal strength will never fall below this level).

To begin with, the claim \eqref{eq:induction-main} for the base case $t=\varsigma$ has already been validated in expression~\eqref{eq:xi-stage1}; 
the condition \eqref{eq:condition-alpha-t-lower-bound} also holds with high probability when $t=\varsigma$. 
Next, assuming that the claim~\eqref{eq:induction-main}) holds up till iteration $t-1$, 
we would like to establish its validity for time $t$.  
Towards this end, inequality~\eqref{eqn:christmas} together with Lemma~\ref{lem:kappa} tells us that  
\begin{align}
\|\xi_{t}\|_{2} & \leq\sqrt{\kappa_{t}+\sqrt{\frac{t\log^{2}n}{n}}}\|\xi_{t-1}\|_{2}+O\Bigg(\sqrt{\frac{t\log n}{n}}+\sqrt{(1+t\ind_{t\le\varsigma})\log n}\,\|\xi_{t-1}\|_{2}^{2}+\sqrt{\frac{t\log n}{n}}\|\xi_{t-1}\|_{2}\Bigg)\nonumber \\
 & \leq\sqrt{1-\frac{1}{15}(\lambda-1)+\sqrt{\frac{t\log^{2}n}{n}}}\|\xi_{t-1}\|_{2}+O\Bigg(\sqrt{\frac{t\log n}{n}}\Bigg)\nonumber \\
 & \qquad+O\Bigg(\sqrt{(1+t\ind_{t\le\varsigma})\log n}\,\left(\sqrt{\frac{t\ind_{t>\varsigma}\log n}{n(\lambda-1)^{2}}}+\sqrt{\frac{\min\{t,\varsigma\}^{3}\log n}{n}}\right)+\sqrt{\frac{t\log n}{n}}\Bigg)\|\xi_{t-1}\|_{2}\nonumber \\
 & \le\lt(1-\frac{1}{15}(\lambda-1)+O\bigg(\sqrt{\frac{\big(\frac{t}{(\lambda-1)^{2}}+\varsigma^{3}\big)\log^{2}n}{n}}\bigg)\rt)\|\xi_{t-1}\|_{2}+O\Big(\sqrt{\frac{t\log n}{n}}\Big), \label{eqn:error-final-2}
\end{align}
where the second line comes from \eqref{eq:kappa-claim-main} and the induction hypothesis \eqref{eq:induction-main} for $t-1$. 
In addition,  the validity of \eqref{eq:condition-alpha-t-lower-bound} for $\alpha_{t+1}$ in the $(t+1)$-th iteration can be justified as well, 
which we shall detail in Section~\ref{sec:proof:eq:condition-alpha-t-lower-bound}.

With the above recursive relation in mind, recognizing $\sqrt{\frac{(t/(\lambda-1)^{2}+\varsigma^{3})\log^{2}n}{n}} \leq 2c(\lambda - 1)$ for some constant $c > 0$ small enough, 
we can readily derive 
\begin{align}
\notag \|\xi_{t}\|_{2} &\le \lt(1-\frac{1}{20}(\lambda-1)\rt)\|\xi_{t-1}\|_{2}+O\Big(\sqrt{\frac{t\log n}{n}}\Big) \\
\notag &\le \lt(1-\frac{1}{20}(\lambda-1)\rt)^{t - \varsigma}\|\xi_{\varsigma}\|_{2}
+ \sum_{j=0}^{t-\varsigma-1} \lt(1-\frac{1}{20}(\lambda-1)\rt)^{j} O\Big(\sqrt{\frac{(t-j)\log n}{n}}\Big) \\
&\lesssim \sqrt{\frac{\varsigma^3\log n}{n}} +\sqrt{\frac{t\log n}{n(\lambda - 1)^2}}
\end{align}
for all $\varsigma \le t \leq \frac{c n(\lambda - 1)^5}{\log^2 n}$. 
Combining this with the bound \eqref{eq:xi-stage1} (for $t\leq \varsigma$) immediately establishes Claim~\eqref{eq:induction-main} for all $t \leq \frac{c n(\lambda - 1)^5}{\log^2 n}$.

\subsubsection{Analysis for approximate state evolution (Proof of Property~\eqref{eqn:cadence})}
\label{sec:SE}

Once the signal strength $\alpha_t$ reaches the order of $\sqrt{\lambda^2 - 1}$, AMP enters the stage of local refinement. According to \eqref{eq:induction-main} and \eqref{claim:alpha}, for any $t \leq \frac{cn(\lambda-1)^5}{\log^2 n}$, the AMP iterate $x_t$ admits the decomposition~\eqref{eqn:xt-decomposition} with the error term bounded by 
\begin{align} 
	\|\xi_{t}\|_2 \lesssim \sqrt{\frac{(t + \frac{\log^3 n}{\lambda-1})\log n}{n(\lambda-1)^2}}.
\end{align}
In the meantime, to describe how $\alpha_{t}$ evolves, we bound $\Delta_{\alpha,t}$ based on the relation~\eqref{eqn:delta-alpha-general} as follows: 
\begin{align}
\label{eqn:stage3-error}
	\vert\Delta_{\alpha,t}\vert&
	\lesssim B_t + \ltwo{\xi_{t-1}}\lesssim \sqrt{\frac{(t+ \frac{\log^3n}{\lambda-1})\log n}{n(\lambda-1)^2}},
\end{align}
where the last inequality arises from \eqref{eqn:ABD}. 
Combining this with relation~\eqref{eqn:alpha-t-genearl} leads to 
\begin{align*}
	\alpha_{t+1} 
	&= \lambda v^{\star \top} \int{\eta}_{t}\lt(\alpha_t \vstar + \frac{1}{\sqrt{n}}x\rt)\varphi_n(\dx) 
	+ O\lt(\sqrt{\frac{(t+ \frac{\log^3n}{\lambda-1})\log n}{n(\lambda-1)^2}}\rt).
\end{align*}

Next, we shall characterize the distance between $\alpha_{t+1}$ and its asymptotic counterpart to further understand the evolution of $\alpha_{t+1}$. 
More specifically, recall that the asymptotic state evolution is defined as
\begin{align}
	\alpha_{t+1}^{\star} &= \lambda \lt[\int \tanh\big(\alpha_t^{\star}\lt(\alpha_t^{\star} + x\rt)\big)\varphi(\dx)\rt]^{1/2},
	\label{eq:analysis-SE-star}
\end{align}
assuming we start from $\alpha_{\varsigma}^{\star} = |\alpha_{\varsigma}|$
for some $\varsigma = O\big(\frac{\log n}{\lambda-1}\big)$. 
We aim to control the difference between $\alpha_{t+1}$ and $\alpha_{t+1}^{\star}$. 
To simplify the presentation, we assume without loss of generality that $\alpha_t >0$, and employ the notation $$\tau_t \defn \big(\alpha_t^{\star}\big)^2.$$

To begin with, the same analysis as in the proof of claim~\eqref{eq:alpha-global} (with different error bound~\eqref{eqn:stage3-error} here) gives 
\begin{align}
\label{eqn:golden}
\notag \alpha_{t+1} 
	&= \lambda \lt[\int \tanh\lt(\alpha_t^2 + \alpha_tx\rt)\varphi(\dx)\rt]^{1/2} + O\lt(\lt|\frac{\pi_t^2}{\alpha_t^2n} - 1\rt|\alpha_t^3 + \sqrt{\frac{(t+ \frac{\log^3n}{\lambda-1})\log n}{n(\lambda-1)^2}}\rt) \\
	&= \lambda \lt[\int \tanh\lt(\alpha_t^2 + \alpha_tx\rt)\varphi(\dx)\rt]^{1/2} + O\lt(\sqrt{\frac{(t+ \frac{\log^3n}{\lambda-1})\log n}{n(\lambda-1)^2}}\rt). 
\end{align}
Here, the last line follows from inequality~\eqref{eqn:pi_t} which indicates 
\begin{align*}
\bigg(\frac{\pi_t}{\sqrt{n}}\bigg)^2 &
= \alpha_t^2 + O\bigg(\|\xi_{t-1}\|_2 + \sqrt{\frac{t\log n}{n}}\bigg)
= \alpha_t^2 + O\lt(\sqrt{\frac{(t+ \frac{\log^3n}{\lambda-1})\log n}{n(\lambda-1)^2}}\rt).
\end{align*}
It then follows from relations~\eqref{eqn:golden} and \eqref{eq:analysis-SE-star} that 
\begin{align}
\label{eqn:recursion-se}
%
\frac{\alpha^2_{t+1} - \tau_{t+1}}{\tau_{t+1}} 
&= 
\frac{\int \lt[\tanh\lt(\alpha_t^2 + \alpha_t x\rt) - \tanh\lt(\tau_t + \sqrt{\tau_t}x\rt)\rt]\varphi(\dx)}{\int \tanh\lt(\tau_t + \sqrt{\tau_t}x\rt)\varphi(\dx)} 
+ O\lt(\sqrt{\frac{(t+ \frac{\log^3n}{\lambda-1})\log n}{n(\lambda-1)^3}}\rt).
\end{align}
Here, we remind the readers that  (see also \eqref{eqn:ut-num} and \citet[Appendix~B.2]{deshpande2017asymptotic})
\begin{align*}
	\int \tanh\lt(\alpha^2 + \alpha x\rt)\varphi(\dx) = \int \tanh^2\lt(\alpha^2 + \alpha x\rt)\varphi(\dx) \asymp \alpha^2,
	\qquad \text{for }\alpha\in (0,\lambda],
\end{align*}
where the last inequality results from relation~\eqref{eqn:ut-num}.
The recursive formula~\eqref{eqn:recursion-se} quantifies how the difference between $\alpha_{t}$ and $\alpha^\star_{t}$ changes over time, which plays a key role in our following analysis.

In order to better understand the above recursion,  let us define --- for every $\tau\in [0,\lambda^2]$ --- that 
\begin{align*}
h(\tau) := \int \tanh\lt(\tau + \sqrt{\tau} x\rt)\varphi(\dx).
\end{align*}
Armed with this function, one can write
\begin{align}
\label{eqn:alpha-in-h}
	\alpha_{t+1}^{\star 2} = \lambda^2 h(\alpha_{t}^{\star 2}).
\end{align}
Also, direct calculations yield
\begin{align*}
h^{\prime}(\tau) := \int \Big(1 + \frac{x}{2\sqrt{\tau}}\Big)\lt(1 - \tanh^2\lt(\tau + \sqrt{\tau} x\rt)\rt)\varphi(\dx) \in (0,1),
\end{align*}
where its range follows from display (263) in~\citet[Section D.3.3]{li2022non}.
We make note of a few direct consequences of the above results. 

\begin{itemize}
\item Recognizing that $h^{\prime}(\tau) > 0$, one has $\alpha_{t+1}^{\star} > \alpha_{t}^{\star} \gtrsim \sqrt{\lambda^2-1}$ for $t \geq \varsigma.$

\item In view of display (264) in~\citet[Section D.3.3]{li2022non}, we have $0 \leq \lambda^2 h^{\prime}(\tau) \leq 1 - (\lambda - 1)$. 
	If we define $\alpha^{\star }$ to be the limiting point of \eqref{eq:analysis-SE-star} (as $t\rightarrow \infty)$),  we can then see that  
\begin{align}
\label{eqn:se-asym}
	|\alpha_{t+1}^{\star 2} - \alpha^{\star 2}| \leq \big(1- (\lambda - 1)\big) \cdot |\alpha_{t}^{\star 2} - \alpha^{\star 2}|, 
	\qquad\text{for } t \ge \varsigma.
\end{align}
In other words, the asymptotic state evolution parameter $\alpha_{t}^{\star 2}$ converges exponentially to some fixed point $\alpha^{\star 2}.$

\item In light of the above notation, we can also write  
\begin{align}
\notag 
\frac{|\alpha^2_{t+1} - \tau_{t+1}|}{\tau_{t+1}} 
&= \frac{|h(\alpha^2_{t}) - h(\tau_{t})|}{h(\tau_{t})} + O\lt(\sqrt{\frac{(t+ \frac{\log^3n}{\lambda-1})\log n}{n(\lambda-1)^3}}\rt)\\
&= \frac{h^{\prime}(\tau)}{h(\tau_t)/\tau_{t}}\cdot\frac{|\alpha^2_{t} - \tau_{t}|}{\tau_{t}} + O\lt(\sqrt{\frac{(t+ \frac{\log^3n}{\lambda-1})\log n}{n(\lambda-1)^3}}\rt)
	\label{eq:recursion-ratio-alpha}
\end{align}
for some $\tau$ satisfying $\min\{\tau_t, \alpha_t^2\} \leq \tau \leq \max\{\tau_t, \alpha_t^2\}$. 
\end{itemize}

We first prove that $\alpha_t^2 = (1+o(1)) \tau_{t}$. 
By definition,   $\alpha_{\varsigma}^{\star} = \alpha_{\varsigma} \gtrsim \sqrt{\lambda^{2}-1}$, and hence this claim holds trivially for $t = \varsigma$. 
Next, assuming the validity of the inductive assumption $\alpha_t^2 = (1+o(1)) \tau_{t}$, 
we would like to prove it for  the $(t+1)$-th step. 
Towards this end, we first claim that there exists some universal constant $c > 0$ small enough such that 
\begin{align}
\label{eqn:derivative}
	\frac{h^{\prime}(\tau)}{h(\tau_t)/\tau_{t}} \leq 1 - c(\lambda - 1),
\end{align}
whose proof is postponed to Section~\ref{sec:proof-eqn:derivative}. 
This further allows us to derive
\begin{align}\label{equ:alpha-se-recursion-claim}
\frac{|\alpha_t^2-\tau_t|}{\tau_t}\lesssim \sqrt{\frac{(t+ \frac{\log^3n}{\lambda-1})\log n}{n(\lambda-1)^5}}, 
\end{align}
which we shall demonstrate via an inductive argument. 
Given that \eqref{equ:alpha-se-recursion-claim} holds trivially when $t=\varsigma$,  we intend to establish~\eqref{equ:alpha-se-recursion-claim} 
for  the $(t+1)$-th iteration, assuming that it holds for all $s\leq t$.
To do so, we combine \eqref{eqn:derivative} and \eqref{eq:recursion-ratio-alpha} to show that   
\begin{align}
\label{eqn:alpha-se-recursion}
\notag \frac{|\alpha^2_{t+1} - \tau_{t+1}|}{\tau_{t+1}}	
&\le \big(1 - c(\lambda-1)\big) \cdot \frac{|\alpha^2_{t} - \tau_{t}|}{\tau_{t}} + O\lt(\sqrt{\frac{(t+ \frac{\log^3n}{\lambda-1})\log n}{n(\lambda-1)^3}}\rt) \\
&\notag=\big(1-c(\lambda-1)\big)^{t+1-\varsigma}\cdot\frac{|\alpha_\varsigma^2-\tau_\varsigma|}{\tau_\varsigma} + O\lt(\sum_{k=\varsigma}^{t-1}\big(1-c(\lambda-1)\big)^{t-k}\sqrt{\frac{(k+ \frac{\log^3n}{\lambda-1})\log n}{n(\lambda-1)^3}}\rt)\\
&\notag\leq \sum_{k=\varsigma}^{t-1}\big(1-c(\lambda-1)\big)^{t-k}\cdot O\lt(\sqrt{\frac{(t+ \frac{\log^3n}{\lambda-1})\log n}{n(\lambda-1)^3}}\rt)\\
&\lesssim \sqrt{\frac{(t+ \frac{\log^3n}{\lambda-1})\log n}{n(\lambda-1)^5}},
\end{align}
where the second line is obtained by applying~\eqref{equ:alpha-se-recursion-claim} recursively. Hence, it also leads to $\alpha_{t+1}^2 = (1+o(1)) \tau_{t+1}$, which concludes the inductive assumption for the $(t+1)$-th iteration. 
Putting the above results together with expression~\eqref{eqn:se-asym} also gives
\begin{align}
\label{eqn:se-non-asym}
	\notag |\alpha^2_{t+1} - \alpha^{\star 2}|
	&=  \big(1- (\lambda - 1)\big)^{t - \varsigma} \cdot |\alpha_{\varsigma}^{\star 2} - \alpha^{\star 2}| + \alpha^{\star 2}_{t+1} O\lt( \sqrt{\frac{(t+ \frac{\log^3n}{\lambda-1})\log n}{n(\lambda-1)^5}} \rt)\\
	&\lesssim  \big(1 - (\lambda - 1)\big)^{t -\varsigma}
	  	+
	 \sqrt{\frac{(t+ \frac{\log^3n}{\lambda-1})\log n}{n(\lambda-1)^5}} .
\end{align}

%

\section{Discussions}

In this paper, we have pinned down the finite-sample convergence behavior of AMP when initialized randomly,  focusing on the prototypical $\mathbb{Z}_{2}$ synchronization problem. 
This algorithm has been shown to enjoy fast global convergence, 
as it takes no more than $O(\frac{\log n}{\lambda - 1})$ iterations to arrive at a point whose risk is
 $O(\sqrt{\frac{\log^4 n}{n(\lambda - 1)^6}})$ close to Bayes-optimal. 
To the best of our knowledge, our theory offers the first rigorous evidence supporting the effectiveness of randomly initialized AMP in low-rank matrix estimation. 
While the present paper concentrates on a specific choice of denoising functions tailored to $\mathbb{Z}_2$ synchronization, we expect our analysis framework to be generalizable to a broader family of separable and Lipschitz-continuous denoising functions.

Moving forward, there is no shortage of research directions worth exploring. One natural extension is concerned with other structural prior about $\vstar$; for instance, it would be interesting to see how randomly initialized AMP performs when $\vstar$ is known to satisfy general cone constraints (see e.g.~\cite{bandeira2019computational,wei2019geometry}). 
Another direction of interest is to go beyond the spiked Wigner model.  
A recent work along this line \citep{wu2022lower}  studied the role of random initialization for power iteration in the problem of tensor decomposition, which leverages upon the AMP-type analysis for analyzing tensor power methods. 
Can we further extend these to understand (randomly initialized) AMP towards solving more challenging problems like low-rank matrix completion and tensor completion? 
Moreover, note that the current theory characterizes the behavior of AMP up to $O(\frac{n}{\textsf{poly}(\log n)})$ number of steps.  
While it is already useful in practice, it would be of interest --- from a technical perspective --- to understand whether AMP eventually converges to a fixed point as the iteration number further increases. 
Finally, while AMP serves as a versatile machinery for understanding various statistical procedures in high dimensions, there are several alternative analysis frameworks like the convex Gaussian minmax theorem (CGMT) \citep{celentano2020lasso,miolane2021distribution,thrampoulidis2018precise} and the leave-one-out analysis \citep{el2018impact,chen2020noisy,abbe2020entrywise} that also prove effective and enjoy their own benifits. Is there any effective way to combine them so as to exploit all of their advantages at once? 
We leave these questions for future investigation.

\section*{Acknowledgment}

This work was partially supported by NSF grants DMS 2147546/2015447 and the NSF CAREER award DMS-2143215. 

\appendix

\section{Proof of auxiliary lemmas and claims}

\subsection{Proof of Lemma~\ref{lem:parameters}}
\label{sec:proof-parameters}


\paragraph{Proof of property~\eqref{eqn:pi_t}.}
Recall that $\pi_t \defn \sqrt{n(\Vert x_t\Vert_2^2-1)\vee 1}$.
To show property~\eqref{eqn:pi_t}, the first step is to calculate $\ltwo{x_t}$. 
Notice that for independent Gaussian vectors $\phi_{k} \sim \mathcal{N}\big(0,\frac{1}{n}I_n\big)$, 
one has $${v^{\star}}^\top\lt[\phi_1,\ldots,\phi_{t-1}\rt]\sim\mathcal{N}\Big(0,\frac{1}{n}I_{t-1}\Big),$$
given that $\|v^{\star}\|_2=1$. 
Therefore, it is easily seen that
\begin{align}
\label{equ:v-t}
\Big|\Big\langle v^{\star},\sum_{k=1}^{t-1}\beta_{t-1}^{k}\phi_k\Big\rangle\Big|
=
\lt|\lt\langle {v^{\star}}^\top\lt[\phi_1,\ldots,\phi_{t-1}\rt],\lt[\beta_{t-1}^1,\ldots,\beta_{t-1}^{t-1}\rt]\rt\rangle\rt|
&\leq \lt\Vert {v^{\star}}^\top\lt[\phi_1,\ldots,\phi_{t-1}\rt]\rt\Vert_2\Vert\beta_{t-1}\Vert_2\notag\\
&\lesssim\sqrt{\frac{t\log n}{n}}
\end{align}
with probability at least $1 - O(n^{-11}).$
Combining inequalities~\eqref{equ:v-t} and \eqref{equ:beta-phi}, we reach 
\begin{align}\label{eqn:ltwo-v-square}
\notag 
\Vert v_t\Vert_2^2 &= \lt\Vert\alpha_tv^{\star} + \sum_{k=1}^{t-1}\beta_{t-1}^{k}\phi_k\rt\Vert_2^2=\Vert\alpha_tv^{\star}\Vert_2^2 + \lt\Vert\sum_{k=1}^{t-1}\beta_{t-1}^{k}\phi_k\rt\Vert_2^2 + 2\lt\langle\alpha_tv^{\star},\sum_{k=1}^{t-1}\beta_{t-1}^{k}\phi_k\rt\rangle\\
&=\alpha_t^2 + 1 + O\lt(\sqrt{\frac{t\log n}{n}}\rt) \asymp 1,
\end{align}
which in turn leads to
\begin{align}
\|x_t\|_2^2 = \big(\|v_t\|_2 + O(\|\xi_{t-1}\|_2)\big)^2 = \alpha_t^2 + 1 + O\bigg(\|\xi_{t-1}\|_2 + \sqrt{\frac{t\log n}{n}}\bigg).
\end{align}

Based on the above properties, 
we can also derive that
\begin{align}
	\sqrt{n(\|x_t\|_2^2 - 1)} &=\sqrt{n\left(\alpha_t^2+O\left(\Vert\xi_{t-1}\Vert_2+\sqrt{\frac{t\log n}{n}}\right)\right)}
\leq \sqrt{n}\,|\alpha_t| + O\lt(\sqrt{n \bigg(\|\xi_{t-1}\|_2 + \sqrt{\frac{t\log n}{n}}\bigg)}\rt) \notag\\
&=\sqrt{n}\lt( |\alpha_t| +O\bigg(\bigg(\|\xi_{t-1}\|_2 + \sqrt{\frac{t\log n}{n}}\bigg)^{1/2}\bigg)\rt). 
\end{align}
where the first inequality based on the fact that $\sqrt{a+b}\leq\sqrt{a}+\sqrt{b}$ for $a,b\geq 0$. 
In particular, if $\alpha_{t}^{2} \lesssim \Vert\xi_{t-1}\Vert_{2}+\sqrt{\frac{t\log n}{n}}$, 
then the basic relation $\sqrt{a+b}=\sqrt{a}+O(\sqrt{b})$ ($0\leq a\lesssim b$) enables us to replace ``$\leq$'' in \eqref{eq:second-bound-xt-minus-1} with ``$=$'' to obtain
\begin{align}
	\sqrt{n(\|x_t\|_2^2 - 1)} 
	=\sqrt{n}\lt(\alpha_t+O\bigg(\bigg(\|\xi_{t-1}\|_2 + \sqrt{\frac{t\log n}{n}}\bigg)^{1/2}\bigg)\rt). 
	\label{eq:second-bound-xt-minus-1}
\end{align}
In contrast, if $\alpha_{t}^{2} \gtrsim \Vert\xi_{t-1}\Vert_{2}+\sqrt{\frac{t\log n}{n}}$, 
then the basic relation $\sqrt{a+b}=\sqrt{a}+O(\frac{b}{\sqrt{a}})$ ($0\leq b\lesssim a$)
yields 
\begin{align}
\label{eqn:msk}
\notag \sqrt{n(\|x_t\|_2^2 - 1)} =\sqrt{n\left(\alpha_t^2+O\left(\Vert\xi_{t-1}\Vert_2+\sqrt{\frac{t\log n}{n}}\right)\right)}
&=\sqrt{n}|\alpha_t| \sqrt{1+\frac{1}{\alpha_t^2}O\left(\Vert\xi_{t-1}\Vert_2+\sqrt{\frac{t\log n}{n}}\right)}\\
 &=\sqrt{n}\left(|\alpha_t| + \frac{1}{|\alpha_t|}O\left(\|\xi_{t-1}\|_2 + \sqrt{\frac{t\log n}{n}}\right)\right).
\end{align}
The preceding two bounds taken collectively demonstrate that
\[
\sqrt{n(\|x_{t}\|_{2}^{2}-1)}=\sqrt{n}|\alpha_t|+\sqrt{n}\left\{ \frac{1}{|\alpha_t|}O\left(\|\xi_{t-1}\|_{2}+\sqrt{\frac{t\log n}{n}}\right)\,\wedge\,O\bigg(\bigg(\|\xi_{t-1}\|_{2}+\sqrt{\frac{t\log n}{n}}\bigg)^{1/2}\bigg)\right\} .
\]


The above bound also leads to the desired form~\eqref{eqn:pi_t} for $\pi_t$ when $\sqrt{n(\Vert x_t\Vert_2^2-1)} \geq 1$. 
Consequently, it remains to examine the case where $\sqrt{n(\Vert x_t\Vert_2^2-1)} < 1$, which clearly can only happen if
\begin{align*}
	c_{10} |\alpha_t| \leq  \bigg(\|\xi_{t-1}\|_2 + \sqrt{\frac{t\log n}{n}}\bigg)^{1/2}
\end{align*}
for some sufficiently small constant $c_{10}>0$.   
But in this situation, we still have
\begin{align*}
\pi_{t}=1 & \lesssim\sqrt{n}\bigg(\|\xi_{t-1}\|_{2}+\sqrt{\frac{t\log n}{n}}\bigg)^{1/2}\asymp\sqrt{n}\lt(|\alpha_t|+O\bigg(\bigg(\|\xi_{t-1}\|_{2}+\sqrt{\frac{t\log n}{n}}\bigg)^{1/2}\bigg)\rt)\\
 & \asymp\sqrt{n}|\alpha_t|+\sqrt{n}\left\{ O\bigg(\bigg(\|\xi_{t-1}\|_{2}+\sqrt{\frac{t\log n}{n}}\bigg)^{1/2}\bigg)\,\wedge\,\frac{1}{|\alpha_t|}O\left(\|\xi_{t-1}\|_{2}+\sqrt{\frac{t\log n}{n}}\right)\right\} ,
\end{align*} 
and hence the claimed bound is still valid.

\paragraph{Proof of property~\eqref{eqn:gamma_t}.}
First recall the definition $\gamma_t^{-2} \defn \Vert\tanh(\pi_t x_t)\Vert^2_2$. 
Towards establishing property~\eqref{eqn:gamma_t}, consider the following difference 
\begin{align}
&\lt\vert\lt\Vert\tanh(\pi_t x_{t})\rt\Vert^2_2-\int\Big\Vert\tanh\lt(\frac{\pi_t}{\sqrt{n}}(\alpha_t + x)\rt)\Big\Vert^2_2\varphi_n(\dx)\rt\vert\notag \\
&\qquad \leq \lt\vert\lt\Vert\tanh(\pi_t x_{t})\rt\Vert^2_2-\lt\Vert\tanh(\pi_t v_t)\rt\Vert^2_2\rt\vert 
+ \lt\vert\lt\Vert\tanh(\pi_t v_{t})\rt\Vert^2_2-\int\Big\Vert\tanh\lt(\frac{\pi_t}{\sqrt{n}}(\alpha_t + x)\rt)\Big\Vert^2_2\varphi_n(\dx)\rt\vert\notag
\end{align}
where  $\varphi_n\sim\mathcal{N}(0,I_n)$.
To bound the right-hand side above, 
note that the Lipschitz property of $\tanh$ gives 
\begin{align*}
  \lt\Vert\tanh(\pi_t x_{t}) - \tanh(\pi_t v_{t})\rt\Vert_2 
	\leq \lt\Vert\tanh(\pi_t (x_{t}-v_t)) \rt\Vert_2 
	\le \pi_t \Vert\xi_{t-1}\Vert_2. 
\end{align*}
This in turn yields 
\begin{align*}
\lt\vert\lt\Vert\tanh(\pi_{t}x_{t})\rt\Vert_{2}^{2}-\lt\Vert\tanh(\pi_{t}v_{t})\rt\Vert_{2}^{2}\rt\vert & =\big|\lt\Vert\tanh(\pi_{t}x_{t})\rt\Vert_{2}+\lt\Vert\tanh(\pi_{t}v_{t})\rt\Vert_{2}\big|\cdot\big|\lt\Vert\tanh(\pi_{t}x_{t})\rt\Vert_{2}-\lt\Vert\tanh(\pi_{t}v_{t})\rt\Vert_{2}\big|\\
 & \le\big(2\lt\Vert\tanh(\pi_{t}v_{t})\rt\Vert_{2}+\pi_{t}\Vert\xi_{t-1}\Vert_{2}\big)\pi_{t}\Vert\xi_{t-1}\Vert_{2}\\
 & \lesssim\big(\pi_{t}\lt\Vert v_{t}\rt\Vert_{2}+\pi_{t}\Vert\xi_{t-1}\Vert_{2}\big)\pi_{t}\Vert\xi_{t-1}\Vert_{2}\asymp\pi_{t}^{2}\Vert\xi_{t-1}\Vert_{2},
\end{align*}
where the last line makes use of the assumption $\ltwo{\xi_{t-1}}\lesssim 1$ and the equation~\eqref{eqn:ltwo-v-square}. 
Thus, we arrive at 
\begin{align}
\label{equ:tanh-pi-gau}
&\lt\vert\lt\Vert\tanh(\pi_t x_{t})\rt\Vert^2_2-\int\Big\Vert\tanh\lt(\frac{\pi_t}{\sqrt{n}}(\alpha_t + x)\rt)\Big\Vert^2_2\varphi_n(\dx)\rt\vert\notag \\
&\qquad \lesssim \pi_t^2\Vert \xi_{t-1}\Vert_2 + \lt\vert\lt\Vert\tanh(\pi_t v_{t})\rt\Vert^2_2-\int\Big\Vert\tanh\lt(\frac{\pi_t}{\sqrt{n}}(\alpha_t + x)\rt)\Big\Vert^2_2\varphi_n(\dx)\rt\vert.
\end{align}

Next we develop a bound on the second term of~\eqref{equ:tanh-pi-gau}, which turns out to be a consequence of the uniform concentration property. 
Specifically, let us consider functions of the following form 
\begin{align*}
  f_{\theta}(\Phi) \defn \lt\Vert\tanh(\pi v)\rt\Vert^2_2-\int\lt\Vert\tanh\lt(\frac{\pi}{\sqrt{n}}(\alpha + x)\rt)\rt\Vert^2_2\varphi_n(\dx),\qquad\text{where }v = \alpha v^{\star} + \sum_{k=1}^{t-1}\beta^k\phi_k; 
\end{align*}
here, we define
\begin{align*}
  \Phi = \sqrt{n}\lt[\phi_1,\phi_2,\ldots,\phi_{t-1}\rt],\quad \beta = \lt[\beta^1,\beta^2,\ldots,\beta^{t-1}\rt],\quad\theta=\lt[\alpha,\beta,\pi\rt]\in\mathbb{R}^{t+1}.
\end{align*}
Then, it suffices to bound $f_{\theta}(\Phi)$ uniformly over all $\theta$ in the following set
\begin{align}
\Theta : = \left\{\theta = (\alpha,\beta,\pi)\;|\;\Vert\beta\Vert_2 =1, \alpha\lesssim \sqrt{\lambda^2-1}, \pi\lesssim\sqrt{n}\right\}.
\end{align}
In order to achieve this, first consider the derivative of $f$ with respect to $\Phi$, which by direct calculations satisfies 
\begin{align}
\notag \Vert \nabla_{\Phi} f_{\theta}(\Phi)\Vert_2 
&\leq \frac{2\pi\Vert\beta\Vert_2}{\sqrt{n}}\Vert\tanh(\pi v)\circ \tanh'(\pi v)\Vert_2\\
&\leq\frac{2\pi\Vert\beta\Vert_2}{\sqrt{n}}\Vert\tanh(\pi v)\Vert_2 \lesssim \frac{\pi^2}{\sqrt{n}}\|v\|_2 \lesssim \frac{\pi^2}{\sqrt{n}},
\end{align}
where we note that $\|v\|_2 \le \alpha + \frac{1}{\sqrt{n}}\|\Phi\| \lesssim 1$.
Additionally, since function $f_{\theta}(\Phi)$ is uniformly bounded by $2$, if we take $\delta=n^{-300}$ for the set $\mathcal{E}$ defined in Lemma~\ref{lem:concentration}, it satisfies (we refer the readers to  \citet[Section D.1.1]{li2022non} for the proof of this property)
\begin{align*}
  \mathbb{E}\lt[|f(\Phi)-f(\mathcal{P}_{\mathcal{E}}(\Phi))|\rt]\lesssim n^{-100}.
\end{align*}
where 
$\mathcal{P}_{\mathcal{E}} (\cdot)$ denotes the Euclidean projection onto the set $\mathcal{E}$. 
Combining the above inequality with the following properties of function $f_{\theta}(\Phi)$,
\begin{enumerate}
  \item $\Vert\nabla_{\theta} f_{\theta}(\Phi)\Vert_2\lesssim n ^{100}$ for any $\Phi\in\mathcal{E}$,
  \item For any fixed $\theta$, one has $\mathbb{E}[f_{\theta}(\Phi)]=0$,
\end{enumerate}
we can apply \citet[Corollary 3]{li2022non} to reach 
\begin{align}\label{equ:bound-f-theta}
\sup_{\theta\in\Theta}\Big\vert\frac{1}{\pi^2} f_{\theta}(\Phi)\Big\vert \lesssim \sqrt{\frac{t\log n}{n}}.
\end{align}

Taking everything collectively, we conclude that
\begin{align}
\label{equ:bound-tanh-diff}
\lt\vert\lt\Vert\tanh(\pi_t x_{t})\rt\Vert^2_2-\int\lt\Vert\tanh\lt(\frac{\pi_t}{\sqrt{n}}(\alpha_t + x)\rt)\rt\Vert^2_2\varphi_n(\dx)\rt\vert 
\
\lesssim \pi_t^2\bigg(\Vert \xi_{t-1}\Vert_2 + \sqrt{\frac{t\log n}{n}}\bigg),
\end{align}
which leads to the property~\eqref{eqn:gamma_t} by recognizing that 
\begin{align*}
  \int\lt\Vert\tanh\lt(\frac{\pi_t}{\sqrt{n}}(\alpha_t + x)\rt)\rt\Vert^2_2\varphi_n(\dx)=n\int\tanh^2\big(\frac{\pi_t}{\sqrt{n}}(\alpha_t+x)\big)\varphi(\dx).
\end{align*}

\paragraph{Proof of property~\eqref{eqn:tanh_pi_t}.}
We are only left with calculating the value of $\int\tanh^2\big(\frac{\pi_t}{\sqrt{n}}(\alpha_t+x)\big)\varphi(\dx)$ which shall be done as follows.
First, by observing that $\tanh(0)=\tanh^{\second}(0)=0$, $\tanh^{\prime}(0)=1$, $|\tanh^{\third}(x)|\leq 4$ for $x\in\mathbb{R}$, we find  
\begin{align*}
\vert\tanh(x)-x\vert \leq \frac{2}{3}|x|^3,
\end{align*}
as a consequence of the mean value theorem. 
Combining this relation with the fact that $|\tanh(x)-x|\leq |x|$, we can further obtain
\begin{align*}
\vert\tanh(x)-x\vert \leq |x|\wedge|x|^3.
\end{align*}
As a result, we see that
\begin{align*}
  \vert\tanh(x)+x\vert=\vert 2x+\tanh(x)-x\vert = 2|x| + O(|x|\wedge |x|^3)\lesssim |x|,
\end{align*}
which in turn leads to 
\begin{align*}
\tanh^2(x)=x^2+O(|\tanh(x)-x||\tanh(x)+x|)=x^2 + O(x^4).
\end{align*}

Now we are ready to compute the value of $\int\tanh^2\big(\frac{\pi_t}{\sqrt{n}}(\alpha_t+x)\big)\varphi(\dx)$.
In view of the expressions obtained above, it follows that 
\begin{align}
\label{eqn:norm_t}
\notag \int\tanh^2\lt(\frac{\pi_t}{\sqrt{n}}(\alpha_t+x)\rt)\varphi(\dx)&=\int\lt(\frac{\pi_t}{\sqrt{n}}(\alpha_t+x)\rt)^2\varphi(\dx)+O\lt(\int\lt(\frac{\pi_t}{\sqrt{n}}
(\alpha_t+x)\rt)^4\varphi(\dx)\rt)\\
&=\frac{\pi_t^2}{n}(\alpha_t^2+1)+O\lt(\frac{\pi_t^4}{n^{2}}\rt).
\end{align}
We thus complete the proof of the advertised result. 

\subsection{Proof of Lemma~\ref{lem:eta_property}}
\label{sec:proof-eta_property}

\paragraph{Proof of property~\eqref{eq:eta_derivative}.}
The property~\eqref{eq:eta_derivative} is concerned with the magnitudes of $\eta_{t}$ and its derivatives. In view of the definition of $\eta_{t}$, we proceed to bounding the parameters $\pi_{t}$ and $\gamma_{t}$ separately. Towards this,  recall from Lemma~\ref{lem:concentration} that: with probability at least $1 - O(n^{-11})$,  
%
\begin{align}\label{equ:beta-phi}
	\left\Vert\sum_{k=1}^{t-1}\beta_{t-1}^{k}\phi_k\right\Vert_2=1+O\left(\sqrt{\frac{t\log n}{n}}\right)
\end{align}
holds simultaneously for all $\beta_{t-1}=[\beta_{t-1}^k]_{1\leq k < t}\in \mathcal{S}^{t-2}$. 
It then follows from this result and the definition \eqref{defn:v-t-thm1} of $v_t$ that
\begin{align}
\label{eqn:v-norm}
\ltwo{v_t} &= \left\Vert\alpha_{t} v^{\star{}} + \sum_{k = 1}^{t-1} \beta_{t-1}^k\phi_k\right\Vert_2
\leq\left\Vert \alpha_tv^{\star}\right\Vert_2+ \left\Vert \sum_{k = 1}^{t-1} \beta_{t-1}^k\phi_k\right\Vert_2
	= |\alpha_t| +1+O\left(\sqrt{\frac{t\log n}{n}}\right)\lesssim 1,
\end{align}
given that $t \lesssim \frac{n}{\log n}$ and
\begin{align}
	|\alpha_t| = \big| \lambda {v^{\star}}^{\top}\eta_{t-1}(x_{t-1}) \big| \leq \lambda \lesssim 1. 
	\label{eq:alphat-UB-135}
\end{align}
These properties together with the assumption on $\xi_t$ enable us to control the $\ell_2$ norm of $x_t$ as follows: 
\begin{align}
\Vert x_{t}\Vert_{2}=\left\Vert v_{t}+\xi_{t-1}\right\Vert _{2} & \leq\left\Vert v_{t}\right\Vert _{2}+\left\Vert \xi_{t-1}\right\Vert _{2}\lesssim1.\label{eq:xt-upper-bound-1}
\end{align}
%
%

Given that $\Vert x_t\Vert_2\lesssim 1$, the value of $\pi_t$ can be controlled as 
\begin{align}
 \label{eqn:pi_t-2}
 \pi_t &\defn \sqrt{n(\Vert x_t\Vert_2^2-1)} \vee 1 \lesssim \sqrt{n}.
\end{align}
Additionally, by observing that $\tanh(0) = 0$, $\tanh^{\prime}(x) = 1 - \tanh^2(x) \in [0,1]$ and $\vert\tanh(x)\vert \le 1$, one has
\begin{align}
\label{eqn:mafia}
|\tanh(x)| \asymp |x| \wedge 1. 
\end{align}
We claim that this leads to the following consequence: 
\begin{align}
	\label{eqn:tanh_bound}
	\big\|\tanh(\pi_tx_t)\big\|_2 
  \asymp \pi_t; 
\end{align}
for the moment, let us first take this as given and we shall come back to its proof after establishing the property~\eqref{eq:eta_derivative}. 
In view of this claim \eqref{eqn:tanh_bound}, we find that 
 \begin{align}
  \gamma_t \defn \|\tanh(\pi_tx_t)\|_2^{-1} \asymp\pi_t^{-1}.\label{eqn:tanh_bound-123}
 \end{align}
In order to prove property~\eqref{eq:eta_derivative}, it suffices to recall the definition of $\eta_t$ as in expression~\eqref{defi:eta}. For any $x\in \real$, direct calculations give 
\begin{subequations}
\label{eqn:super-basic-derivative}
\begin{align}
\eta_t(x) & = \gamma_t \tanh(\pi_t x) \\
  \eta_{t}^{\prime}(x) &= \gamma_t\pi_t\big(1-\tanh^2(\pi_tx)\big) \\
   \eta_{t}^{\second}(x) &= -2\gamma_t\pi_t^2\tanh(\pi_tx)\big(1-\tanh^2(\pi_tx)\big)\\ 
  \eta_{t}^{\third}(x) &= -2\gamma_t\pi_t^3\big(1-\tanh^2(\pi_tx)\big)\big(1-3\tanh^2(\pi_tx)\big).
\end{align}
\end{subequations}
Combining these identities with \eqref{eqn:mafia}, \eqref{eqn:tanh_bound-123} and the fact that $|\tanh(x)| \leq 1$, one can easily validate expression~\eqref{eq:eta_derivative}. 
It then boils down to justifying the claim \eqref{eqn:tanh_bound}, which we accomplish below.

\paragraph{Proof of relation~\eqref{eqn:tanh_bound}.}
Note that from display~\eqref{eqn:mafia}, one has $\big\|\tanh(\pi_tx_t)\big\|_2 \asymp \big \||\pi_tx_t| \wedge 1 \big\|_2$, 
where both operators $|\cdot|$ and $\wedge$ are applied in an entrywise manner and we overlad the notation $1$ to denote an all-one vector. 
To establish the relation~\eqref{eqn:tanh_bound}, it is sufficient to prove that $\||\pi_tx_t| \wedge 1 \big\|_2 \asymp \pi_{t}$. 
Towards this end, first we invoke \eqref{eqn:pi_t-2} to make the observation that 
\begin{align*}
\big \||\pi_tx_t| \wedge 1 \big\|_2 \le 
	\| \pi_tx_t \|_2 \wedge  \| 1 \|_2  =  \|\pi_tx_t \|_2 \wedge \sqrt{n} \asymp \pi_t \wedge \sqrt{n}\asymp \pi_t .
\end{align*}
In addition, let us introduce an index set $\mathcal{I}$ as follows:
\begin{align}
    \mathcal{I} := \bigg\{i \in [n] ~\Big|~ |\xi_{t-1, i}| \le 0.9|v_{t, i}|\text{ and }|v_{t, i}| \lesssim \frac{1}{\sqrt{n}}\bigg\},
	\label{eq:defn-I-auxiliary}
\end{align}
which clearly satisfies 
\begin{align*}
	0.1|\pi_t v_{t,i}|\leq|\pi_t x_{t,i}|\leq 1.9|\pi_t v_{t,i}|, \qquad \forall i\in \mathcal{I}.
\end{align*}
It then follows that:
\begin{align}
\label{eqn:fran}
\big \||\pi_tx_t| \wedge 1 \big\|_2 
\ge \big \||\pi_tx_t \circ \ind_{\mathcal{I}}| \wedge 1 \big\|_2 
\asymp\||\pi_tv_t \circ \ind_{\mathcal{I}}| \wedge 1 \big\|_2. 
\end{align}

To further control the right-hand side of display \eqref{eqn:fran}, we claim that 
\begin{align}
  \big\||\pi_tv_t \circ \ind_{\mathcal{I}}| \wedge 1 \big\|_2 
  \stackrel{(\mathrm{i})}{\asymp}
  \big\|\pi_tv_t \circ \ind_{\mathcal{I}}\big\|_2 
  \stackrel{(\mathrm{ii})}{\asymp} \pi_t.
	\label{eq:pit-i-ii-part}
\end{align}
In order to see this, relation $(\mathrm{i})$ can be verified using expression~\eqref{eqn:pi_t-2} and the definition of \eqref{eq:defn-I-auxiliary}, as they guarantee that
\begin{align*}
  |\pi_t v_{t,i}|\le |\pi_t|\cdot|v_{t,i}|\lesssim 1, \qquad \forall i\in \mathcal{I}.
\end{align*}
To validate $(\mathrm{ii})$, 
note that on the index set $\mathcal{I}^c$, one has $|v_{t,i}|\lesssim |\xi_{t-1,i}|$. Therefore, it holds that 
\begin{align}\label{eq:bound-vt-Ic}
\Vert\pi_t v_t\circ\ind_{\mathcal{I}^c}\Vert_2\lesssim \Vert\pi_t \xi_{t-1}\circ\ind_{\mathcal{I}^c}\Vert_2\lesssim \pi_t\ltwo{\xi_{t-1}}\lesssim \pi_t \sqrt{\frac{1}{\log n}}, 
\end{align}
where we recall our assumption 
$\ltwo{\xi_{t-1}}\lesssim 
\sqrt{\frac{1}{\log n}}.$
Equipped with the above calculation, we can recall $\ltwo{v_t}\asymp 1$ from expression~\eqref{eqn:v-norm} to obtain 
\begin{align}
\label{eq:bound-vt-I}
\Vert\pi_t v_t\circ\ind_{\mathcal{I}}\Vert_2^2 = \Vert\pi_t v_t\Vert^2_2 - \Vert\pi_t v_t\circ\ind_{\mathcal{I}^c}\Vert^2_2
\asymp
\bigg(1-O\bigg(\frac{1}{\log n}\bigg)\bigg)\pi^2_t\asymp \pi^2_t,
\end{align}
as claimed in Part (ii) of \eqref{eq:pit-i-ii-part}. 

\paragraph{Proof of property~\eqref{eq:eta_derivative_small}.} 
To study the derivatives of $\eta_{t}(x_t)$, we first consider the parameters $\pi_{t}$ and $\gamma_{t}$. 
Given that $\Vert\xi_{t-1}\Vert_2$ satisfies the expression~\eqref{eq:induction-lemma-xi}, when $t\lesssim \frac{\log n}{\lambda-1}$,  
it holds true that  
\begin{align}
\label{eq:xi-varsigma}
  \Vert\xi_{t-1}\Vert_2\lesssim
  \sqrt{\frac{t^3\log n}{n}}\lesssim    \sqrt{\frac{\log^4 n}{n(\lambda-1)^3}}.
\end{align}
Under the assumption $\lambda-1\gtrsim n^{-1/9}\log n$, we see that $\ltwo{\xi_{t-1}}$ also satisfies 
\begin{align*}
  \ltwo{\xi_{t-1}} \lesssim (\lambda-1)\cdot\sqrt{\frac{\log^4 n}{n}\cdot \frac{1}{(\lambda-1)^5}}\lesssim(\lambda-1)\cdot\sqrt{\frac{\log^4 n}{n}\cdot \frac{n^{5/9}}{\log^5 n}}\lesssim (\lambda-1)n^{-0.2}.
\end{align*}
Similarly, it is also easily seen that $$\Big(\frac{t\log n}{n}\Big)^{1/2}\lesssim (\lambda-1)n^{-0.2}.$$
Taking these together with the relation~\eqref{eq:second-bound-xt-minus-1} in the proof of Lemma~\ref{lem:parameters} and the assumption $\alpha_{t} \lesssim \sqrt{\lambda-1} \, n^{-0.1}$ yields 
\begin{align}
	\pi_{t} & \leq\sqrt{n}\alpha_{t}+O\bigg(\sqrt{n}\Big(\|\xi_{t-1}\|_{2}+\sqrt{\frac{t\log n}{n}}\Big)^{1/2}\bigg)\lesssim\sqrt{\lambda-1}\,n^{0.4}.
\end{align}
Similarly, the relation~\eqref{eqn:gamma_t_arrange} of Lemma~\ref{lem:parameters} combined with the assumption $\alpha_t\leq \sqrt{\lambda-1}n^{-0.1}$ and the above bounds leads to  
\begin{align*}
  (\gamma_t\pi_t)^{-2} = \alpha_t^2 + 1 + O\bigg(\frac{\pi_t^2}{n} + \Vert \xi_{t-1}\Vert_2 + \sqrt{\frac{t\log n}{n}}\bigg) = 1 + O\lt((\lambda-1)n^{-0.2}\rt),
\end{align*}
which by direct calculation also gives $\gamma_t\pi_t=1+O((\lambda-1)n^{-0.2})$ under our assumption on $\lambda-1$. 

Armed with the above properties, some algebra together with \eqref{eqn:super-basic-derivative} further results in   
\begin{align*}
\eta_t'(x) &= \gamma_t\pi_t(1-\tanh^2(\pi_t x) ) 
= 1+O\big((\lambda-1)n^{-0.2}\log n \big),\\
|\eta_t''(x)| &\lesssim  \pi_t\cdot |\pi_t x|\lesssim (\lambda-1)n^{0.8}|x|,
\end{align*}
where we invoke the relation $|1-\tanh^2(\pi_t x)|\asymp 1$ for any $|x|\lesssim\sqrt{\frac{\log n}{n}}$.

\paragraph{Proof of property~\eqref{eq:eta_small}.}
Again when $t\lesssim \frac{\log n}{\lambda-1}$, $\Vert\xi_{t-1}\Vert_2$ satisfies inequality~\eqref{eq:xi-varsigma}.
Taking this collectively with the assumption $\alpha_t\lesssim\sqrt{\lambda-1}n^{-1/4}$ and the relation~\eqref{eq:second-bound-xt-minus-1} in the proof of Lemma~\ref{lem:parameters}, we obtain 
\begin{align}
\pi_{t}&\leq\sqrt{n}\alpha_{t}+O\bigg(\sqrt{n}\Big(\|\xi_{t-1}\|_{2}+\sqrt{\frac{t\log n}{n}}\Big)^{1/2}\bigg)\lesssim\sqrt{\lambda-1}\,n^{1/4}+(\lambda-1)^{-3/4}n^{1/4}\log n \notag\\
	& \asymp(\lambda-1)^{-3/4}n^{1/4}\log n.
\end{align}
Similarly, the relation~\eqref{eqn:gamma_t_arrange} of Lemma~\ref{lem:parameters} together with our assumption on $\lambda - 1$ yields 
\begin{align}
\gamma_t\pi_t = 1 + O\lt(\frac{\log^2 n}{\sqrt{n(\lambda-1)^3}}\rt).
\end{align}
To derive property~\eqref{eq:eta_small}, we make note of some simple facts that $\tanh(0)=\tanh^{\second}(0)=\tanh^{\prime\prime\prime\prime}(0)=0$, $\tanh^{\prime}(0)=1$, $\tanh^{\third}(0)=-2$ and $|\tanh^{(5)}(x)|\leq K$ for some constant $K$. As a result, the mean value theorem ensures that for any $x$, there exists a quantity $c$ such that
\begin{align*}
  \tanh(x) = x-\frac{1}{3}x^3+cx^5\qquad\text{for some }\quad0\leq c\leq K', 
\end{align*}
for $K'=K/120$. 
Based on the calculations above, we can conclude that
\begin{align*}
\eta_t(x) = \gamma_t\tanh(\pi_t x)=(1-c_0)\pi_t^{-1}\tanh(\pi_t x) = (1-c_0)\Big(x-\frac{1}{3}\pi_t^2 x^3 + c_x\Big),
\end{align*}
where $c_0$ and $c_x$ are some quantities obeying 
\begin{align*}
|c_0| \lesssim \frac{\log^2 n}{\sqrt{n(\lambda-1)^3}} 
\qquad\text{and}\qquad
|c_x|\lesssim \pi_t^4 |x|^5 \lesssim \frac{n|x|^5\log^4 n}{(\lambda-1)^3}. 
\end{align*}
This completes the proof of the desired property. 


\subsection{Proof of Lemma~\ref{lem:kappa}}
\label{sec:pf-kappa}

Without loss of generality, throughout this proof, we assume $\alpha_{t} > 0.$ 
Before we begin to bound $\kappa_{t}$, let us simplify the term of interest a little bit.  
First of all, as each entry follows $\vstar_i \stackrel{\mathrm{i.i.d.}}{\sim}\mathsf{Unif}\{\frac{1}{\sqrt{n}}, -\frac{1}{\sqrt{n}}\}$, one can easily derive 
\begin{align*}
&\bigg\langle\int\left[ x\eta_t^{\prime}\left(\alpha_tv^\star+\frac{1}{\sqrt{n}}x\right)-\frac{1}{\sqrt{n}}\eta_t^{\prime\prime}\left(\alpha_tv^{\star}+\frac{1}{\sqrt{n}}x\right)\right]^2\varphi_n(\dx)\bigg\rangle\\
&\qquad=\int\left[ x\eta_t^{\prime}\left(\alpha_tv^\star_1+\frac{1}{\sqrt{n}}x\right)-\frac{1}{\sqrt{n}}\eta_t^{\prime\prime}\left(\alpha_tv^{\star}_1+\frac{1}{\sqrt{n}}x\right)\right]^2\varphi(\dx)\\
&\qquad =\int\left[ x\eta_t^{\prime}\left(\frac{1}{\sqrt{n}}(\alpha_t+x)\right)-\frac{1}{\sqrt{n}}\eta_t^{\prime\prime}\left(\frac{1}{\sqrt{n}}(\alpha_t+x)\right)\right]^2\varphi(\dx),
\end{align*}
where $\varphi(\cdot)$ is the p.d.f.~of $\mathcal{N}\left(0,1\right)$, and we have used the fact that $\eta_t(\cdot)$ is symmetric about $0$ and the integration is over the distribution $\mathcal{N}(0, I_n)$. Similarly, we obtain 
\begin{align*}
\bigg\langle\int\left[\eta_t^{\prime}\left(\alpha_tv^\star+\frac{1}{\sqrt{n}}x\right)\right]^2\varphi_n(\dx)\bigg\rangle=\int\left[\eta_{t}^{\prime}\Big(\frac{1}{\sqrt{n}}(\alpha_t + x)\Big)\right]^2\varphi(\dx).
\end{align*}
Therefore, it holds that
\begin{align}
\kappa_t^2 = 
\max\bigg\{\int |I_1(x)|^2\varphi(\dx), ~
\int |I_2(x)|^2\varphi(\dx)\bigg\},
	\label{eq:kappa-t-simplified}
\end{align}
where we define
\begin{align*}
I_1(x) &\defn x\eta_{t}^{\prime}\Big(\frac{1}{\sqrt{n}}(\alpha_t + x)\Big)  - \frac{1}{\sqrt{n}}\eta_{t}^{\second}\Big(\frac{1}{\sqrt{n}}(\alpha_t + x)\Big),\\
I_2(x) &\defn\eta_{t}^{\prime}\Big(\frac{1}{\sqrt{n}}(\alpha_t + x)\Big).
\end{align*}

We now proceed to the proof of Lemma~\ref{lem:kappa}, and begin by restricting our attention to the range $t < \varsigma \wedge \frac{\log n}{c(\lambda - 1)}$. We divide into two cases depending on the value of $\alpha_{t}$.

\paragraph{Case I: $\alpha_t\lesssim\sqrt{\lambda^2-1} \, n^{-0.1}$.} 
Let us introduce an event $\mathcal{A} \defn \left\{x: \vert x\vert\leq\sqrt{24\log n}\right\}$. 
For $x\sim\mathcal{N}\left(0,1\right)$, it is easily verified that $P(\mathcal{A})=1-O(n^{-12})$. 
Conditional on $\mathcal{A}$ and assuming $\alpha_t \lesssim \sqrt{\lambda^2-1}\, n^{-0.1}$, one has 
\begin{align*}
\Big\vert \frac{1}{\sqrt{n}}(\alpha_t + x)\Big\vert\lesssim\sqrt{\frac{\log n}{n}}.
\end{align*}
Meanwhile, recall that $\sqrt{\lambda+1}\asymp 1$, and therefore $\alpha_t\lesssim \sqrt{\lambda^2-1}\,n^{-0.1}$ is equivalent to $\alpha_t\lesssim \sqrt{\lambda-1}\,n^{-0.1}$. As a result, according to the property~\eqref{eq:eta_derivative_small} established in Lemma~\ref{lem:eta_property}, we see that: when $|z| \lesssim \sqrt{\frac{\log n}{n}}$, one has  
\begin{align}
	\eta_t^{\prime}(z)= 1+ O\big((\lambda-1)n^{-0.2}\log n\big) \qquad \text{and}\qquad |\eta_t''(z)|\lesssim (\lambda-1)n^{0.8} |z|.
\end{align}
Hence, for all $x$ residing within $\mathcal{A}$, we can bound the difference between $I_1(x)$ and $x$ uniformly as follows: 
\begin{align}
|I_{1}(x)-x| & \notag\leq\lt|\eta_{t}^{\prime}\Big(\frac{1}{\sqrt{n}}(\alpha_{t}+x)\Big)-1\rt||x|+\frac{1}{\sqrt{n}}\lt|\eta_{t}^{\second}\Big(\frac{1}{\sqrt{n}}(\alpha_{t}+x)\Big)\rt|\\
 & \notag\lesssim\lt((\lambda-1)n^{-0.2}\log n\rt)|x|+(\lambda-1)n^{0.3}\lt|\frac{1}{\sqrt{n}}(\alpha_{t}+x)\rt|\\
 & \notag\lesssim\lt((\lambda-1)n^{-0.2}\log n\rt)|x|+(\lambda-1)n^{0.3}\lt(\frac{1}{\sqrt{n}}+\frac{1}{\sqrt{n}}|x|\rt)\\
 & \lesssim\lt((\lambda-1)n^{-0.2}\log n\rt)(1+|x|).  
\end{align}
In addition, consider any $x\in\mathbb{R}$. Recalling the relation~\eqref{eq:eta_derivative} of Lemma~\ref{lem:eta_property} --- which reveals that $|\eta_t'(x)|\lesssim 1$, $|\eta_t''(x)|\lesssim \sqrt{n}$ --- we observe the crude bound that $|I_1(x)|\lesssim |x|+1$. 
Putting the preceding two bounds together, we arrive at
\begin{align}
\label{eqn:I-1}
\notag \int |I_1(x)|^2\varphi(\dx) &= \int x^2\varphi(\dx) + \int (I_1(x)-x)(I_1(x)+x) \varphi(\dx)  \\
\notag&= 1 + \int_{\mathcal{A}} (I_1(x)-x)(I_1(x)+x) \varphi(\dx) + \int_{\mathcal{A}^{\mathsf{c}}} (I_1(x)-x)(I_1(x)+x) \varphi(\dx) \\
\notag&= 1 + O\lt(\int_\mathcal{A} (\lambda-1)n^{-0.2}\log n\cdot (|x|+1)^2\varphi(\dx) + \int_{\mathcal{A}^c} (|x|+1)^2\varphi(\dx)\rt) \\
&= 1 + O\lt((\lambda-1)n^{-0.2}\log n\rt),
\end{align}
where the last equality utilizes the fact that for $c_n=\sqrt{24\log n}$, one has 
\begin{align*}
\int_{\mathcal{A}^c}x^2\varphi(\dx)
= 2 \int_{c_n}^{\infty} x^2\varphi(\dx) \lesssim \int_{c_n}^{\infty} x^2 \exp\lt(-\frac{1}{2}x^2\rt)\dx 
%
%
	\lesssim \frac{\log n}{n^{12}}.
\end{align*} 

Regarding the other term $I_2(x)$, relation~\eqref{eq:eta_derivative} of Lemma~\ref{lem:eta_property} implies that $|I_2(x)|\lesssim 1$. 
Furthermore, if $|x|\leq\sqrt{24\log n}$, in view of the relation~\eqref{eq:eta_derivative_small} we have $|I_2(x)| = 1+ O((\lambda-1)n^{-0.2}\log n)$. Putting these together, we obtain
\begin{align}\begin{aligned}\label{eqn:I-2}
\int |I_2(x)|^2\varphi(\dx) &= \int_\mathcal{A} |I_2(x)|^2\varphi(\dx) +\int_{\mathcal{A}^c} |I_2(x)|^2\varphi(\dx) \\ 
	&= 1+ O\big((\lambda-1)n^{-0.2}\log n \big) + O\big( P(\mathcal{A}^c) \big) \\
&= 1+ O\big((\lambda-1)n^{-0.2}\log n \big).
\end{aligned}\end{align}

Combining inequalities \eqref{eqn:I-1} and~\eqref{eqn:I-2} then leads to 
\begin{align*}
\kappa_t^2 = 1 + O\lt((\lambda-1)n^{-0.2}\log n\rt)= 1+o\lt(\frac{\lambda-1}{\log n}\rt).
\end{align*}

\paragraph{Case II: $\sqrt{\lambda^2-1}\,n^{-0.1}\lesssim\alpha_t\lesssim\sqrt{\lambda^2-1}$.} 
We first note that under the assumption~\eqref{eq:induction}, the following relation holds for $\ltwo{\xi_{t-1}}$ when $t\lesssim\frac{\log n}{\lambda-1}$: 
\begin{align}\label{eqn:bound-xi-t}
\ltwo{\xi_{t-1}}\lesssim \sqrt{\frac{t^3\log n}{n}}\lesssim\sqrt{\frac{\log^4 n}{n(\lambda-1)^3}}.
\end{align}
We recall the basic facts obtained in property~\eqref{eqn:super-basic-derivative}, and as a result, 
\begin{align}
\label{eqn:I1-allegro}
\int |I_1(x)|^2\varphi(\dx) 
%
 & =\int\lt[\lt(\gamma_{t}\pi_{t}x+\frac{2}{\sqrt{n}}\gamma_{t}\pi_{t}^{2}\tanh\Big(\frac{\pi_{t}}{\sqrt{n}}\big(\alpha_{t}+x\big)\Big)\rt)\cdot\Big(1-\tanh^{2}\Big(\frac{\pi_{t}}{\sqrt{n}}\big(\alpha_{t}+x\big)\Big)\Big)\rt]^{2}\varphi(\dx) \\
 %
 \label{eqn:I2-allegro}
\int |I_2(x)|^2\varphi(\dx) &= \int \Big[\gamma_t\pi_t\big(1-\tanh^{2}\Big(\frac{\pi_{t}}{\sqrt{n}}\big(\alpha_{t}+x\big)\Big)\Big)\Big]^2 \varphi(\dx).
\end{align}
It then comes down to controlling the right-hand side of the above two expressions, under the condition that  
$\sqrt{\lambda^2-1}n^{-0.1}\lesssim \alpha_t \le \sqrt{\lambda^2-1}.$

For notational convenience, let us introduce additional shorthand notation as follows:   
\begin{align*}
J_1(x) &\defn \gamma_{t}\pi_{t}x+\frac{2}{\sqrt{n}}\gamma_{t}\pi_{t}^{2}\tanh\Big(\frac{\pi_{t}}{\sqrt{n}}\big(\alpha_{t}+x\big)\Big), \\
J_2(x) &\defn \mu_{t}\alpha_{t}x+2\mu_{t}\alpha_{t}^{2}\tanh\Big(\alpha_{t}\big(\alpha_{t}+x\big)\Big),\\
K_1(x)   &\defn 1-\tanh^{2}\Big(\frac{\pi_{t}}{\sqrt{n}}\big(\alpha_{t}+x\big)\Big),\\
K_2(x)   &\defn 1-\tanh^{2}\Big(\alpha_t\big(\alpha_{t}+x\big)\Big)\\
\mu_{t} &\defn \lt[\int\tanh^{2}\big(\alpha_{t} (\alpha_{t}+x)\big)\varphi(\dx)\rt]^{-\frac{1}{2}}.
\end{align*}
With this set of notation in place, it follows from \eqref{eqn:I1-allegro} and a little algebra that 
\begin{align}
\notag \int |I_1(x)|^2\varphi(\dx) 
&= \int [J_1(x)K_1(x)]^2\varphi(\dx) \\
\notag &= \int [J_2(x)K_2(x)]^2\varphi(\dx) +O\lt(\int \lt\vert J_2(x)\rt\vert^2\lt\vert K_1(x)-K_2(x)\rt\vert\lt\vert K_1(x)+K_2(x)\rt\vert\varphi(\dx)\rt)\\
&\qquad\qquad+O\lt(\int\vert J_1(x)-J_2(x)\vert\vert J_1(x)+J_2(x)\vert\vert K_1(x)\vert^2\varphi(\dx)\rt).
%
\label{eqn:I1-adagio}
\end{align}
To bound $\int |I_1(x)|^2\varphi(\dx)$, it is thus sufficient to control these three terms on the right-hand side of \eqref{eqn:I1-adagio} separately. 
Before proceeding, we find it helpful to make note of several preliminary properties. 

\begin{itemize}
\item First, we would like to show that
\begin{align}
\label{eqn:ut-num}
  \mu_t\asymp \alpha_t^{-1}.
\end{align}
In order to see this, note that the elementary fact  $\tanh^2(x) \leq x^{2}$ implies that 
\begin{align*}
  \int\tanh^2(\alpha_t(\alpha_t+x))\varphi(\dx)\leq \int\alpha_t^2(\alpha_t+x)^2\varphi(\dx)\asymp \alpha_t^2, \qquad \text{for } \alpha_t \leq \lambda.
\end{align*}
On the other hand, when $|x| \leq 1/2$, $\alpha_{t}\in (0,\lambda]$ and $\lambda \in (1,1.2]$, one has $\alpha_t(\alpha_t+x) \in [-0.0625, 2.04]$. 
		Clearly, for any $z\in [-0.0625, 2.04]$, we have $\tanh^2(z)/z^2 \geq  0.22$, and as a consequence,  
\begin{align*}
\int\tanh^2(\alpha_t(\alpha_t+x))\varphi(\dx)&\geq \int\tanh^2(\alpha_t(\alpha_t+x))\textbf{1}(\vert x\vert\leq 0.5)\varphi(\dx)\\
&\gtrsim\int\alpha_t^2(\alpha_t+x)^2\textbf{1}(\vert x\vert\leq 0.5)\varphi(\dx)
\asymp \alpha_t^2.
\end{align*}
The preceding two bounds taken collectively justify the claim \eqref{eqn:ut-num}. 

\item Additionally, based on equation~\eqref{eqn:pi_t} of Lemma~\ref{lem:parameters} and the bound \eqref{eqn:bound-xi-t}, we have 
	\begin{align}
\pi_t &= \sqrt{n}\alpha_t+\frac{\sqrt{n}}{\alpha_t}O\bigg(\Vert\xi_{t-1}\Vert_2+\sqrt{\frac{t\log n}{n}}\bigg)=\alpha_t\sqrt{n} +O\lt(\sqrt{\frac{\log^4 n}{\alpha_t^2(\lambda-1)^3}}\rt) \label{equ:pi-alpha-1}\\
&=\alpha_t\sqrt{n}\lt(1+ O\lt(\sqrt{\frac{\log^4 n}{n\alpha_t^4(\lambda-1)^3}}\rt)\rt)\asymp \alpha_t\sqrt{n}  \label{equ:pi-alpha}
\end{align}
as long as $t\lesssim\frac{\log n}{\lambda-1}$. 
Here, the last inequality holds since 
\begin{align*}
\sqrt{\frac{\log^4 n}{(\lambda-1)^3}}\lesssim(\lambda-1)\cdot \frac{\log^2 n}{(\lambda-1)^{5/2}}=o\big((\lambda-1)n^{0.3}\big)=o\big(\alpha_t^2\sqrt{n} \big), 
\end{align*}
provided that $\alpha_t\gtrsim\sqrt{\lambda^2-1}n^{-0.1}$ and $\lambda-1\gtrsim n^{-1/9}\log n$.

\item Moreover, by combining equation~\eqref{eqn:gamma_t} of Lemma~\ref{lem:parameters} with \eqref{equ:pi-alpha-1}, \eqref{equ:pi-alpha} and \eqref{eqn:bound-xi-t}, we see that 
\begin{align*}
\gamma_t^{-2}&
	= n\int\tanh^2\lt(\frac{\pi_t}{\sqrt{n}}(\alpha_t +x)\rt)\varphi(\dx) +  O\lt(\alpha_t^2\sqrt{\frac{n\log^4 n}{(\lambda-1)^3}}\rt)\\
&=n\int\tanh^2\lt(\alpha_t(\alpha_t +x)\rt)\varphi(\dx) + O\lt(\sqrt{\frac{n\log^4 n}{(\lambda-1)^3}}\rt)\\
&=n\mu_t^{-2} + O\lt(\sqrt{\frac{n\log^4 n}{(\lambda-1)^3}}\rt). 
\end{align*}
		Here, the second equality arises from the following fact (see also \eqref{eq:integral-pi}): 
\begin{align}\label{equ:tanh-pi-alpha}
&\lt\vert\int\tanh^2\lt(\frac{\pi_t}{\sqrt{n}}(\alpha_t +x)\rt)-\tanh^2\lt(\alpha_t(\alpha_t +x)\rt)\varphi(\dx)\rt\vert \notag\\
&\le \lt\vert1 - \frac{\alpha_t^2n}{\pi_t^2}\rt\vert\int\tanh^2\lt(\frac{\pi_t}{\sqrt{n}}(\alpha_t +x)\rt)\varphi(\dx) + \lt\vert\int\frac{\alpha_t^2n}{\pi_t^2}\tanh^2\lt(\frac{\pi_t}{\sqrt{n}}\lt(\alpha_t + x\rt)\rt) - \tanh^2\lt(\alpha_t\lt(\alpha_t + x\rt)\rt)\varphi(\dx)\rt\vert \notag\\
&\lesssim \lt\vert1 - \frac{\alpha_t^2n}{\pi_t^2}\rt\vert\alpha_t^2 + \lt|\frac{\pi_t^2}{\alpha_t^2n} - 1\rt|\alpha_t^4 
=O\lt(\sqrt{\frac{\log^4 n}{n(\lambda-1)^3}}\rt).
\end{align}
As a result, we can express $\gamma_t$ in term of $\mu_t$ as
\begin{align}\label{eqn:gamma-mu}
\gamma_t &= \lt[n\mu_t^{-2} + O\lt(\sqrt{\frac{n\log^4 n}{(\lambda-1)^3}}\rt)\rt]^{-1/2}
=\frac{\mu_t}{\sqrt{n}}\lt(1 + O\lt(\mu_t^2\sqrt{\frac{\log^4 n}{n(\lambda-1)^3}}\rt)\rt)\asymp\frac{\mu_t}{\sqrt{n}}.
\end{align}

\item With \eqref{equ:pi-alpha} and~\eqref{eqn:gamma-mu} in place, a little algebra further leads to 
\begin{align}
\lt\vert\gamma_t\pi_t-\mu_t\alpha_t\rt\vert &\lesssim \mu_t\alpha_t\sqrt{\frac{\log^4 n}{n(\lambda-1)^3}}\left(\frac{1}{\alpha_t^2}+\mu_t^2\right)
\lesssim \frac{1}{\alpha_t^2}\sqrt{\frac{\log^4 n}{n(\lambda-1)^3}}, \label{equ:bound-kappa-t-1}\\
\lt\vert\frac{1}{\sqrt{n}}\gamma_t\pi_t^2-\mu_t\alpha_t^2\rt\vert &\lesssim \mu_t\alpha_t^2\sqrt{\frac{\log^4 n}{n(\lambda-1)^3}}\left(\frac{1}{\alpha_t^2}+\mu_t^2\right)
\lesssim \frac{1}{\alpha_t}\sqrt{\frac{\log^4 n}{n(\lambda-1)^3}}. \label{equ:bound-kappa-t-2}
\end{align}
In addition, according to the relation~\eqref{eqn:gamma_t_arrange} in Lemma~\ref{lem:parameters}, we have
\begin{align*}
\gamma_t^2\pi_t^2 = (\alpha_t^2+1)^{-1}+O\bigg(\frac{\pi_t^2}{n}+\Vert\xi_{t-1}\Vert_2+\sqrt{\frac{t\log n}{n}}\bigg),
\end{align*}
and using similar analysis as for~\eqref{eqn:tanh_pi_t} yields
\begin{align*}
\mu_t^{-2} = \int \tanh^2\lt(\alpha_t(\alpha_t + x)\rt)\varphi(\dx) &= \alpha_t^2(\alpha_t^2 + 1) + O\lt(\alpha_t^4\rt).
\end{align*}
These two bounds taken together with a little algebra leads to
\begin{align}
\lt\vert\gamma_t\pi_t-\mu_t\alpha_t\rt\vert \lesssim \frac{\pi_t^2}{n}+\alpha_t^2+\Vert\xi_{t-1}\Vert_2+\sqrt{\frac{t\log n}{n}} \lesssim \alpha_t^2 + \sqrt{\frac{\log^4 n}{n(\lambda-1)^3}}
\lesssim \alpha_t^2,
\end{align}
recognizing the range $\sqrt{\lambda^2-1}\,n^{-0.1}\lesssim\alpha_t\lesssim\sqrt{\lambda^2-1}$
and our assumption $\lambda - 1 \gtrsim n^{-1/9}\log n.$
Combining the above bound with~\eqref{equ:bound-kappa-t-1} gives
\begin{align}
\lt\vert\gamma_t\pi_t-\mu_t\alpha_t\rt\vert &\lesssim \alpha_t^2 \wedge \frac{1}{\alpha_t^2}\sqrt{\frac{\log^4 n}{n(\lambda-1)^3}} \lesssim \left(\frac{\log^4 n}{n(\lambda-1)^3}\right)^{1/4}, \label{eqn:bound-gamma-pi}
\end{align}
where the last inequality follows from the elementary fact that $\min \{a, b\} \leq \sqrt{ab}.$

\item
	Taking inequalities~\eqref{equ:pi-alpha},~\eqref{equ:bound-kappa-t-1} and~\eqref{equ:bound-kappa-t-2} collectively with the fact that $\mu_t\alpha_t^2\asymp \alpha_t$ (cf.~\eqref{eqn:ut-num}) yields
\begin{align}
\label{bound-K1-K2}
\lt\vert K_1(x) - K_2(x)\rt\vert \lesssim & \lt\vert\frac{\pi_t}{\sqrt{n}}-\alpha_t\rt\vert\lt\vert\alpha_t +x\rt\vert \lesssim \frac{1}{\alpha_t}\sqrt{\frac{\log^4 n}{n(\lambda-1)^3}}\cdot (\vert x\vert + \alpha_t),
\end{align}
and 
\begin{align}\label{bound-J1-J2}
\lt\vert J_1(x) - J_2(x)\rt\vert &\leq  \lt\vert\gamma_t\pi_t-\mu_t\alpha_t\rt\vert \vert x\vert + 2\lt\vert\frac{1}{\sqrt{n}}\gamma_t\pi_t^2-\mu_t\alpha_t^2\rt\vert\lt\vert\tanh\lt(\frac{\pi_t}{\sqrt{n}}(\alpha_t+x)\rt)\rt|\notag \\
&\qquad\qquad + \mu_t\alpha_t^2\lt\vert\tanh\lt(\frac{\pi_t}{\sqrt{n}}(\alpha_t+x)\rt) - \tanh\lt(\alpha_t(\alpha_t+x)\rt)\rt\vert \notag \\
&\lesssim \left(\frac{\log^4 n}{n(\lambda-1)^3}\right)^{1/4}\vert x\vert +  \frac{1}{\alpha_t}\sqrt{\frac{\log^4 n}{n(\lambda-1)^3}} + \sqrt{\frac{\log^4 n}{n(\lambda-1)^3}} (\vert x\vert+\alpha_t)\notag\\
&\lesssim  \left(\frac{\log^4 n}{n(\lambda-1)^3}\right)^{1/4}\vert x\vert +  \frac{1}{\alpha_t}\sqrt{\frac{\log^4 n}{n(\lambda-1)^3}}.
\end{align}
\end{itemize}

Equipped the above relations, we are positioned to control the right-hand side of expression~\eqref{eqn:I1-adagio}. 
Combining~\eqref{bound-K1-K2} and \eqref{bound-J1-J2} directly yields 
\begin{align}
\notag \int |I_1(x)|^2\varphi(\dx) 
&= \int [J_2(x)K_2(x)]^2\varphi(\dx) + O\left( \frac{1}{\alpha_t}\sqrt{\frac{\log^4 n}{n(\lambda-1)^3}}\right) 
%
+ O\left(\left(\frac{\log^4 n}{n(\lambda-1)^3}\right)^{1/4} +  \frac{1}{\alpha_t}\sqrt{\frac{\log^4 n}{n(\lambda-1)^3}}\right) \\
&=  \int \Big[\mu_{t}\alpha_{t}x+2\mu_{t}\alpha_{t}^{2}\tanh\Big(\alpha_{t}\big(\alpha_{t}+x\big)\Big)\Big]^2\Big[1-\tanh^{2}\Big(\alpha_t\big(\alpha_{t}+x\big)\Big)\Big]^2\varphi(\dx)\notag\\
&\qquad\qquad
+ O\left( \left(\frac{\log^4 n}{n(\lambda-1)^3}\right)^{1/4} +  \frac{1}{\alpha_t}\sqrt{\frac{\log^4 n}{n(\lambda-1)^3}}\right),
\label{eqn:I1-rondo}
\end{align}
where the first line invokes the simple observations that $\vert J_i(x)\vert\lesssim |x|+\alpha_t$ and $\vert K_i(x)\vert\lesssim 1$ for $i=1,2$.
Similarly, one can derive in the same manner that  
\begin{align}
\int|I_2(x)|^2\varphi(\dx)
&= \mu_{t}^2\alpha_{t}^2\int \lt[1-\tanh^2\Big(\alpha_{t}\big(\alpha_{t}+x\big)\Big)\rt]^2 \varphi(\dx) + O\left( \left(\frac{\log^4 n}{n(\lambda-1)^3}\right)^{1/4} +  \frac{1}{\alpha_t}\sqrt{\frac{\log^4 n}{n(\lambda-1)^3}}\right). 
\end{align}
As it turns out, the main terms in the above two identities satisfy 
\begin{subequations}
\label{eqn:I12-finale}
\begin{align}
\int \Big[\mu_{t}\alpha_{t}x+2\mu_{t}\alpha_{t}^{2}\tanh\Big(\alpha_{t}\big(\alpha_{t}+x\big)\Big)\Big]^2\Big[1-\tanh^{2}\Big(\alpha_t\big(\alpha_{t}+x\big)\Big)\Big]^2\varphi(\dx) &\le 1 - c\alpha_t^2, \\
\mu_{t}^2\alpha_{t}^2\int \lt[1-\tanh^2\Big(\alpha_{t}\big(\alpha_{t}+x\big)\Big)\rt]^2 \varphi(\dx) &\le 1 - c\alpha_t^2,
\end{align}
\end{subequations}
which we shall justify momentarily. 
Combine these results with \eqref{eq:kappa-t-simplified} to conclude that: 
if $\sqrt{\lambda^2-1}\,n^{-0.1}\lesssim\alpha_t\lesssim\sqrt{\lambda^2-1}$, then
\begin{align}
   \kappa_t^2 
  \leq 1 - c \alpha_t^2 + O\left( \left(\frac{\log^4 n}{n(\lambda-1)^3}\right)^{1/4} +  \frac{1}{\alpha_t}\sqrt{\frac{\log^4 n}{n(\lambda-1)^3}}\right) 
  &= 1 + o\left(\frac{\lambda - 1}{\log n}\right),
\end{align}
where the last inequality follows by recognizing that 
\begin{align*}
  \left(\frac{\log^4 n}{n(\lambda-1)^3}\right)^{1/4} +  \frac{1}{\sqrt{\lambda-1}n^{-0.1}}\sqrt{\frac{\log^4 n}{n(\lambda-1)^3}}
  = 
  o\left(\frac{\lambda - 1}{\log n}\right)
\end{align*}
under our assumption $\lambda - 1 \gtrsim n^{-1/9}\log n$.


\paragraph{Proof of relation \eqref{eqn:I12-finale}.} 
To proceed, consider the problem of estimating $\vstar$ (which obeys $\vstar_{i} \sim \mathsf{Unif}\{\pm\frac{1}{n}\}$) from the noisy observation $Y = \alpha_t \vstar + g$, 
where $g \sim \mathcal{N}(0,\frac{1}{n}I_n)$. As alluded to previously, the Bayes-optimal estimate (or minimum mean square estimator (MMSE)) is given by 
\begin{align*}
	\Exs[\vstar \mid Y] = \tanh(\sqrt{n}\alpha_t Y), 
\end{align*}
which satisfies (due to its optimality)  
\begin{align}
  \operatorname{Cor}\lt(\vstar, f(Y)\rt)\leq\operatorname{Cor}\lt(\vstar, \Exs[\vstar \mid Y]\rt),
	\label{eq:corr-MMSE}
\end{align}
for any measurable function $f$; here $\operatorname{Cor}(\cdot,\cdot)$ denotes the correlation of two random vectors.  
In particular, the Bayes-optimal estimator outperforms the identity estimator (i.e., $f(Y)=Y$), so that \eqref{eq:corr-MMSE} translates to
\begin{align}
	\notag \frac{\Exs [\inprod{\vstar}{\alpha_t \vstar + g}]}{\sqrt{\Exs [\ltwo{\alpha_t \vstar + g}^2]}}
   \leq 
	\frac{\Exs [\inprod{\vstar}{\tanh(\sqrt{n}\alpha_t (\alpha_t \vstar + g))}]}{\sqrt{\Exs [\ltwo{\tanh(\sqrt{n}\alpha_t (\alpha_t \vstar + g)}^2]}}
   &=
   \frac{\int \tanh(\alpha_t (\alpha_t + x)) \varphi(dx)}{\sqrt{\int \tanh^2(\alpha_t (\alpha_t + x)) \varphi(\dx)}} \\
   &=
   \sqrt{\int \tanh^2(\alpha_t(\alpha_t + x)) \varphi(\dx)},
\end{align}
where the fist equality holds due to the symmetry  of $\varphi(\cdot)$, and
the second equality holds since $\int \tanh(\alpha^2 + \alpha x) \varphi(\dx) = \int \tanh^2(\alpha^2 + \alpha x) \varphi(\dx)$ (see \citet[Appendix B.2]{deshpande2017asymptotic}). 
As a consequence, the above relation implies that 
\begin{align}
\label{eqn:bayes-opt}
\frac{\alpha_t}{\sqrt{\alpha^2_t + 1 }} = 
	\frac{\Exs [\inprod{\vstar}{\alpha_t \vstar + g}]}{\sqrt{\Exs [\ltwo{\alpha_t \vstar + g}^2]}} \leq \sqrt{\int \tanh^2(\alpha_t(\alpha_t + x)) \varphi(\dx)} = 
  \frac{1}{\mu_t},
\end{align}
which in turns reveals that $$\mu_t\alpha_t \le \sqrt{\alpha_t^2 + 1} = : \gamma.$$ 
Armed with this relation, we can conclude that 
\begin{align*}
&\max\bigg\{\int \Big[\mu_{t}\alpha_{t}x+2\mu_{t}\alpha_{t}^{2}\tanh\Big(\alpha_{t}\big(\alpha_{t}+x\big)\Big)\Big]\Big[1-\tanh^{2}\Big(\alpha_t\big(\alpha_{t}+x\big)\Big)\Big]\varphi(\dx), \\
&\qquad\qquad\qquad\qquad\mu_{t}^2\alpha_{t}^2\int \lt[1-\tanh^2\Big(\alpha_{t}\big(\alpha_{t}+x\big)\Big)\rt]^2 \varphi(\dx)\bigg\} \\
&\le \gamma^2\max\bigg\{\int \Big[x+2\alpha_{t}\tanh\Big(\alpha_{t}\big(\alpha_{t}+x\big)\Big)\Big]\Big[1-\tanh^{2}\Big(\alpha_t\big(\alpha_{t}+x\big)\Big)\Big]\varphi(\dx), \\
&\qquad\qquad\qquad\qquad\int \lt[1-\tanh^2\Big(\alpha_{t}\big(\alpha_{t}+x\big)\Big)\rt]^2 \varphi(\dx)\bigg\} \\
&=: \kappa^2(\gamma, \alpha_t^2).
\end{align*}
As it turns out, this function $\kappa^2(\cdot, \cdot)$ has been studied in \citet{li2022non}; 
more specifically,    \citet[relation~(272)]{li2022non} together with $\gamma \defn \sqrt{\alpha^2_{t}+1}$ indicates 
that $\kappa(\gamma, \alpha_t^2) \le 1 - \frac{\gamma - 1}{12}$, and hence
\begin{align}
	\kappa^2(\gamma, \alpha_t^2) \le 1 - \frac{\gamma - 1}{12} \le 1 - c\alpha_t^2
\end{align}
for some constant $c>0$. 
Here, notice that we view $\gamma$ and $\alpha_t^2$ as $\lambda$ and $\tau$ respectively in \citet[relation~(272)]{li2022non}. Putting everything together completes the proof of the required relation~\eqref{eqn:I12-finale}.

\paragraph{Case III: $\alpha_t\gtrsim \sqrt{\lambda^2-1}$ and $t<\frac{cn(\lambda-1)^5}{\log^2 n}$.}
The calculation of $\kappa_t$ in this case follows from similar arguments as in \citet[Section D.3.4.]{li2022non}. 
The only difference lies in the computing the parameters $\pi_{t}$ and $\gamma_{t}$,
which was done in \citet[Lemma~14]{li2022non} therein but requires a different proof here. Specifically, we aim to show that 
\begin{align}
\label{eqn:happiness}
  \pi_t=\big(1+o(\lambda-1)\big)\alpha_t\sqrt{n}
  \qquad \text{ and } \qquad \gamma_t^{-2}=\big(1+o(\lambda-1)\big)n\int\tanh(\alpha_t(\alpha_t+x))\varphi(\dx).
\end{align}
If these two relations were valid, then one could follow the argument in \cite[Section D.3.4.]{li2022non} verbatim to demonstrate that 
\begin{align*}
  \kappa_t \leq 1 - \frac{1}{15}(\lambda - 1),
\end{align*}
as claimed.

We now present how to prove relation~\eqref{eqn:happiness}. In view of the equation~\eqref{eqn:pi_t} of Lemma~\ref{lem:parameters},  we have 
\begin{align*}
\pi_t&=\alpha_t\sqrt{n}\bigg(1+\frac{1}{\alpha_t^2}O\Big(\ltwo{\xi_{t-1}}+\sqrt{\frac{t\log n}{n}}\Big)\bigg)
=\alpha_t\sqrt{n}\bigg(1+O\bigg(\sqrt{\frac{(t+\log^3 n/(\lambda-1))\log n}{n(\lambda-1)^2}}\bigg)\bigg)\\
&=\alpha_t\sqrt{n}\bigg(1+O\bigg(\sqrt{\frac{(\lambda-1)^2}{\log n}}\bigg)\bigg)=\alpha_t\sqrt{n}(1+o(\lambda-1))
\end{align*}
under the condition $t\lesssim\frac{n(\lambda-1)^5}{\log^2 n}$ and the assumption~\eqref{eq:induction-main}. 
In addition, with the same analysis as inequality~\eqref{equ:tanh-pi-alpha}, we can guarantee that 
\begin{align*}
&\int\tanh^2\lt(\frac{\pi_t}{\sqrt{n}}(\alpha_t+x)\rt)\varphi(\dx)= \int\tanh^2\lt(\alpha_t(\alpha_t+x)\rt)\varphi(\dx) + O\lt(\sqrt{\frac{(t+\log^3 n/(\lambda-1))\log n}{n\alpha_t^2(\lambda-1)^2}}\rt)\\
&=\int\tanh\lt(\alpha_t(\alpha_t+x)\rt)\varphi(\dx) +   O\lt(\sqrt{\frac{(\lambda-1)^2}{\log n}}\rt)=\int\tanh\lt(\alpha_t(\alpha_t+x)\rt)\varphi(\dx) + o(\lambda-1). 
\end{align*}
Then according to equation~\eqref{eqn:gamma_t}, we can reach 
\begin{align*}
\gamma_t^{-2}
&=n\lt(\int\tanh(\alpha_t(\alpha_t+x))\varphi(\dx)+o(\lambda-1)\rt) + n\alpha^2_t\big(1+o(\lambda-1)\big)O\bigg(\sqrt{\frac{(t+\log^3 n/(\lambda-1))\log n}{n}}\bigg)\\
&=\big(1+o(\lambda-1)\big)n\int\tanh(\alpha_t(\alpha_t+x))\varphi(\dx),
\end{align*}
where the last equality follows from the fact that $\int\tanh(\alpha_t(\alpha_t+x))\varphi(\dx)\asymp\alpha_t^2\asymp1$ (see relation~\eqref{eqn:ut-num}). 
This establishes the claim \eqref{eqn:happiness}.

\subsection{Proof of Claim~\eqref{eq:alpha_t}}
\label{sec:claim-alpha-t}

This subsection aims to establish the advertised decomposition~\eqref{eq:alpha_t}. 
To do so, recall that $\{\eta_i(x_i)\}_{1 \le i \le t-1}$ spans the same linear space as $\{z_i\}_{1\le i\le t-1}$ (see \eqref{eq:span-Ut-eta-t-equal} and \eqref{eq:defn-Ut-minus-1}). It is important to notice that $\{\eta_1(x_1),\ldots,\eta_{t-1}(x_{t-1})\}$ are almost orthogonal to each other,  thus forming a set of near-orthonormal basis; 
this property is summarized in the lemma below, whose proof is provided in Section~\ref{proof-linear-ind}. 
\begin{lemma} 
\label{lem:linear-ind} 
Suppose that the assumptions of Theorem~\ref{thm:main} hold.  
With probability at least $1 - O(n^{-11})$,  we have 
\begin{align} 
\label{eq:linear-ind}
	\Big\|\sum_{i=1}^{t} w_i \eta_i(x_i)\Big\|_2 = \big(1+o(1)\big) \|w\|_2
\end{align}
simultaneously for all $t \le \tau_{0}$ and all $w = [w_i]_{1\leq i\leq t} \in \mathbb{R}^{t}$, where $\tau_{0}$ is defined in \eqref{eqn:tau0}.
\end{lemma}

In view of Lemma~\ref{lem:linear-ind} and the fact that $\xi_{t} \in \mathsf{span}(U_{t-1}) = \mathsf{span}\{\eta_1(x_1),\ldots,\eta_{t-1}(x_{t-1})\}$ (cf.~\eqref{eq:span-Ut-eta-t-equal}),
one can write $\xi_{t}$ as a linear combination of $\{\eta_i(x_i)\}_{1 \le i \le t-1}$ as follows: 
\begin{align}
	\xi_{t} = \sum_{k=1}^{t-1} \gamma_t^k\eta_k(x_k),\qquad\text{with }\gamma_t=[\gamma_t^k]_{1\leq k < t}\in\real^{t-1} \text{ obeying }\|\gamma_t\|_2 \asymp \|\xi_t\|_2. \label{eq:xi_decomposition}
\end{align}
Armed with this decomposition, we intend to prove that  
\begin{align}
\label{eqn:alphat-decomp}
\notag \alpha_{t+1} 
= \lambda v^{\star\top}\eta_{t}\lt(x_t\rt) &= \lambda v^{\star\top}\eta_{t}\Big(v_t + \sum _{k=1}^{t-1}\gamma_{t-1}^k\eta_k(x_k)\Big) \\
&=\lambda v^{\star\top}\eta_{t}\lt(v_t\rt) + O\lt(\frac{\log^{2.5} n}{n^{3/4}(\lambda-1)^{1.5}}\rt),
\end{align}
which shall be done as follows. 
\begin{itemize}
	\item In order to see this, first note that $\eta_{t}(\cdot)$ is a Lipschitz function with Lipschitz constant  $O(1)$ (see Lemma~\ref{lem:eta_property}).
		Therefore, for every $t\lesssim\frac{\log n}{\lambda-1}$ we have 
\begin{align}
\notag & \Bigg| v^{\star\top}\eta_{t}\Big(v_t + \sum _{k=1}^{t-1}\gamma_{t-1}^k\eta_k(x_k)\Big) -   v^{\star\top}\eta_{t}\lt(v_t + \sum_{k=1}^{t-1} \gamma_{t-1}^k\eta_k(v_k)\rt) \Bigg| \\
\notag &\leq \Big\|\eta_{t}\Big(v_t + \sum_{k=1}^{t-1} \gamma_{t-1}^k\eta_k(x_k)\Big) - \eta_{t}\Big(v_t + \sum_{k=1}^{t-1}\gamma_{t-1}^k\eta_k(v_k)\Big)\Big\|_2 
\lesssim \sum_{k=1}^{t-1} \left\vert\gamma_{t-1}^k\rt\vert\lt\|\eta_k(x_k) - \eta_k(v_k)\rt\|_2. 
\end{align}
In view of the decomposition~\eqref{eq:xi_decomposition} and the Cauchy-Schwarz inequality, we can further obtain 
\begin{align}
\notag & \Bigg|v^{\star\top}\eta_{t}\Big(v_{t}+\sum_{k=1}^{t-1}\gamma_{t-1}^{k}\eta_{k}(x_{k})\Big)-v^{\star\top}\eta_{t}\Big(v_{t}+\sum_{k=1}^{t-1}\gamma_{t-1}^{k}\eta_{k}(v_{k})\Big)\Bigg|\lesssim\sum_{k=1}^{t-1}\lt\vert\gamma_{t-1}^{k}\rt\vert\|\eta_{k}(x_{k})-\eta_{k}(v_{k})\|_{2}\\
\notag & \quad\lesssim\sum_{k=1}^{t-1}\lt\vert\gamma_{t-1}^{k}\rt\vert\|\xi_{k-1}\|_{2}\leq\|\gamma_{t-1}\|_{2}\Big(\sum_{k=1}^{t-1}\Vert\xi_{k-1}\Vert_{2}^{2}\Big)^{1/2}\\
\notag & \quad\asymp\Vert\xi_{t-1}\Vert_{2}\Big(\sum_{k=1}^{t-1}\Vert\xi_{k-1}\Vert_{2}^{2}\Big)^{1/2}\lesssim\sqrt{\frac{t^{3}\log n}{n}}\cdot\lt(\sum_{k=1}^{t-1}\frac{k^{3}\log n}{n}\rt)^{1/2}\lesssim\frac{\log^{4.5}n}{n(\lambda-1)^{3.5}},
\end{align}
where the last line invokes $\ltwo{\xi_t} \lesssim \sqrt{\frac{t^3\log n}{n}}$ (cf.~\eqref{eq:xi-stage1}) and $t\lesssim\frac{\log n}{\lambda-1}$.

\item In addition, when $|\alpha_t| \lesssim \sqrt{\lambda-1}n^{-1/4}$, 
	we know that $\Big\Vert\sum_{k=1}^{t-1}\beta_{t-1}^{k}\phi_{k}\Big\Vert_\infty \lesssim \sqrt{\frac{t\log n}{n}}$
		conditioned on the event $\{\phi_k\}_{k=1}^{t-1} \in \mathcal{E}$ (defined in Lemma~\ref{lem:concentration} with $\delta = O(n^{-10})$).
It therefore guarantees that 
\begin{subequations}
\label{eqn:v-inf-3424} 
\begin{align}
\lt\vert v_{t,i}\rt\vert 
&= \lt\vert\alpha_t v_i^{\star} + \sum_{k=1}^{t-1}\beta_{t-1}^{k}\phi_{k,i}\rt\vert 
\leq \frac{|\alpha_t|}{\sqrt{n}} + \Big|\sum_{k=1}^{t-1}\beta_{t-1}^{k}\phi_{k,i}\Big| \lesssim  \sqrt{\frac{t\log n}{n}}, 
	\label{eqn:v-inf-342478}\\
\lt\vert v_{t,i} + \sum_{k=1}^{t-1} \gamma_{t-1}^k\eta_k(v_{k, i})\rt\vert &\lesssim \lt\vert v_{t,i}\rt\vert + \sum_{k=1}^{t-1} \lt\vert \gamma_{t-1}^k\rt\vert\lt\vert v_{k, i}\rt\vert \lesssim \lt\vert v_{t,i}\rt\vert + \ltwo{\gamma_{t-1}}\ltwo{\widetilde{v}_{t-1,i}}\notag\\
&\lesssim  \sqrt{\frac{t\log n}{n}} + \sqrt{\frac{t^3\log n}{n}}\cdot \sqrt{\frac{t^2\log n}{n}}\lesssim\sqrt{\frac{t\log n}{n}},
\label{eqn:v-inf-342456}
\end{align}
\end{subequations}
		for every $1 \le i < t$, where we denote $\widetilde{v}_{t-1,i}\defn (v_{1,i},v_{2,i},\ldots,v_{t-1,i})$. 
		To see why \eqref{eqn:v-inf-342456} is valid, 
		we note that the first inequality applies Lemma~\ref{lem:eta_property}, the second inequality results from the Cauchy-Schwarz inequality, 
		whereas the last line makes use of \eqref{eqn:v-inf-342478} and the fact $\ltwo{\gamma_{t-1}}\asymp \ltwo{\xi_{t-1}}\lesssim \sqrt{\frac{t^3 \log n}{n}}$ (cf.~\eqref{eq:xi_decomposition} and \eqref{eq:xi-stage1}).
In addition, given that $t \lesssim \log n/(\lambda - 1)$ and $\lambda-1\gtrsim n^{-1/9}\log n$, we have $t^8\lesssim n/\log n$.
Repeating the argument for inequality~\eqref{eqn:ny} in the proof of Lemma~\ref{lem:linear-ind}, 
one can ensure that for any $2\le k\le 14$,
\begin{subequations}
\label{enq:kpower}
\begin{align}
\sum_{i = 1}^n \lt\vert v_{t,i}\rt\vert^k 
= \sum_{i = 1}^t \lt\vert v_{t,(i)}\rt\vert^k 
+ \sum_{i = t+1}^n \lt\vert v_{t,(i)}\rt\vert^k 
\lesssim
\lt(\frac{\log n}{n}\rt)^{k/2-1}, 
\end{align}
and 
\begin{align}
\sum_{i = 1}^n \lt\vert v_{t,i} + \sum_{k=1}^{t-1} \gamma_{t-1}^k\eta_k(v_{k, i})\rt\vert^k &\lesssim  
\left(1+t^{2k}\lt(\frac{\log n}{n}\rt)^{k/2}\right) \lt(\frac{\log n}{n}\rt)^{k/2-1}\lesssim \lt(\frac{\log n}{n}\rt)^{k/2-1}.
\end{align}
\end{subequations}
Here, we have made the observation that
\begin{align*}
  &\sum_{i = 1}^n \lt\vert \sum_{k=1}^{t-1} \gamma_{t-1}^k\eta_k(v_{k, i})\rt\vert^k \\
  &\lesssim 
   \sum_{i = 1}^n \|\gamma_{t-1}\|_2^k \Big(\sum_{k=1}^{t-1} \eta_k^2(v_{k, i})\Big)^{k/2}
  \lesssim 
  \sum_{i = 1}^n \left(\frac{t^3\log n}{n}\right)^{k/2}   \Big(\sum_{k=1}^{t-1} v^2_{k, i}\Big)^{k/2}\\
  &\lesssim \left(\frac{t^3\log n}{n}\right)^{k/2} \sum_{i = 1}^n \Big(\sum_{k=1}^{t-1} v^2_{k, (i)}\Big)^{k/2} \\
  &\lesssim \left(\frac{t^3\log n}{n}\right)^{k/2}\cdot t\left(\frac{t^2\log n}{n}\right)^{k/2} + \left(\frac{t^3\log n}{n}\right)^{k/2} \Big(\sum_{k=1}^{t-1} v^2_{k, (t+1)}\Big)^{k/2-1}\cdot\sum_{i = t+1}^n \Big(\sum_{k=1}^{t-1} v^2_{k, (i)}\Big) \\
  &\lesssim t^{\frac{5k}{2}+1}\lt(\frac{\log n}{n}\rt)^{k} + t^{2k}\lt(\frac{\log n}{n}\rt)^{k-1}\\
  &\lesssim t^{2k}\lt(\frac{\log n}{n}\rt)^{k-1},
\end{align*}
where the second line uses the fact $\ltwo{\gamma_{t-1}}\asymp \ltwo{\xi_{t-1}}\lesssim \sqrt{\frac{t^3 \log n}{n}}$ (cf.~\eqref{eq:xi_decomposition} and \eqref{eq:xi-stage1}) and 
Lemma~~\ref{lem:eta_property}, 
the ante-penultimate line invokes inequality~\eqref{eqn:v-inf-3424};
the penultimate line follows from the fact that $\ltwo{v_k}\lesssim 1$ (see e.g.~\eqref{eqn:v-norm}) and  
conditional on event $\{\phi_k\}_{k=1}^{t-1}\in\mathcal{E}$, 
\begin{align*}
	|v_{k,(t+1)}| \leq \frac{|\alpha_k|}{\sqrt{n}} + \Big|\sum_{i=1}^{k-1}\beta_{k-1}^i\phi_i\Big|_{(t+1)} 
  \lesssim \sqrt{\frac{\log n}{n}};
\end{align*}
and the last line follows from the fact that $t^{k/2+1}\lesssim t^8\lesssim n/\log n$.
Therefore, combining Lemma~\ref{lem:eta_property} with expression~\eqref{eqn:v-inf-3424} gives 
\begin{align}
&\notag \eta_{t}\lt(v_t + \sum_{k=1}^{t-1} \gamma_{t-1}^k\eta_k(v_k)\rt) - \eta_{t}\lt(v_t\rt) \\
&= \lt(1-c_0\rt)\sum_{k=1}^{t-1} \gamma_{t-1}^k\eta_k(v_k) + 
	O\lt(\pi_t^2\rt) \cdot \lt[\lt(v_t + \sum_{k=1}^{t-1} \gamma_{t-1}^k\eta_k(v_k)\rt)^3 - \lt(v_t\rt)^3\rt]+ 
	c_{x} 
\end{align}
for some vectors $c_{x}\in \mathbb{R}^n$, 
where the parameters obey $\pi_t \lesssim (\lambda-1)^{-3/4}n^{1/4}\log n$ and $c_{0} \lesssim \frac{\log^4 n}{\sqrt{n(\lambda - 1)^3}}$.
Here, the last equation makes use of the fact that
\begin{align}
\label{eqn:tmp-cx}
\notag \|c_{x}\|_{2} & \lesssim\lt\|\frac{n\big|v_{t}+\sum_{k=1}^{t-1}\gamma_{t-1}^{k}\eta_{k}(v_{k})\big|^{5}\log^{4}n}{(\lambda-1)^{3}}\rt\|_{2}+\lt\|\frac{n|v_{t}|^{5}\log^{4}n}{(\lambda-1)^{3}}\rt\|_{2}\\
& \lesssim\frac{n\log^{4}n}{(\lambda-1)^{3}}\sqrt{\lt(\frac{\log n}{n}\rt)^{4}}\lesssim\frac{\log^{6}n}{n(\lambda-1)^{3}},
\end{align} 
where the property \eqref{enq:kpower} is invoked with $k = 10.$
Next, observe that 
\begin{align*}
& \Bigg\|\lt(v_t + \sum_{k=1}^{t-1} \gamma_{t-1}^k\eta_k(v_k)\rt)^3 - \lt(v_t\rt)^3\Bigg\|_2 \\
&=
\Bigg\|\sum_{k=1}^{t-1} \gamma_{t-1}^k\eta_k(v_k) \circ 
\left((v_t + \sum_{k=1}^{t-1} \gamma_{t-1}^k\eta_k(v_k))^2 + \lt(v_t\rt)^2 + (v_t + \sum_{k=1}^{t-1} \gamma_{t-1}^k\eta_k(v_k))\circ v_t\right) \Bigg\|_2  \\
&\lesssim  
\Big\|\max_{1\leq k\leq t-1} \eta_k(v_k) \circ \Big((v_t + \sum_{k=1}^{t-1} \gamma_{t-1}^k\eta_k(v_k))^2 + \lt(v_t\rt)^2\Big) \Big\|_2 \cdot \sqrt{t} \|\xi_{t-1}\|_2 \\
&\lesssim \max_{1\leq k\leq t-1}\left\|\eta_k(v_{k})\right\|_{\infty}\cdot\Big(\big\|(v_t + \sum_{k=1}^{t-1} \gamma_{t-1}^k\eta_k(v_k))^2\big\|_2 + \big\|\lt(v_t\rt)^2\big\|_2\Big)\cdot \sqrt{t} \|\xi_{t-1}\|_2\\
&\lesssim {\frac{t\log n}{n} \cdot \|\xi_{t-1}\|_2} \lesssim {\frac{\log^{4} n}{n^{1.5}(\lambda-1)^{2.5}}},
\end{align*}
where the last line can be obtained by invoking property~\eqref{eqn:v-inf-3424} and~\eqref{enq:kpower} with $k=4$. 
Here, we have used the facts that $\ltwo{\xi_t} \lesssim \sqrt{\frac{t^3\log n}{n}}$ (cf.~\eqref{eq:xi-stage1}), $t \lesssim \frac{\log n}{\lambda - 1}$ and 
\begin{align*}
  \Big\|v_t + \sum_{k=1}^{t-1} \gamma_{t-1}^k\eta_k(v_k)\Big\|_2 
  \leq 
  \ltwo{v_t} + \Big\|\sum_{k=1}^{t-1} \gamma_{t-1}^k\eta_k(v_k)\Big\|_2
  \leq
  1 + \sqrt{t}\ltwo{\gamma_{t-1}} \leq 1+\sqrt{\frac{t^4\log n}{n}}
  \lesssim 1.
\end{align*}
Putting these together, we arrive at 
\begin{align*}
 \lt\| \eta_{t}\lt(v_t + \sum_{k=1}^{t-1} \gamma_{t-1}^k\eta_k(v_k)\rt) - \eta_{t}\lt(v_t\rt) - (1-c_0)\sum_{k=1}^{t-1} \gamma_{t-1}^k\eta_k(v_k) \rt\|_2 
	&=   O\lt(\frac{\log^{6} n}{n(\lambda-1)^{3}}\rt).
\end{align*}

\item Finally, it is sufficient for us to consider $v^{\star\top}\sum_{k=1}^{t-1} \gamma_{t-1}^k\eta_k(v_k)$ which shall be controlled as follows:  
\begin{align*}
 \lt| \sum_{k=1}^{t-1}\gamma_{t-1}^{k}v^{\star\top}\eta_{k}\lt(v_k\rt) \rt|
&\le \lt| \sum_{k=1}^{t-1}\gamma_{t-1}^{k} \lt[v^{\star\top}\eta_{k}\lt(x_k\rt) + O(\|\xi_{k-1}\|_2)\rt] \rt| \\
&=  \lt| 
\sum_{k=1}^{t-1}\gamma_{t-1}^{k}\lt(\frac{\alpha_{k+1}}{\lambda} + O\lt(\sqrt{\frac{k^3\log n}{n}}\rt)\rt) \rt| \\ 
&\lesssim \sqrt{t}\|\gamma_{t-1}\|_2\cdot\frac{\sqrt{\lambda-1}}{n^{1/4}} 
\asymp \sqrt{t}\|\xi_{t-1}\|_2\cdot\frac{\sqrt{\lambda-1}}{n^{1/4}}
\lesssim \frac{\log^{2.5} n}{n^{3/4}(\lambda-1)^{1.5}}.
\end{align*}
\end{itemize}
Putting the above three inequalities together yields the desired bound \eqref{eqn:alphat-decomp}.


Built upon expression~\eqref{eqn:alphat-decomp}, we now proceed to the proof of claim~\eqref{eq:alpha_t}. 
To begin with, let us recall that $v_t \defn \alpha_t v^{\star} + \sum_{k=1}^{t-1}\beta_{t-1}^{k}\phi_k$.
If we define $g_{t-1} \defn v^{\star\top}\phi_{t-1}\sim\mathcal{N}\lt(0,\frac{1}{n}\rt)$ and $\widetilde{\beta}_{t-1} \defn [\beta_{t-1}^1,\ldots,\beta_{t-1}^{t-2}]$, some direct algebra thus leads to  
\begin{align*}
\lt|v^{\star\top}\lt(v_{t}-\alpha_{t}v^{\star}-\phi_{t-1}\rt)\rt| & =\lt\vert v^{\star\top}\Big[\sum_{k=1}^{t-2}\beta_{t-1}^{k}\phi_{k}-(1-\beta_{t-1}^{t-1})\phi_{t-1}\Big]\rt\vert\\
 & \le\Big|\sum_{k=1}^{t-2}\beta_{t-1}^{k}g_{k}\Big|+|1-\beta_{t-1}^{t-1}|\vert g_{t-1}\vert\\
 & =\Big|\sum_{k=1}^{t-2}\beta_{t-1}^{k}g_{k}\Big| + \frac{1-\big(\beta_{t-1}^{t-1}\big)^{2}}{1+|\beta_{t-1}^{t-1}|}\vert \, g_{t-1}\vert\\
 & \lesssim\|\widetilde{\beta}_{t-1}\|_{2}\sqrt{\sum_{k=1}^{t-2}(g_{k})^{2}}+\|\widetilde{\beta}_{t-1}\|_{2}^{2}\vert g_{t-1}\vert\lesssim\|\widetilde{\beta}_{t-1}\|_{2}\sqrt{\frac{t\log n}{n}}\\
 & \lesssim\sqrt{\frac{t\log n}{n}}\cdot\lt(\frac{t\sqrt{\lambda-1}}{n^{1/4}}+\frac{t\log^{4}n}{\sqrt{n(\lambda-1)^{3}}}\rt)\lesssim\frac{\log^{2}n}{n^{3/4}(\lambda-1)}+\frac{\log^{6}n}{n(\lambda-1)^{3}},
\end{align*}
where the last inequality comes from the bound~\eqref{eq:beta-tilde} in the proof of Lemma~\ref{lem:linear-ind} and the condition $t \lesssim \frac{\log n}{\lambda - 1}$.
By virtue of the above calculations, we can deduce that 
\begin{align*}
  v^{\star\top} v_t =  \alpha_t + g_{t-1} + O\lt(\frac{\log^{2} n}{n^{3/4}(\lambda-1)} + \frac{\log^{6} n}{n(\lambda-1)^3}\rt).
\end{align*}
In fact, a direct application of Lemma~\ref{lem:eta_property} further leads to the following claim:  
\begin{align}
\label{eqn:mahler3}
\notag  \lt| v^{\star\top}\big(\eta_t(v_t)-v_t\big) \rt|
&\lesssim  \lt| v^{\star\top}\lt[c_0 v_t + \pi_t^2(v_t\circ v_t\circ v_t) + O\lt(\frac{n\log^{4} n}{(\lambda-1)^{3}}|v_t|^5\rt) \rt] \rt| \\
&\lesssim \frac{\log^4 n}{n^{3/4}(\lambda-1)},
\end{align}
whose the proof of the last inequality is postponed to the end of this subsection. 

To summarize, taking the above results collectively and using the relation~\eqref{eqn:alphat-decomp}, we arrive at  
\begin{align}
\label{eq:stage-1-recursion}
\notag \alpha_{t+1} &= \lambda v^{\star\top}v_t+\lambda v^{\star\top}\big(\eta_{t}\lt(v_t\rt)-v_t \big) + O\lt(\frac{\log^{2.5} n}{n^{3/4}(\lambda-1)^{1.5}}\rt) \\
&=\lambda \alpha_t + \lambda g_{t-1} + O\lt(\frac{\log^4 n}{n^{3/4}(\lambda-1)^{1.5}}\rt). 
\end{align}
Therefore, invoking the above relation recursively leads to our desired decomposition: 
\begin{align*}
\alpha_{t+1} = \lambda^{t-k+1}\alpha_k + \sum_{i = 1}^{t-k+1} \lambda^ig_{t-i} + O\lt(\sum_{i = 1}^{t-k+1} \lambda^i\frac{\log^4 n}{n^{3/4}(\lambda-1)^{1.5}}\rt)
\end{align*}
for any  $1 \le k \le t$.

\paragraph{Proof of inequality~\eqref{eqn:mahler3}.}
In order to establish inequality~\eqref{eqn:mahler3}, let us first make note of the following simple  properties: with probability at least $1 - O(n^{-11})$,  
\begin{align*}
v^{\star\top}(v^{\star}\circ v^{\star}\circ v^{\star}) &= \frac{1}{n}; \\
v^{\star\top}\Big(v^{\star}\circ v^{\star}\circ \sum_{k=1}^{t-1}\beta_{t-1}^{k}\phi_{k}\Big) &= \frac{1}{n}v^{\star\top}\Big(\sum_{k=1}^{t-1}\beta_{t-1}^{k}\phi_{k}\Big) \lesssim \frac{1}{n}; \\
v^{\star\top}\Big(v^{\star}\circ \sum_{k=1}^{t-1}\beta_{t-1}^{k}\phi_{k}\circ \sum_{k=1}^{t-1}\beta_{t-1}^{k}\phi_{k}\Big) &= \frac{1}{n}\Big\|\sum_{k=1}^{t-1}\beta_{t-1}^{k}\phi_{k}\Big\|_2^2 \asymp \frac{1}{n}; \\
v^{\star\top}\Big(\sum_{k=1}^{t-1}\beta_{t-1}^{k}\phi_{k}\circ \sum_{k=1}^{t-1}\beta_{t-1}^{k}\phi_{k}\circ \sum_{k=1}^{t-1}\beta_{t-1}^{k}\phi_{k}\Big) &\lesssim \frac{t\log n}{n}v^{\star\top}\Big(\sum_{k=1}^{t-1}\beta_{t-1}^{k}\phi_{k}\Big) \lesssim \sqrt{\frac{t^3\log^3 n}{n^3}}.
\end{align*}
We remind the readers that $\vstar_{i} \sim \mathsf{Unif}(\pm \frac{1}{\sqrt{n}})$ and we have invoked Lemma~\ref{lem:concentration}.

Next, recall $v_t \defn \alpha_t v^{\star} + \sum_{k=1}^{t-1}\beta_{t-1}^{k}\phi_k$ to obtain 
\begin{align*}
\lt| v^{\star\top}(v_t\circ v_t\circ v_t) \rt| &= \lt| \alpha_t^3v^{\star\top}(v^{\star}\circ v^{\star}\circ v^{\star}) + 3\alpha_t^2v^{\star\top}\Big(v^{\star}\circ v^{\star}\circ \sum_{k=1}^{t-1}\beta_{t-1}^{k}\phi_{k}\Big)  \rt. \\
&\qquad \lt. + 3\alpha_tv^{\star\top}\Big(v^{\star}\circ \sum_{k=1}^{t-1}\beta_{t-1}^{k}\phi_{k}\circ \sum_{k=1}^{t-1}\beta_{t-1}^{k}\phi_{k}\Big)+ v^{\star\top}\Big(\sum_{k=1}^{t-1}\beta_{t-1}^{k}\phi_{k}\circ \sum_{k=1}^{t-1}\beta_{t-1}^{k}\phi_{k}\circ \sum_{k=1}^{t-1}\beta_{t-1}^{k}\phi_{k}\Big) \rt| \\
&\lesssim \frac{\alpha_t^3}{n} + \frac{\alpha_t^2}{n} + \frac{\alpha_t}{n} + \sqrt{\frac{t^3\log^3 n}{n^3}}\\ 
&\lesssim \frac{1}{n^{5/4}},
\end{align*}
where the last line holds as long as $\alpha_t \lesssim \sqrt{\lambda-1}\,n^{-1/4}$. 
Consequently, in order to derive~\eqref{eqn:mahler3}, it suffices to notice \eqref{enq:kpower}, $c_0\lesssim \frac{\log^4 n}{\sqrt{n(\lambda-1)^3}}$, $\pi_t\lesssim \frac{n^{1/4}}{(\lambda-1)^{3/4}}$ and
\begin{align*}
  \lt| v^{\star\top} v_t \rt| =  \lt| \alpha_t + v^{\star\top}\Big(\sum_{k=1}^{t-2}\beta_{t-1}^{k}\phi_k\Big)\rt|
  \lesssim \frac{\sqrt{\lambda-1}}{n^{1/4}}
  +
  \sqrt{\frac{t\log n}{n}}  \lesssim \frac{\sqrt{\lambda-1}}{n^{1/4}}.
\end{align*}


\subsection{Proof of Claim~\eqref{eq:alpha-global}}
\label{sec:proof-alpha-global}

For notational simplicity, we assume without loss of generality that $\alpha_t>0$ throughout this proof. 
Before delving into the proof of claim~\eqref{eq:alpha-global}, let us recall Lemma~\ref{lem:parameters} to obtain 
\begin{align*}
\frac{\pi_t}{\alpha_t\sqrt{n}} 
&= 1 + O\bigg(\frac{1}{\alpha_t^2}\lt(\|\xi_{t-1}\|_2 + \sqrt{\frac{t\log n}{n}}\rt)\wedge \frac{1}{\alpha_t}\lt(\|\xi_{t-1}\|_2 + \sqrt{\frac{t\log n}{n}}\rt)^{1/2}\bigg) \\
&= 1 + O\bigg(\frac{\log^2 n}{\alpha_t^2\sqrt{(\lambda-1)^3n}}\wedge \lt(\frac{\log^2 n}{\alpha_t^2\sqrt{(\lambda-1)^3n}}\rt)^{1/2}\bigg),
\end{align*}
where we have used  $\ltwo{\xi_t} \leq \sqrt{\frac{t^3\log n}{n}}$ (see \eqref{eq:xi-stage1}) and $t\lesssim \frac{\log n}{\lambda - 1}$. 
In turn, this implies 
\begin{align}
\label{eqn:pit-stage2}
\lt|\frac{\pi_t^2}{\alpha_t^2n} - 1\rt|\alpha_t^2 \lesssim \frac{\log^2 n}{\sqrt{(\lambda-1)^3n}}.
\end{align}

Now, let us move on to establish a recursive relation of $\alpha_{t}$. 
Recalling the definition \eqref{defi:eta} of $\eta_{t}$ and Theorem~\ref{thm:main-AMP}, one sees that 
\begin{align}
\label{eqn:alpha-t-merry}
\notag \alpha_{t+1} 
  &= \lambda v^{\star \top} \int{\eta}_{t}\lt(\alpha_t v^{\star} + \frac{1}{\sqrt{n}}x\rt)\varphi_n(\dx) + \Delta_{\alpha,t} \\
\notag  &=\lambda\gamma_t v^{\star\top}\int\tanh\lt(\pi_t\big(\alpha_t v^{\star}+\frac{1}{\sqrt{n}}x\big)\rt)\varphi_n(\dx)+\Delta_{\alpha,t} \\
&= \lambda\gamma_t\sqrt{n} \int \tanh\lt(\frac{\pi_t}{\sqrt{n}}\lt(\alpha_t + x\rt)\rt)\varphi(\dx) + \Delta_{\alpha,t},
\end{align} 
where the last equality holds by symmetry of $\varphi(\cdot)$, namely,  
\begin{align*}
\frac{1}{\sqrt{n}}\int\tanh\lt(\frac{\pi_t}{\sqrt{n}}(\alpha_t + x)\rt)\varphi(\dx) = -\frac{1}{\sqrt{n}}\int\tanh\lt(\frac{\pi_t}{\sqrt{n}}(-\alpha_t + x)\rt)\varphi(\dx).
\end{align*}
We note that  similar analysis as for relation~\eqref{eqn:ut-num} leads to $\int \tanh^2\lt(\frac{\pi_t}{\sqrt{n}}(\alpha_t + x)\rt)\varphi(\dx)\asymp\frac{\pi_t^2}{n}$. Combining this result with Lemma~\ref{lem:parameters} and \eqref{eq:xi-stage1}, we arrive at
\begin{align}
\label{eqn:gammat-ii}
  \gamma_t^{-2} = n\lt(1 + O\lt(\sqrt{\frac{t^3\log n}{n}}\rt)\rt)\int \tanh^2\lt(\frac{\pi_t}{\sqrt{n}}\lt(\alpha_t + x\rt)\rt)\varphi(\dx). 
\end{align}
Taking \eqref{eqn:alpha-t-merry} and \eqref{eqn:gammat-ii} together, we arrive at 
\begin{align}
\label{eqn:alphat-mid}
\alpha_{t+1} 
    &= 
    \Big(1 + O\Big(\sqrt{\frac{t^3\log n}{n}}\Big)\Big)\frac{\lambda\int \tanh\lt(\frac{\pi_t}{\sqrt{n}}\lt(\alpha_t + x\rt)\rt)\varphi(\dx)}{\lt[\int \tanh^2\lt(\frac{\pi_t}{\sqrt{n}}\lt(\alpha_t + x\rt)\rt)\varphi(\dx)\rt]^{1/2}} + \Delta_{\alpha,t}. 
%
\end{align}

To prove claim~\eqref{eq:alpha-global}, it then suffices to control $\int \tanh\lt(\frac{\pi_t}{\sqrt{n}}\lt(\alpha_t + x\rt)\rt)\varphi(\dx)$ and $\int \tanh^2\lt(\frac{\pi_t}{\sqrt{n}}\lt(\alpha_t + x\rt)\rt)\varphi(\dx)$.
Towards this goal, we find it helpful to first make several observations. 
Define two functions:  
\begin{align*}
f(z) &:= \frac{1}{z}\tanh(zy) - \tanh(y), \\
g(z) &:= \frac{1}{z^2}\tanh^2(zy) - \tanh^2(y). 
\end{align*}
The Taylor expansion of $\tanh(zy)$ gives
\begin{align*}
f^{\prime}(z) &= -\frac{1}{z^2}\big[\tanh(zy) - zy + zy\tanh^2(zy)\big] = -\frac{2}{3}zy^3 + O(z^3y^5), \\
g^{\prime}(z) &= -\frac{2\tanh(zy)}{z^3}\big[\tanh(zy) - zy + zy\tanh^2(zy)\big] = \frac{1}{3}zy^4 + O(z^3y^6),
\end{align*}
which leads the following relation by direct calculation
\begin{align*}
f(z)&=\int_{1}^{z} f^{\prime}(t)dt = -\frac{1}{3}(z^2-1)y^3 + O((z^4-1)y^5),\\
g(z)&=\int_{1}^{z} g^{\prime}(t)dt = \frac{1}{6}(z^2-1)y^4 + O((z^4-1)y^6).
\end{align*}
By taking $z=\frac{\pi_t}{\alpha_t\sqrt{n}}$, $y=\alpha_t(\alpha_t + x)$, we can see that 
\begin{align*}
\frac{\alpha_t\sqrt{n}}{\pi_t}\tanh\lt(\frac{\pi_t}{\sqrt{n}}\lt(\alpha_t + x\rt)\rt) - \tanh\lt(\alpha_t\lt(\alpha_t + x\rt)\rt)
&= -\frac{1}{3}\lt(\frac{\pi_t^2}{\alpha_t^2n} - 1\rt)\alpha_t^3(\alpha_t + x)^3 + O\lt(\frac{\pi_t^4}{\alpha_t^4n^2} - 1\rt)\alpha_t^5(\alpha_t + x)^5, \\
\frac{\alpha_t^2n}{\pi_t^2}\tanh^2\lt(\frac{\pi_t}{\sqrt{n}}\lt(\alpha_t + x\rt)\rt) - \tanh^2\lt(\alpha_t\lt(\alpha_t + x\rt)\rt)
&= \frac{1}{6}\lt(\frac{\pi_t^2}{\alpha_t^2n} - 1\rt)\alpha_t^4(\alpha_t + x)^4 + O\lt(\frac{\pi_t^4}{\alpha_t^4n^2} - 1\rt)\alpha_t^6(\alpha_t + x)^6.
\end{align*}
Hence, we can conclude that
\begin{align*}
&\lt|\int \frac{\alpha_t\sqrt{n}}{\pi_t}\tanh\lt(\frac{\pi_t}{\sqrt{n}}\lt(\alpha_t + x\rt)\rt) - \tanh\lt(\alpha_t\lt(\alpha_t + x\rt)\rt)\varphi(\dx)\rt| \\
&=
\lt|\int \frac{1}{3}\lt(\frac{\pi_t^2}{\alpha_t^2n} - 1\rt)\alpha_t^3(\alpha_t + x)^3 \varphi(\dx)\rt| + \int O\lt(\frac{\pi_t^4}{\alpha_t^4n^2} - 1\rt)\alpha_t^5(\alpha_t + x)^5 \varphi(\dx)\\
&\lesssim \lt|\frac{\pi_t^2}{\alpha_t^2n} - 1\rt|\alpha_t^4. 
\end{align*}
Similarly, we can show that  
\begin{align}
\lt|\int\frac{\alpha_t^2n}{\pi_t^2}\tanh^2\lt(\frac{\pi_t}{\sqrt{n}}\lt(\alpha_t + x\rt)\rt) - \tanh^2\lt(\alpha_t\lt(\alpha_t + x\rt)\rt)\varphi(\dx)\rt| 
&\lesssim \lt|\frac{\pi_t^2}{\alpha_t^2n} - 1\rt|\alpha_t^4. \label{eq:integral-pi}
\end{align}

Substituting these relations into \eqref{eqn:alphat-mid}, we arrive at  
\begin{align}
\label{eqn:miss}
  \notag &\alpha_{t+1} =\lambda\Big(1 + O\Big(\sqrt{\frac{t^3\log n}{n}}\Big)\Big)\frac{\int \tanh\lt(\alpha_t\lt(\alpha_t + x\rt)\rt)\varphi(\dx) + O\lt(\lt|\frac{\pi_t^2}{\alpha_t^2n} - 1\rt|\alpha_t^4\rt)}{\lt[\int \tanh^2\lt(\alpha_t\lt(\alpha_t + x\rt)\rt)\varphi(\dx) + O\lt(\lt|\frac{\pi_t^2}{\alpha_t^2n} - 1\rt|\alpha_t^4\rt)\rt]^{1/2}} + \Delta_{\alpha,t}\\
  \notag &=\lambda\Big(1 + O\Big(\sqrt{\frac{t^3\log n}{n}}\Big)\Big)\frac{\int \tanh\lt(\alpha_t\lt(\alpha_t + x\rt)\rt)\varphi(\dx)}{\lt[\int \tanh^2\lt(\alpha_t\lt(\alpha_t + x\rt)\rt)\varphi(\dx)\rt]^{1/2}}\cdot
  \frac{1+ O\lt(\lt|\frac{\pi_t^2}{\alpha_t^2n} - 1\rt|\alpha_t^2\rt)}{1+ O\lt(\lt|\frac{\pi_t^2}{\alpha_t^2 n}-1\rt|\alpha_t^2\rt)}+\Delta_{\alpha,t}\\
  &= \lambda \lt[\int \tanh^2\lt(\alpha_t^2 + \alpha_tx\rt)\varphi(\dx)\rt]^{-1/2}\int \tanh\lt(\alpha_t^2 + \alpha_tx\rt)\varphi(\dx) + O\lt(\lt|\frac{\pi_t^2}{\alpha_t^2n} - 1\rt|\alpha_t^3 + \sqrt{\frac{t^3\log n}{n}}\alpha_t + \Delta_{\alpha,t}\rt), 
\end{align}
where we make use of the fact that (see, \eqref{eqn:ut-num})
\begin{align*}
\int \tanh\lt(\alpha_t^2 + \alpha_tx\rt)\varphi(\dx) = \int \tanh^2\lt(\alpha_t^2 + \alpha_tx\rt)\varphi(\dx) \asymp \alpha_t^2. 
\end{align*}
In addition, \eqref{eqn:bayes-opt} ensures that
\begin{align*}
   \frac{\alpha_t}{\sqrt{\alpha^2_t + 1 }} 
   \leq \sqrt{\int \tanh^2(\alpha_t(\alpha_t + x)) \varphi(\dx)}
   =
   \lt[\int \tanh^2\lt(\alpha_t^2 + \alpha_tx\rt)\varphi_n(\dx)\rt]^{-1/2}\int \tanh\lt(\alpha_t^2 + \alpha_tx\rt)\varphi(\dx),
\end{align*}
which in turn gives 
\begin{align*}
 \alpha_{t+1} \geq \frac{\lambda \alpha_t}{\sqrt{\alpha^2_t + 1 }} + O\lt(\lt|\frac{\pi_t^2}{\alpha_t^2n} - 1\rt|\alpha_t^3 + \sqrt{\frac{t^3\log n}{n}}\alpha_t + |\Delta_{\alpha,t}| \rt). 
\end{align*}
To finish up, putting the above results together with \eqref{eqn:pit-stage2} 
leads to 
\begin{align}
\label{eqn:alpha-recursion}
  \alpha_{t+1} 
  &\ge \frac{\lambda\alpha_t}{\sqrt{\alpha^2_t + 1 }} 
  + o\big( (\lambda-1)\alpha_t \big) + O(\Delta_{\alpha,t})
\end{align}
where we again invoke the assumption that $\lambda - 1 \gtrsim n^{-1/9}\log n.$
This concludes the proof of claim~\eqref{eq:alpha-global}.

\subsection{Proof of Claim~\eqref{eqn:delta-alpha}}
\label{sec:proof-claim:eqn:delta-alpha}

Consider the regime where 
$$ |\alpha_t| < (\lambda-1)^{-3/4}n^{-1/4} \lesssim \sqrt{\lambda-1} \, n^{-0.1}.$$
First, invoke property~\eqref{eq:eta_derivative_small} in Lemma~\ref{lem:eta_property} to ensure that
\begin{align}
	\big|v^{\star\top}\eta_{t}(x_t) - v^{\star\top}\eta_{t}(v_t)\big| \lesssim |v^{\star\top}\xi_{t-1}| + (\lambda-1)n^{-0.2}(\log n) \|\xi_{t-1}\|_2.
	\label{eq:diff-etat-vstar-135}
\end{align}
To further bound \eqref{eq:diff-etat-vstar-135}, 
 note that  $\xi_{t-1}$ admits the following decomposition in terms of $\{\eta_{k}(x_k)\}$:
\begin{align*}
	\xi_{t-1} = \sum_{k=1}^{t-2} \gamma_{t-1}^k\eta_k(x_k),\qquad\text{with }\gamma_{t-1}=[\gamma_{t-1}^k]_{1\leq k \leq t-2}\in\real^{t-2} \text{ obeying }\|\gamma_{t-1}\|_2 \lesssim \|\xi_{t-1}\|_2; 
\end{align*}
the proof of this claim  can be found in Section~\ref{sec:claim-alpha-t} (see Lemma~\ref{lem:linear-ind} therein and its proof). 
In view of this relation, we can apply \eqref{eq:diff-etat-vstar-135} and the Cauchy-Schwarz inequality to reach 
\begin{align*}
\big|v^{\star\top}\eta_{t}(x_{t})-v^{\star\top}\eta_{t}(v_{t})\big| & \lesssim\bigg|\sum_{k=1}^{t-2}\gamma_{t-1}^{k}v^{\star\top}\eta_{k}(x_{k})\bigg|+(\lambda-1)n^{-0.2}(\log n)\|\xi_{t-1}\|_{2}\\
 & \lesssim\|\gamma_{t-1}\|_{2}\bigg(\sum_{k=1}^{t-2}\big(v^{\star\top}\eta_{k}(x_{k})\big)^{2}\bigg)^{1/2}+(\lambda-1)n^{-0.2}(\log n)\|\xi_{t-1}\|_{2}\\
 & \lesssim\sqrt{t}\|\xi_{t-1}\|_{2}\max_{\tau_0\leq s\leq t}|\alpha_{s}|+(\lambda-1)n^{-0.2}(\log n)\|\xi_{t-1}\|_{2}\lesssim\sqrt{\frac{t\log n}{n}}, 
\end{align*}
provided that $t \lesssim {\frac{\log n}{\lambda - 1}}$ and $\displaystyle\max_{\tau_0\leq s\leq t}|\alpha_{s}|\lesssim\sqrt{\lambda-1}n^{-0.1}$. Here the last inequality invokes 
$\|\xi_{t}\|_{2} \lesssim \sqrt{\frac{t^3\log n}{n}}$ (see \eqref{eq:xi-stage1}). 
Substitution into \eqref{eq:delta-alpha-general-579} yields
\begin{align}
|\Delta_{\alpha,t}| \lesssim \sqrt{\frac{t\log n}{n}}
\ll (\lambda - 1)|\alpha_t|,
\end{align}
given $|\alpha_t|\gtrsim \sqrt{\lambda-1}n^{-1/4}$ and $\lambda-1 \gtrsim n^{-1/9}\log n$.
It thus completes the proof of the relation~\eqref{eqn:delta-alpha}.

\subsection{Proof of Claim~\eqref{eq:condition-alpha-t-lower-bound}}
\label{sec:proof:eq:condition-alpha-t-lower-bound}

Throughout this section, we assume without loss of generality that $\alpha_{t} > 0.$
As computed in Section~\ref{sec:proof-alpha-global} for relation~\eqref{eq:alpha-global}, applying Lemma~\ref{lem:parameters} reveals that 
\begin{align*}
\lt|\frac{\pi_t^2}{\alpha_t^2n} - 1\rt|\alpha_t^2 &\lesssim \lt(\|\xi_{t-1}\|_2 + \sqrt{\frac{t\log n}{n}}\rt)\wedge \alpha_t \lt(\|\xi_{t-1}\|_2 + \sqrt{\frac{t\log n}{n}}\rt)^{1/2} \lesssim \alpha_t \lt(\frac{(\lambda - 1)^3}{\log n}\rt)^{1/4}, \\
\gamma_t^{-2} &= n\lt(1 + O\lt(\sqrt{\frac{(\lambda - 1)^3}{\log n}}\rt)\rt)\int \tanh^2\lt(\frac{\pi_t}{\sqrt{n}}\lt(\alpha_t + x\rt)\rt)\varphi(\dx), \\
|\Delta_{\alpha,t}| &\lesssim \|\xi_{t-1}\|_2 + \sqrt{\frac{t\log n}{n}} \lesssim \sqrt{\frac{(\lambda - 1)^3}{\log n}},
\end{align*}
given the inductive assumptions $\|\xi_{t-1}\|_2 \lesssim \sqrt{\frac{(\lambda - 1)^3}{\log n}}$ when $t \lesssim \frac{n(\lambda - 1)^5}{\log^2 n}$ and  $\alpha_t \gtrsim \sqrt{\lambda^2 - 1}$. 
Now in view of similar calculations for \eqref{eqn:miss}, we can deduce that 
\begin{align}
\label{eqn:alpha-tmp}
\alpha_{t+1} 
&= \lambda \lt[\int \tanh\lt(\alpha_t^2 + \alpha_tx\rt)\varphi_n(\dx)\rt]^{1/2} + O\lt(\lt|\frac{\pi_t^2}{\alpha_t^2n} - 1\rt|\alpha_t^3 + \frac{\lambda - 1}{\sqrt{\log n}}\alpha_t + \Delta_{\alpha,t}\rt) \notag\\
&\ge \frac{\lambda\alpha_t}{\sqrt{1+\alpha_t^2}} + o\lt((\lambda-1)^{3/4}\alpha_t^2 + (\lambda-1)\alpha_t\rt)
\end{align}
where we have made use of the fact that $\sqrt{\frac{(\lambda - 1)^{3}}{\log n}} \ll (\lambda-1)\alpha_t$.

We then demonstrate that  this relation~\eqref{eqn:alpha-tmp} together with a little algebra indicates that $\alpha_{t+1} \geq\frac{1}{2}\sqrt{\lambda^2 - 1}$. 
Specifically,  consider the following two cases separately. 
\begin{itemize}
\item First, consider the case where $\alpha_t \le \frac{2}{3}\sqrt{\lambda^2 - 1}$. Akin to inequality~\eqref{eqn:algebra-sim}, relation~\eqref{eqn:alpha-tmp} implies the existence of some constant $c > 0$ such that 
\begin{align*}
\alpha_{t+1} &\ge \big(1 + c(\lambda-1)\big)\alpha_t + o\big((\lambda-1)\alpha_t\big) \ge \alpha_t \ge \frac{1}{2}\sqrt{\lambda^2 - 1}.
\end{align*}

\item Otherwise, consider the case where $\alpha_t > \frac{2}{3}\sqrt{\lambda^2 - 1}$.  
Recognizing the fact that $\frac{\lambda\alpha_t}{\sqrt{1+\alpha_t^2}}$ is monotonically increasing in $\alpha_t$, we arrive at 
\begin{align*}
\alpha_{t+1} &\ge \big(1 + c(\lambda-1)\big)\frac{2}{3}\sqrt{\lambda^2 - 1} + o\big((\lambda-1)^{3/4}\big) \ge \frac{1}{2}\sqrt{\lambda^2 - 1}.
\end{align*}
\end{itemize}
Thus, this completes the proof of our desired bound \eqref{eq:induction-main}.

\subsection{Proof of Lemma~\ref{lem:linear-ind}}
\label{proof-linear-ind}

Throughout the proof, we work with the event that $\left\{\phi_k\right\}_{k=1}^{t-1}\in \mathcal{E}$ (defined in Lemma~\ref{lem:concentration}  with $\delta = O(n^{-11})$), which holds true with probability at least $1 - O(n^{-11}).$
On this event, one has $\|\phi_{t}\|_\infty \lesssim \sqrt{\frac{\log n}{n}}$ and 
\begin{align}
\|x_{t}\|_\infty 
&\leq \|\alpha_t \vstar\|_\infty + \Big\|\sum_{k=1}^{t-1}\beta_{t-1}^k \phi_k\Big\|_\infty + \ltwo{\xi_{t-1}}
	\lesssim \frac{|\alpha_t|}{\sqrt{n}}+ \sqrt{\frac{t\log n}{n}} + \sqrt{\frac{t^3\log n}{n}}  \lesssim \sqrt{\frac{\log^4 n}{n(\lambda - 1)^3}}
	\label{eq:xt-inf-norm-bound}
\end{align}
for any $t\leq \tau_0$, 
where we remind the readers that (see \eqref{eq:xi-stage1}, the property $\|\beta_{t-1}\|_2=1$, the definition of $\mathcal{E}_1$, and the definition \eqref{eq:defn-tau0-threshold} of $\tau_0$)
\begin{align}
|\alpha_t| \lesssim \sqrt{\lambda - 1}\, n^{-1/4}, \qquad 
    \Big\|\sum_{k=1}^{t-1}\beta_{t-1}^k \phi_k\Big\|_\infty \lesssim \sqrt{\frac{t\log n}{n}} 
    \qquad \text{and }~
    \ltwo{\xi_{t-1}} \lesssim \sqrt{\frac{t^3\log n}{n}}
	\label{eq:three-properties-lem}
\end{align}
as long as $t \lesssim \frac{\log n}{\lambda  - 1}$ (see \eqref{eqn:tau0}).

To show that $\{\eta_1(x_1),\ldots,\eta_{t}(x_{t})\}$ forms a near-orthogonal basis, 
one strategy is to show that  $\eta_i(x_i) \approx \phi_{i-1}$ for each $1\leq i \leq t \leq \tau_0$,  which in turn implies that
\begin{align}
\lt\|\sum_{i=1}^{t} w_i \eta_i(x_i)\rt\|_2 
\approx & \Big\|\sum_{i=1}^{t} w_i \phi_{i-1} \Big\|_2 
\asymp \ltwo{w},
\end{align}
given that $\phi_{k}\sim \mathcal{N}(0,\frac{1}{n} I_n)$ are independent Gaussian vectors; 
here, we introduce $\phi_0 := x_1 \sim \mathcal{N}(0, \frac{1}{n} I_n)$ for notational convenience.
Guided by this intuition, we first use the triangle inequality to derive that 
\begin{align}
\label{eqn:linear-decomp}
\notag 
&\Bigg|\Big\|\sum_{i=1}^{t} w_i \eta_i(x_i)\Big\|_2 - \Big\|\sum_{i=1}^{t} w_i \beta_{i-1}^{i-1}\phi_{i-1} \Big\|_2\Bigg|\\
&\qquad \leq 
\Big\|\sum_{i=1}^{t} w_i \Big[\eta_i(x_i) - \eta_i(\beta_{i-1}^{i-1}\phi_{i-1})\Big]\Big\|_2 
+ \Big\|\sum_{i=1}^{t} w_i \Big[\eta_i(\beta_{i-1}^{i-1}\phi_{i-1})-\beta_{i-1}^{i-1}\phi_{i-1}\Big]\Big\|_2
\end{align}
for any vector $w \in \mathbb{R}^{t}$. 
In order to bound these terms, we proceed with the following three steps. 

 \begin{itemize}
 \item 
Let us first consider the difference between $\eta_t(x_t)$ and $x_t$. 
Invoking property~\eqref{eq:eta_small} in Lemma~\ref{lem:eta_property} allows us to express
\[
\eta_{t}(x_{t})=(1-c_{0})\bigg(x_{t}-\frac{\pi_{t}^{2}}{3}x_{t}\circ x_{t}\circ x_{t}+c_{x_{t}}\bigg),
\]
where $c_0 \lesssim \frac{\log^2 n}{\sqrt{n(\lambda-1)^{3}}}$, $\pi_{t} \lesssim \frac{n^{1/4}\log n}{(\lambda - 1)^{3/4}}$, 
		 and $c_{x_{t}}$ is a vector obeying (cf.~\eqref{eq:xt-inf-norm-bound})
		 \begin{equation}
			 \|c_{x_{t}}\|_{\infty}\lesssim\frac{n\|x_{t}\|_{\infty}^{5}\log^{4}n}{(\lambda-1)^{3}}\lesssim\frac{\log^{14}n}{n^{3/2}(\lambda-1)^{10.5}}.	 
			 \label{eq:UB-c-xt-inf}
		 \end{equation}
Then one has
\begin{equation}
	\left\Vert \eta_t(x_t)-x_t\right\Vert_2 \lesssim c_0 \Vert x_t\Vert_2 + \frac{1}{3}\pi_t2 \Vert x_t\circ x_t\circ x_t\Vert_2 + \sqrt{n}\,\|c_{x_t}\|_{\infty} .
	\label{eqn:bound-eta-x-123}
\end{equation}
We then claim that 
\begin{align}
\label{eqn:bound-eta-x}
\left\Vert \eta_t(x_t)-x_t\right\Vert_2
 &\lesssim \frac{{\log^3 n}}{\sqrt{n(\lambda-1)^3}};  
\end{align}
to streamline the presentation,  the proof of this claim is deferred to the end of this section. 
Applying the same argument once gain also leads to 
\begin{align}
\label{eqn:bound-phi}
\left\Vert\eta_t(\beta_{t-1}^{t-1}\phi_{t-1})-\beta_{t-1}^{t-1}\phi_{t-1}\right\Vert_2\lesssim \frac{\log^3 n}{\sqrt{n(\lambda-1)^3}}.  
\end{align}
%

 \item 
Combining relations~\eqref{eqn:bound-eta-x} and \eqref{eqn:bound-phi} and invoking the triangle inequality, we obtain 
\begin{align}
\label{eqn:eta-xi-allegro}
\lt\|\eta_{t}\lt(x_t\rt) - \eta_{t}\lt(\beta_{t-1}^{t-1}\phi_{t-1}\rt)\rt\|_2 &\leq \left\Vert \eta_t(x_t)-x_t\right\Vert_2 + \left\Vert x_t-\beta_{t-1}^{t-1}\phi_{t-1}\right\Vert_2 + \left\Vert\eta_t(\beta_{t-1}^{t-1}\phi_{t-1})-\beta_{t-1}^{t-1}\phi_{t-1}\right\Vert_2\notag\\
& =\lt\|x_t - \beta_{t-1}^{t-1}\phi_{t-1}\rt\|_2 + O\lt(\frac{\log^3 n}{\sqrt{n(\lambda-1)^3}}\rt)\notag \\
& = \Big\Vert\alpha_t v^{\star} + \sum_{k=1}^{t-2}\beta_{t-1}^{k}\phi_k +\xi_t\Big\Vert_2 + O\lt(\frac{\log^3 n}{\sqrt{n(\lambda-1)^3}}\rt)\notag \\
\notag & \leq \left\Vert\alpha_t v^{\star}\right\Vert_2 + \Big\Vert\sum_{k=1}^{t-2}\beta_{t-1}^{k}\phi_k\Big\Vert_2 + \ltwo{\xi_t} + O\lt(\frac{\log^3 n}{\sqrt{n(\lambda-1)^3}}\rt)\\
& \leq \Big\Vert\sum_{k=1}^{t-2}\beta_{t-1}^{k}\phi_k\Big\Vert_2 
+ |\alpha_t| + O\left(\sqrt{\frac{t^3\log n}{n}}\right) + O\lt(\frac{\log^3 n}{\sqrt{n(\lambda-1)^3}}\rt),
\end{align}
		 where the last line makes use of \eqref{eq:three-properties-lem}. 
		 Let us denote $\widetilde{\beta}_{t-1} \defn (\beta_{t-1}^1,\ldots,\beta_{t-1}^{t-2}) \in \real^{t-2}$, 
		 obtained by removing the last entry of $\beta_{t-1}$.
		 According to Lemma~\ref{lem:concentration}, we note that with probability $1 - O(n^{-11})$,  
\begin{align*}
\bigg\vert\Big\Vert\sum_{k=1}^{t-2}\beta_{t-1}^{k}\phi_k\Big\Vert_2-\ltwo{\widetilde{\beta}_{t-1}}\bigg\vert\leq\ltwo{\widetilde{\beta}_{t-1}}\cdot \sup_{a=[a_k]_{1\leq k< t-2}\in\mathcal{S}^{t-3}}\bigg\vert\Big\Vert\sum_{k=1}^{t-2}a_k\phi_k\Big\Vert_2-1\bigg\vert\lesssim \sqrt{\frac{t\log n}{n}}\ltwo{\widetilde{\beta}_{t-1}}, 
\end{align*}
which in turn implies that
\begin{align*}
    \bigg(1 - O\Big(\sqrt{\frac{t\log n}{n}}\Big)\bigg) \ltwo{\widetilde{\beta}_{t-1}} \leq \Big\Vert\sum_{k=1}^{t-2}\beta_{t-1}^{k}\phi_k\Big\Vert_2
    \leq \bigg(1 + O\Big(\sqrt{\frac{t\log n}{n}}\Big)\bigg) \ltwo{\widetilde{\beta}_{t-1}}.
\end{align*}
As a consequence, we can further control the right-hand side of \eqref{eqn:eta-xi-allegro} by 
\begin{align}
\label{eqn:bound-eta-xt-phi}
\lt\|\eta_{t}\lt(x_t\rt) - \eta_{t}\lt(\beta_{t-1}^{t-1}\phi_{t-1}\rt)\rt\|_2 
%
&\leq \bigg(1 + O\Big(\sqrt{\frac{t\log n}{n}}\Big)\bigg)\|\widetilde{\beta}_{t-1}\|_2 + O\lt(|\alpha_t| + \frac{\log^3 n}{\sqrt{n(\lambda-1)^3}}\rt).
\end{align}

 \item 
Our next step is concerned with bounding the term $\|\widetilde{\beta}_{t-1}\|_2.$
First, recall that $\beta_{t}$ corresponds to the linear coefficients of $\eta_{t}(x_t)$ when projected to the linear space $U_t$ (see \eqref{eq:eta-t-expand-zt} and \eqref{eq:defn-Ut-minus-1}). We can thus write 
\begin{align*}
	\|\widetilde{\beta}_{t}\|_2 =& \lt\|U_{t-1}^{\top}\eta_{t}\lt(x_t\rt)\rt\|_2 \\
	\le& \lt\|U_{t-1}^{\top}\lt[\eta_{t}\lt(x_t\rt) - \eta_{t}\lt(\beta_{t-1}^{t-1}\phi_{t-1}\rt)\rt]\rt\|_2 + \lt\|U_{t-1}^{\top}\lt[\eta_{t}\lt(\beta_{t-1}^{t-1}\phi_{t-1}\rt)-\beta_{t-1}^{t-1}\phi_{t-1}\rt]\rt\|_2 + \left\Vert U_{t-1}^{\top}\big(\beta_{t-1}^{t-1}\phi_{t-1}\big)\right\Vert_2\\
	\le& \bigg(1 + O\Big(\sqrt{\frac{t\log n}{n}}\Big)\bigg)\|\widetilde{\beta}_{t-1}\|_2 + O\lt(|\alpha_t| + \frac{\log^3 n}{\sqrt{n(\lambda-1)^3}}\rt) + O\lt(\frac{\log^3 n}{\sqrt{n(\lambda-1)^3}}\rt) + \left\Vert U_{t-1}^{\top}\big(\beta_{t-1}^{t-1}\phi_{t-1}\big)\right\Vert_2\\
=& \bigg(1 + O\Big(\sqrt{\frac{t\log n}{n}}\Big)\bigg)\|\widetilde{\beta}_{t-1}\|_2 + O\lt(|\alpha_t| + \frac{\log^3 n}{\sqrt{n(\lambda-1)^3}}\rt) + O\lt(\sqrt{\frac{t\log n}{n}}\rt)\\
=& \bigg(1 + O\Big(\sqrt{\frac{t\log n}{n}}\Big)\bigg)\|\widetilde{\beta}_{t-1}\|_2 + O\lt(|\alpha_t| + \frac{\log^3 n}{\sqrt{n(\lambda-1)^3}}\rt). 
\end{align*}
Here, the third line follows from~\eqref{eqn:bound-phi} and~\eqref{eqn:bound-eta-xt-phi};
the penultimate line holds due to the independence between $\phi_{t-1}$ and $U_{t-1}$ (see the properties below display~\eqref{eqn:nutcracker}) and hence $U_{t-1}^{\top}\phi_{t-1} \sim \mathcal{N}(0, \frac{1}{n}I_t)$; and the last line holds as long as $t\lesssim\frac{\log n}{\lambda-1}$. 
Recognizing that $\ltwo{\widetilde{\beta}_1}=\sqrt{1-\ltwo{\beta_{1}}^2} = 0$, we can apply the above relation recursively to yield
\begin{align*}
\ltwo{\widetilde{\beta}_{t}} &\leq \sum_{\tau=1}^{t} \bigg(1 + O\Big(\sqrt{\frac{t\log n}{n}}\Big)\bigg)^{t-1-\tau}\cdot O\lt(|\alpha_\tau| + \frac{\log^3 n}{\sqrt{n(\lambda-1)^3}}\rt)\notag\\
	&\lesssim \sum_{\tau=1}^{t} \lt( |\alpha_\tau| + \frac{\log^3 n}{\sqrt{n(\lambda-1)^3}} \rt) 
\lesssim \frac{t \sqrt{\lambda - 1}}{n^{1/4}} + \frac{t\log^3 n}{\sqrt{n(\lambda-1)^3}},
\end{align*}
where the penultimate inequality follows from the fact $\sqrt{\frac{t\log n}{n}}\lesssim\frac{1}{t}$ as
$t\lesssim \frac{\log n}{\lambda-1}$, and the last inequality uses $|\alpha_t| \lesssim \sqrt{\lambda - 1}\,n^{-1/4}$
for $t \leq \tau_0$ (see the definition of $\tau_{0}$ in \eqref{eq:defn-tau0-threshold}). 
Under the assumption that $\lambda - 1 \gtrsim n^{-1/9}\log n$, we can further obtain 
\begin{align}
\label{eq:beta-tilde}
\ltwo{\widetilde{\beta}_{t}} &\lesssim \frac{t \sqrt{\lambda - 1}}{n^{1/4}}.
\end{align}
Plugging this relation into \eqref{eqn:bound-eta-xt-phi} and using the condition $|\alpha_t| \lesssim \sqrt{\lambda - 1}\,n^{-1/4}$
($\forall t \leq \tau_0$) give 
\begin{align}
\lt\|\eta_{t}\lt(x_t\rt) - \eta_{t}\lt(\beta_{t-1}^{t-1}\phi_{t-1}\rt)\rt\|_2 
\lesssim \frac{t \sqrt{\lambda - 1}}{n^{1/4}}.
\end{align}
It is worth noting that the inequality~\eqref{eq:beta-tilde} also implies 
\begin{align*}
|\beta_t^t| = \sqrt{\|{\beta}_{t}\|_2^2 - \|\widetilde{\beta}_{t}\|_2^2}= \sqrt{1 - \|\widetilde{\beta}_{t}\|_2^2}= 1 - O\lt(\frac{t^{2}(\lambda - 1)}{n^{1/2}}\rt).
\end{align*} 

 \end{itemize}

\noindent  
To finish up, putting the above bounds together with expression~\eqref{eqn:linear-decomp}, we conclude that
\begin{align*}
&\Bigg|\Big\|\sum_{i=1}^{t} w_i \eta_i(x_i)\Big\|_2  - \Big\|\sum_{i=1}^{t} w_i \phi_{i-1} \Big\|_2 \Bigg|\\
 &\leq
(1-\beta_{t-1}^{t-1})\Big\|\sum_{i=1}^{t} w_i \phi_{i-1} \Big\|_2 + 
\Big\|\sum_{i=1}^{t} w_i \Big[\eta_i(x_i) - \eta_i(\beta_{i-1}^{i-1}\phi_{i-1})\Big]\Big\|_2 
+ \Big\|\sum_{i=1}^{t} w_i \Big[\eta_i(\beta_{i-1}^{i-1}\phi_{i-1})-\beta_{i-1}^{i-1}\phi_{i-1}\Big]\Big\|_2 \\
 &= 
O\lt(\frac{t^{2}(\lambda - 1)}{n^{1/2}}\rt)\Big\|\sum_{i=1}^{t} w_i \phi_{i-1} \Big\|_2  
+ O\lt(\frac{t \sqrt{\lambda - 1}}{n^{1/4}}\rt)\cdot\sqrt{t}\Vert w\Vert_2 + O\lt(\frac{\log^3 n}{\sqrt{n(\lambda-1)^3}}\rt)\cdot\sqrt{t}\Vert w\Vert_2 \\
 &= O\lt(\frac{t^{3/2} \sqrt{\lambda - 1}}{n^{1/4}} + \sqrt{\frac{t\log^6 n}{n(\lambda-1)^3}}\rt) \|w\|_2 
= o(\|w\|_2),
\end{align*}
where the last relation results from the facts $t < \tau_0 \lesssim \log n/(\lambda - 1)$ and the assumption $\lambda - 1\geq n^{-1/9}$.
Finally, observing $\|\sum_{i=1}^{t} w_i \phi_{i-1}\|_2 = (1 + O(\sqrt{\frac{t\log n}{n}}))\ltwo{w}$, we reach
\begin{align*}
	\Big\|\sum_{i=1}^{t} w_i \eta_i(x_i)\Big\|_2  = \big( 1+o(1) \big) \|w\|_2,
\end{align*}
which completes the proof of our desired bound.

\paragraph{Proof of inequality~\eqref{eqn:bound-eta-x}.} 
For any fixed integer $k \ge 2$ that does not scale with $n$, 
we can write 
\begin{align*}
\sum_{i=1}^n |x_{t, i}|^k \lesssim 
\sum_{i=1}^n |\alpha_t \vstar_{i}|^k + \sum_{i=1}^n |u_{t,i}|^k + \sum_{i=1}^n |\xi_{t-1, i}|^k, 
\qquad \text{with }~ u_t \defn \sum_{k=1}^{t-1}\beta_{t-1}^k \phi_k. 
\end{align*}
Let us bound each term separately. Firstly, recalling that
 $\ltwo{\xi_{t-1}} \lesssim  \sqrt{\frac{t^3\log n}{n}}$ (cf.~\eqref{eq:xi-stage1}) gives 
\begin{align}
\sum_{i=1}^n |\xi_{t-1, i}|^k \le \|\xi_{t-1}\|_2^k \lesssim \lt(\frac{t^3 \log n}{n}\rt)^{k/2}.
\end{align} 
Secondly, on the event $\{\phi_k\}_{k=1}^{t-1}\in\mathcal{E}$ (see Lemma~\ref{lem:concentration}), we see that $\ltwo{u_t} \lesssim 1$ and $\|u_{t}\|_{\infty} \leq \sqrt{(t\log n)/n}$. This in turn gives 
\begin{align}
\label{eqn:messi}
\sum_{i=1}^n |u_{t, i}|^k 
= \sum_{i \le t} |u_{t, (i)}|^k + \sum_{i > t} |u_{t, (i)}|^k 
\notag 
&\stackrel{(*)}{\lesssim} t\lt(\frac{t \log n}{n}\rt)^{k/2} + \lt(\frac{\log n}{n}\rt)^{k/2-1}\sum_{i > t} |u_{t, (i)}|^2 \\
&\lesssim t\lt(\frac{t \log n}{n}\rt)^{k/2} + \lt(\frac{\log n}{n}\rt)^{k/2-1} \notag\\
&\lesssim \lt(1+\frac{t^{k/2+1}\log n}{n}\rt)\lt(\frac{\log n}{n}\rt)^{k/2-1},
\end{align}
with $x_{(i)}$ denoting the $i$-th largest entry of $x$ (in magnitude). 
Here, to see why inequality $(*)$ holds, we recall that on the event $\{\phi_k\}_{k=1}^{t-1}\in\mathcal{E}$ (see Lemma~\ref{lem:concentration}), one has 
\begin{align*}
    \sup_{a\in \mathcal{S}^{t-2}} \sum_{i=1}^t\left|\sum_{k=1}^{t-1} a_k \phi_k\right|_{(i)}^2
    \lesssim \frac{t\log n}{n},
\end{align*}
which also guarantees that $|u_{t, (t+1)}| \lesssim \sqrt{\log n/n}.$
Putting the above pieces together, we obtain 
\begin{align}
\label{eqn:ny}
    \sum_{i=1}^n |x_{t, i}|^k &\lesssim 
    n\left(\frac{1}{n(\lambda-1)}\right)^{3k/4} + \lt(1+\frac{t^{k/2+1}\log n}{n}\rt)\lt(\frac{\log n}{n}\rt)^{k/2-1}+ \lt(\frac{t^3 \log n}{n}\rt)^{k/2}  
    &\lesssim
    \left(\frac{\log n}{n}\right)^2, 
\end{align}
where the last inequality is valid if we take $k = 6$, $t \lesssim \log n/(\lambda-1)$ and assume $\lambda - 1\gtrsim n^{-1/9}\log n.$
It therefore leads to 
\begin{align*}
    \pi_t^2 \Vert x_t\circ x_t\circ x_t\Vert_2
    = 
    \pi_t^2 \Big(\sum_{i=1}^n |x_{t,i}|^6\Big)^{1/2}
    \lesssim  
    \frac{n^{1/2}\log^2 n}{(\lambda - 1)^{3/2}} \cdot \frac{\log n}{n}
    = \frac{\log^3 n}{\sqrt{n (\lambda - 1)^3}},
\end{align*}
where we have used the bound on $\pi_t$ in Lemma~\ref{lem:eta_property} (the 3rd case). 
This together with \eqref{eq:UB-c-xt-inf}, \eqref{eqn:bound-eta-x-123} and the fact $c_0\lesssim \frac{\log^2n}{\sqrt{n(\lambda-1)^3}}$ concludes the proof of inequality~\eqref{eqn:bound-eta-x}.

\subsection{Proof of inequality~\eqref{eqn:derivative}}
\label{sec:proof-eqn:derivative}
Before proceeding, let us make several observations about $\frac{\tau_t}{h(\tau_t)}$. 
As discussed around \citet[display~(254)]{li2022non}, the sequence $\tau_{t}$ with $\tau_{t+1} = \lambda^{2}h(\tau_t)$ is monotonically increasing, which implies that $\frac{\tau_t}{h(\tau_t)} \leq \lambda^2$.
In addition, the optimality of the Bayes estimator (cf.~\eqref{eqn:bayes-opt}) implies that $\frac{\tau_t}{h(\tau_t)} \leq \tau_t+1$. Combining these two observations, we obtain  
\begin{align*}
	\frac{\tau_t}{h(\tau_t)} \leq (\tau_t+1) \wedge \lambda^2.
\end{align*}

In view of the inductive assumption, we have $\alpha_t^2 = (1+o(1)) \tau_{t}$ for $t \gtrsim \varsigma$. Hence, for every $\tau$ obeying $\min\{\tau_t, \alpha_t^2\} \leq \tau \leq \max\{\tau_t, \alpha_t^2\}$, it holds that $\tau = (1+o(1)) \tau_{t}$ with $\tau_{t} \gtrsim \lambda^{2} - 1$.
Define $\mathcal{T}_2$ as in display (263) of~\citet{li2022non} such that 
\begin{align}
	\mathcal{T}_2(s, \tau) \defn s^2 h'(\tau) = 
	s^2 \int \lt(1 + \frac{x}{2\sqrt{\tau}}\rt)\lt(1 - \tanh^2\lt(\tau + \sqrt{\tau} x\rt)\rt)\varphi(\dx).
\end{align}
Armed with this notation, we can bound the target quantity as 
\begin{align*}
\frac{\tau_t}{h(\tau_t)} h^{\prime}(\tau) \le \mathcal{T}_2(\sqrt{\tau_t+1} \wedge \lambda, \tau).
\end{align*}
Therefore, it suffices to upper bound the right-hand side of  the above inequality by $1 - c(\lambda - 1).$

Towards this end, direct calculations yield 
\begin{align}
\label{eqn:T2-lip}
	\frac{\mathcal{T}_2(\sqrt{\tau_t+1} \wedge \lambda, \tau) 
	- \mathcal{T}_2(\sqrt{\tau+1} \wedge \lambda, \tau)}{\mathcal{T}_2(\sqrt{\tau+1} \wedge \lambda, \tau)}
	= 
	(\sqrt{\tau+1} \wedge \lambda)^2 - (\sqrt{\tau_t+1} \wedge \lambda)^2 
	= o(\lambda-1).
\end{align}
Moreover, it has been proved numerically (see Figure 1 in~\citet{li2022non}) that
\begin{align*}
	\mathcal{T}_2(\lambda, \tau) \leq 1 - (\lambda - 1),
	\quad \text{ for } \lambda \in (0,1.2] \text{ and } \tau > \sqrt{\lambda^2-1}.
\end{align*}
Recognizing that $\tau = (1+o(1)) \tau_{t} \gtrsim \lambda^{2} - 1$, we can deduce from the relation above that 
\begin{align}
\label{eqn:T2-range}
	\mathcal{T}_2(\sqrt{\tau+1} \wedge \lambda, \tau) \leq 1 - \big((\sqrt{\tau+1}\wedge\lambda) - 1\big)
	=
	1 - c_1(\lambda - 1)
\end{align}
for some universal constant $c_1 > 0$.
Finally, putting  relations \eqref{eqn:T2-lip} and \eqref{eqn:T2-range} together,  
we arrive at 
\begin{align*}
\frac{\tau_t}{h(\tau_t)} h^{\prime}(\tau) \le \mathcal{T}_2(\sqrt{\tau_t+1} \wedge \lambda, \tau) = (1+o(\lambda-1))\mathcal{T}_2(\sqrt{\tau+1} \wedge \lambda, \tau) \le 1 - c(\lambda-1) 
\end{align*}
for some universal constant $c > 0$. 
We have thus finished the proof of relation~\eqref{eqn:derivative}. 

\section{Proof of expression~\eqref{eqn:opt} and Corollary~\ref{cor:opt}}
\label{sec:pf-cor}


To begin with, by definition \eqref{eq:defn-AMP-based-estimator} of $u_{t}$, one has
$$\ltwo{u_t} = \lt\| \frac{1}{\lambda\sqrt{n(\alpha_t^2+1)}} \tanh(\pi_t x_t) \rt\|_2 \leq \frac{1}{\lambda},$$ 
where we have used the fact that $|\tanh(\pi_t x_t)| < 1$.
Therefore the quantity of interest $\|\vstar v^{\star \top} - u_t u_t^\top\|_{\mathsf{F}}^2$ is uniformly upper bounded by a constant $1 + \frac{1}{\lambda^4}.$
In addition, we find it helpful to observe that 
\begin{align}
\label{eqn:rachmaninoff}
\notag \|\vstar v^{\star \top} - u_t u_t^\top\|_{\mathsf{F}}^2
	& = \|\vstar\|_2^4 - 2 (v^{\star \top} u_t)^2 + \|u_t\|_2^4\\
\notag &= 1 - 2 \frac{1}{n\lambda^2(\alpha_t^2+1)} \big(v^{\star \top} \tanh(\pi_t x_t) \big)^2+ \frac{1}{n^2\lambda^4(\alpha_t^2+1)^2} \ltwo{\tanh(\pi_t x_t)}^4 \\
	&= 1 -  \frac{1}{\lambda^4}\Bigg(\frac{2\alpha^2_{t+1}}{(\alpha_t^2+1)n\gamma_t^2} + \frac{1}{(\alpha_t^2+1)^2n^2\gamma_t^4}\Bigg).
	%
	%
\end{align}

To validate expression~\eqref{eqn:opt}, it is sufficient to notice that $\gamma^{-2}_{t} - n\alpha_t^2(\alpha_t^2+1) = o(1)$ with probability at least $1 - O(n^{-10})$ (according to Lemma~\ref{lem:parameters}), which in turn leads to 
\begin{align*}
	\lim_{t\to \infty} \lim_{n\to \infty} \Exs\big[\|\vstar v^{\star \top} - u_t u_t^\top\|_{\mathsf{F}}^2\big]
	=   1 - \frac{\alpha^{\star 4}}{\lambda^4}.
\end{align*}

To prove Corollary~\ref{cor:opt}, we again invoke Lemma~\ref{lem:parameters} to demonstrate that  
\begin{align}
\label{eqn:gamma_t_final}
\notag \frac{1}{n}\gamma_t^{-2} 
& = \frac{1}{n} \pi_t^2\bigg(\alpha_t^2+1+O\bigg(\frac{\pi_t^2}{n}+\Vert\xi_{t-1}\Vert_2+\sqrt{\frac{t\log n}{n}}\bigg)\bigg) \\
\notag &= \bigg(\alpha^2_t + O\bigg(\bigg(\|\xi_{t-1}\|_2 + \sqrt{\frac{t\log n}{n}}\bigg)^{1/2}\bigg)\bigg) \bigg(\alpha_t^2+1+O\bigg(\frac{\pi_t^2}{n}+\Vert\xi_{t-1}\Vert_2+\sqrt{\frac{t\log n}{n}}\bigg)\bigg)\\
\notag &= (\alpha^2_t + O(\delta^{1/2})) (\alpha_t^2+1+O(\delta))\\
&=\alpha^2_t (\alpha_t^2+1) + O(\delta^{1/2})
\end{align}
with $\delta \defn \sqrt{\frac{t\log n}{n(\lambda-1)^2}} + \sqrt{\frac{\log^4 n}{n(\lambda-1)^3}}$, 
where we plug in the bound on $\|\xi_t\|_2$ as in expression~\eqref{eqn:xi-norm-fin}. 
Substituting the expression~\eqref{eqn:gamma_t_final} into \eqref{eqn:rachmaninoff} yields 
\begin{align}
\label{eqn:ut-error}
	 \|\vstar v^{\star \top} - u_t u_t^\top\|_{\mathsf{F}}^2
	 =  1 -  \frac{\alpha^2_{t}}{\lambda^4}\Big(2\alpha^2_{t+1} - \alpha^2_{t} +O(\delta^{1/2})\Big).
\end{align}
After an order of $\frac{\log n}{\lambda -1}$ iterations,  property~\eqref{eqn:alpha-stage-3} ensures that 
$\alpha^2_{t} - \alpha^{\star 2} = O(\sqrt{\frac{\log^4 n}{n(\lambda - 1)^6}})$.
Putting everything together, we arrive at 
\begin{align*}
	 \|\vstar v^{\star \top} - u_t u_t^\top\|_{\mathsf{F}}^2
	 =  1 -  \frac{\alpha^{\star 4}}{\lambda^4}
	 +  O\bigg(\sqrt{\frac{\log^4 n}{n(\lambda - 1)^6}}\bigg),
\end{align*}
which holds true with probability at least $1 - O(n^{-10}).$

\bibliographystyle{apalike}
\bibliography{reference-amp}

\end{document}